\documentclass[12pt]{article}

\usepackage[a4paper, margin=2.5cm]{geometry}

\usepackage[nosumlimits]{amsmath}
\usepackage{amssymb}
\usepackage{amsfonts}
\usepackage{lmodern}
\usepackage{amsthm}

\usepackage{graphicx}
\usepackage{booktabs}
\usepackage{setspace}
\usepackage{titlesec}
\usepackage{needspace}
\usepackage{ragged2e}

\usepackage[authoryear,round]{natbib}
\setcitestyle{aysep={},yysep={; },notesep={, }}
\usepackage{hyperref}
\usepackage{xurl}

\hypersetup{
    colorlinks=true,
    linkcolor=blue,
    citecolor=blue,
    urlcolor=blue
}
\makeatletter
\renewcommand{\tagform@}[1]{\maketag@@@{(\textcolor{blue}{#1})}}
\makeatother
\justifying
\titleformat{\section}[hang]
  {\bfseries\fontsize{13}{18}\selectfont}
  {\makebox[1cm][l]{\thesection.}}
  {0pt}
  {}
\titlespacing*{\section}{0pt}{6pt}{6pt}
\titleformat{\subsection}[hang]
  {\normalfont\fontsize{13}{19.5}\selectfont}
  {\makebox[1cm][l]{\thesubsection}}
  {0pt}
  {}
\titlespacing*{\subsection}{0pt}{6pt}{6pt}
\titleformat{\subsubsection}[hang]
  {\normalfont\itshape\fontsize{12}{18}\selectfont}
  {\makebox[1cm][l]{\thesubsubsection}}
  {0pt}
  {}
\titlespacing*{\subsubsection}{0pt}{6pt}{6pt}

\makeatletter
\pretocmd{\section}{\Needspace{6\baselineskip}}{}{}
\pretocmd{\subsection}{\Needspace{4\baselineskip}}{}{}
\pretocmd{\subsubsection}{\Needspace{4\baselineskip}}{}{}
\makeatother
\everydisplay{\fontsize{12}{18}\selectfont}
\numberwithin{equation}{section}

\let\oldsum\sum
\renewcommand{\sum}{\oldsum\nolimits}
\newcommand{\vect}[1]{\boldsymbol{#1}}
\newcommand{\mat}[1]{\boldsymbol{#1}}
\newcommand{\R}{\mathbb{R}}

\newcommand{\1}{\mathbf{1}}
\newcommand{\0}{\mathbf{0}}

\newcommand{\prooftitle}[3]{%
  \par\vspace{3pt}%
  \phantomsection\label{#3}%
  \Needspace{4\baselineskip}%
  {\noindent\fontsize{12}{18}\selectfont\textbf{#1} #2\par}%
  \vspace{3pt}%
}

\begin{document}
{\noindent\bfseries\fontsize{15}{18}\selectfont Kling-Gupta linear regression\par}\vspace{6pt}

{\noindent\fontsize{12}{18}\selectfont
Hristos Tyralis\textsuperscript{1*}, Georgia Papacharalampous\textsuperscript{2} \par \vspace{3pt}
\textsuperscript{1}Support Command, Hellenic Air Force, Elefsina Air Base, 19 200, Elefsina, Greece \par \vspace{3pt}
(\href{mailto:montchrister@gmail.com}{\texttt{montchrister@gmail.com}}, \href{mailto:hristos@itia.ntua.gr}{\texttt{hristos@itia.ntua.gr}}, \url{https://orcid.org/0000-0002-8932-4997}) \par \vspace{3pt}
\textsuperscript{2}Department of Land, Environment, Agriculture and Forestry, University of Padova, Viale dell'Università 16, 35020, Legnaro, Italy \par \vspace{3pt}
(\href{mailto:papacharalampous.georgia@gmail.com}{\texttt{papacharalampous.georgia@gmail.com}}, \href{mailto:georgia.papacharalampous@unipd.it}{\texttt{georgia.papacharalampous@unipd.it}}, \url{https://orcid.org/0000-0001-5446-954X}) \par \vspace{3pt}
\textsuperscript{*}Corresponding author\par}

{\noindent\fontsize{12}{18}\selectfont\textbf{Abstract}: Kling-Gupta efficiency ($\mathrm{KGE}$) is a model performance evaluation metric widely used in hydrology, but its properties as a statistical estimator have remained unexplored. We formalize the Kling-Gupta loss $L_\mathrm{KG} = (1 - \mathrm{KGE})^2$ in an extremum estimation framework (maximizing $\mathrm{KGE}$) for multiple linear regression. We give explicit formulas showing that Kling-Gupta regression scales the ordinary least squares (OLS) coefficient vector by a variance-inflation factor depending on sample variances and covariances. Its predictions reproduce the training set response variance, unlike OLS's variance reduction, while both maintain the response mean and achieve the same sample correlation. We prove that no estimator simultaneously maximizes Nash-Sutcliffe efficiency ($\mathrm{NSE}$) and $\mathrm{KGE}$: OLS maximizes $\mathrm{NSE}$ but not $\mathrm{KGE}$, whereas Kling-Gupta regression maximizes $\mathrm{KGE}$ at the expense of $\mathrm{NSE}$. We establish almost-sure convergence of the Kling-Gupta estimator to well-defined population limits. The training and test set performance metrics for both estimators converge asymptotically to identical limits (different for OLS vs. Kling-Gupta). In a single-predictor model with fixed intercept, we identify conditions where a global minimum of $L_\mathrm{KG}$ does not exist because of discontinuity at zero slope. This work establishes a mathematical foundation for $\mathrm{KGE}$-based estimation and clarifies its effects on predictive performance in hydrologic modeling.
\par \vspace{3pt}
\textbf{Keywords}: calibration; extremum estimator; hydrologic modeling; Kling-Gupta efficiency; mean squared error; Nash-Sutcliffe efficiency \par}
\section{Introduction}
\label{sec:introduction}
Model estimation (calibration, training) and evaluation (validation, testing) lie at the heart of hydrologic research \citep{beven2025, klemes1986, moriasi2012, moriasi2015b}. For the purpose of estimation and evaluation, a critical step is the selection of a suitable loss function, often referred to as an objective function. The squared error loss is one such function that sees widespread application across hydrology \citep{bennett2013, biondi2012, jackson2019, krause2005, moriasi2007, moriasi2015a}:
\begin{equation}
\label{eq:squared_error_loss}
L_\mathrm{SE}(z, y) := (z - y)^2
\end{equation}
This loss penalizes the error of a prediction $z$ for a random variable $\underline{y}$ when $y$ realizes. In practice, the average empirical counterpart of the squared error loss, the mean squared error ($\mathrm{MSE}$), is computed in the validation (test) set, to compare predictions:
\begin{equation}
\label{eq:mean_squared_error}
\mathrm{MSE}(\vect{z}_n, \vect{y}_n) := (1/n) \sum_{i = 1}^n L_\mathrm{SE}(z_i, y_i) = (1/n) \sum_{i = 1}^n (z_i - y_i)^2
\end{equation}
Here and in the following, $\vect{y}_n = (y_1, \ldots, y_n)^\mathsf{T}$ represents the observations of the random variable $\underline{y}$ (with vector formulation defined in eq.~\eqref{eq:vecdef}) and $\vect{z}_n = (z_1, \ldots, z_n)^\mathsf{T}$ represents the model's predictions. Minimizing the $\mathrm{MSE}$ also functions as a consistent estimator of semiparametric regression models that predict the conditional mean \citep{dimitriadis2024, gneiting2011}.

To compare model predictions across multiple catchments, research has favored homogeneous performance metrics of degree zero (often called dimensionless) whose ranges facilitate straightforward interpretation. The Nash-Sutcliffe efficiency ($\mathrm{NSE}$), introduced by \citep{nash1970}, is arguably the earliest metric in the field to exhibit these properties, mapping to the interval $(-\infty, 1]$:
\begin{equation}
\label{eq:Nash_Sutcliffe_Efficiency}
\mathrm{NSE}(\vect{z}_n, \vect{y}_n) := 1 - \frac{\mathrm{MSE}(\vect{z}_n, \vect{y}_n)}{\mathrm{MSE}(\mu(\vect{y}_n) \1_n, \vect{y}_n)}
\end{equation}
where $\mu(\vect{y}_n)$ represents the sample mean of the vector $\vect{y}_n$, as defined in eq.~\eqref{eq:samplemean} and $\1_n$ represents the all-ones vector, as defined in eq.~\eqref{eq:onesvec}. Interpreting $\mathrm{NSE}$ is straightforward; higher values correspond to better performance. An efficiency of unity indicates perfect performance, whereas an efficiency of zero suggests the model performs as well as the constant mean climatology, with respect to the $\mathrm{MSE}$.

Building upon the legacy of $\mathrm{NSE}$ \citep{melsen2025}, research in the discipline strived to develop performance metrics that retain both homogeneity of degree zero and an upper bound of unity. To achieve this, one might select a loss function bounded below by zero, render its average empirical counterpart homogeneous of degree zero (representing this homogeneous form as $\overline{L}$), for instance by dividing the function by a suitable transformation of the observations and subsequently subtract this quantity from unity, to formulate $1 - \overline{L}$. In hydrology, such formulations are often referred to as efficiencies. Following $\mathrm{NSE}$, the most prominent example is the Kling-Gupta efficiency ($\mathrm{KGE}$), established by \citep{gupta2009}:
\begin{equation}
\label{eq:kge}
\begin{aligned}
\mathrm{KGE}(\vect{z}_n, \vect{y}_n) &:= 1 - \sqrt{(1 - \frac{\mu(\vect{z}_n)}{\mu(\vect{y}_n)})^2 + (1 - \frac{\sigma(\vect{z}_n)}{\sigma(\vect{y}_n)})^2 + (1 - \rho(\vect{z}_n, \vect{y}_n))^2},\\
\mu(\vect{y}_n), \sigma(\vect{y}_n), \sigma(\vect{z}_n) & \in \R\backslash\{0\}
\end{aligned}
\end{equation}
where $\sigma(\vect{y}_n)$ and $\sigma(\vect{z}_n)$ represent the sample standard deviations of the corresponding vectors, defined in eq.~\eqref{eq:samplestd} and $\rho(\vect{z}_n, \vect{y}_n)$ is the Pearson sample correlation between $\vect{z}_n$ and $\vect{y}_n$, defined in eq.~\eqref{eq:pearson}. Following the usual convention, which we adopt throughout the manuscript, we refer to $(1 - \frac{\mu(\vect{z}_n)}{\mu(\vect{y}_n)})^2$ as the bias term, $(1 - \frac{\sigma(\vect{z}_n)}{\sigma(\vect{y}_n)})^2$ as the variability term and $(1 - \rho(\vect{z}_n, \vect{y}_n))^2$ as the correlation term \citep{kling2012}.

Regardless of the reasoning behind the development of the $\mathrm{KGE}$ from the $\mathrm{NSE}$ (detailed in the work by \citep{gupta2009}) and evaluating strictly its mathematical structure in eq.~\eqref{eq:kge}, the metric rewards predictions where $\mu(\vect{z}_n)$ converges to $\mu(\vect{y}_n)$, $\sigma(\vect{z}_n)$ converges to $\sigma(\vect{y}_n)$ and $\rho(\vect{z}_n, \vect{y}_n)$ approaches unity. Expressed differently, $\mathrm{KGE}$ favors predictions that align closely with the observed data in terms of their sample mean and their sample variance (or spread) and that are perfectly linearly correlated with them.

The original paper examined the $\mathrm{KGE}$ in the contexts of both model estimation and prediction evaluation. Although existing literature presents some theoretical discourse on its application as an evaluation metric, a theoretical assessment of its use in estimation remains absent. As a direct result, the wider literature frequently labels the $\mathrm{KGE}$ as an ``informal'' metric. This represents a knowledge gap, particularly given the established link between consistent $M$-estimators for a specific model \citep{huber1964, huber1967, newey1994} and strictly consistent loss functions for a specified functional. By definition, a strictly consistent loss function for a functional achieves its expected minimum when a prediction matches that functional \citep{gneiting2011}; \citep{dimitriadis2024} demonstrate this equivalence between consistent $M$-estimation and strictly consistent loss functions; the latter in the context of prediction evaluation. In practice, when applied to fit a correctly specified and point-identified semiparametric regression model (e.g., minimizing the $\mathrm{MSE}$), the $M$-estimator determines which functional of the dependent variable's distribution (such as the mean) the model predicts, and vice versa \citep{dimitriadis2024}.

Here, we aim to establish some formal theoretical properties of the estimator that minimizes the \textit{Kling-Gupta loss} for a linear model. Focusing on the linear model is the canonical starting point for any new estimation theory. Nearly every major statistical method, from OLS to ridge regression, was first fully understood in the linear setting before being extended to nonlinear systems. To this end, we consider the loss defined by $(1 - \mathrm{KGE})^2$, the square of a negatively oriented version of the $\mathrm{KGE}$. We refer to this formulation as the \textit{Kling-Gupta loss}, formally defining it in eq.~\eqref{eq:Kling_Gupta_loss}. Therefore, we analyze the estimator that minimizes this loss as an extremum estimator. Literature on extremum estimators includes works by \citep{amemiya1973, amemiya1985} and \citep{newey1994}. The class of $M$-estimators is a subset of the class of extremum estimators, which includes, among others, the minimizer of the $\mathrm{MSE}$ (i.e., the least squares estimator). The estimator that minimizes the Kling-Gupta loss is not an $M$-estimator, but it does belong to the wider class of extremum estimators.

Analyzing hydrology-specific metrics as extremum estimators is not a new concept, though it remains rare. \citep{tyralis2025ioa} previously employed this framework to analyze the index of agreement \citep{willmott1981}. We restrict our attention to the case of estimation and prediction evaluation for a single time series (or catchment); even though practical applications frequently span multiple catchments. Applying the $\mathrm{NSE}$ or $\mathrm{KGE}$ to evaluate predictions for multiple time series simultaneously alters the interpretation of these metrics, as empirically demonstrated by \citep{williams2025}. For the $\mathrm{NSE}$, \citep{tyralis2026nse} contribute additional theoretical results within a decision-theoretic setting, illustrating property shifts during multi-series prediction evaluations; however, mirroring this analysis for the $\mathrm{KGE}$ falls outside the scope of this manuscript.

Sections~\ref{sec:kling_gupta_loss_sec} and \ref{sec:kling_gupta_regressions} detail the contributions of this paper, which are as follows:

\noindent\makebox[1cm][l]{(i)} We analyze the extremum estimator defined by minimizing the Kling-Gupta loss function, specifically for linear models (Section~\ref{sec:kling-gupta-linear-regression}), a formulation which we refer to as \textit{Kling-Gupta linear regression}. We formulate explicit, closed-form parameter estimates and demonstrate that Kling-Gupta linear regression operates as a variance-inflation procedure relative to OLS, ensuring the model predictions reproduce the marginal variance of the observations on the training set. We also prove that both estimators have the same sample mean (equal to that of the observations) and the same correlation between predictions and observations on the training set.

\noindent\makebox[1cm][l]{(ii)} We characterize the predictive performance of the estimator across multiple metrics, establishing the exact theoretical trade-offs between $\mathrm{MSE}$, $\mathrm{NSE}$, $\mathrm{KGE}$, and the Kling-Gupta loss within the training set (Section~\ref{sec:comparative_training_data}). We extend these findings to asymptotic settings, addressing both the infinite training data limit (Section~\ref{sec:comparative_infinite_training_data}) and the associated infinite test set data performance (Section~\ref{sec:comparative_infinite_training_test_data}).

\noindent\makebox[1cm][l]{(iii)} We prove the uniqueness of the global minimizer for the Kling-Gupta loss (Section~\ref{sec:kling_gupta_loss_sec}). Furthermore, we examine constrained estimation settings, estimating parameters in models with fixed regression coefficients (Section~\ref{sec:fixed_regression_coefficient}) and in single-predictor linear models with fixed intercepts (Section~\ref{sec:fixed_intercept_kg_linear_regression}). In the latter case, we show that a global minimum may fail to exist for the slope parameter, a consequence of the loss function's domain restrictions and its inherent discontinuity at a slope of zero.

The remainder of the manuscript is structured as follows. Section~\ref{sec:theoretical_background} outlines the prerequisite theoretical background. Section~\ref{sec:kling_gupta_regressions} details the primary statistical contributions, while Section~\ref{sec:applications} demonstrates practical hydrologic applications for these theoretical findings. Section~\ref{sec:discussion} synthesizes the results, situating them within the context of extant literature in the discipline and discusses their implications. We conclude in Section~\ref{sec:conclusions}. Appendix~\ref{app:notation} establishes the notation, Appendix~\ref{app:proofs} details all mathematical proofs and Appendix~\ref{app:software} outlines the statistical software used. All computations required to reproduce this manuscript are fully accessible as supplementary information.
\section{Theoretical background}
\label{sec:theoretical_background}
This section establishes the theoretical concepts required to analyze the Kling-Gupta loss as an extremum estimator. The exposition proceeds from strictly consistent loss functions (Section~\ref{sec:strictly_consistent_loss_functions}) through the framework of extremum and $M$-estimators (Section~\ref{sec:extremum_and_m_estimators}) to the equivalence between strictly consistent loss functions and consistent $M$-estimators (Section~\ref{sec:estimation_and_evaluation}), a correspondence that underpins the analysis throughout. Hydrologic models are then situated within the semiparametric regression framework (Section~\ref{sec:hydrologic_models}), where model predictions are understood to represent specific statistical functionals of the conditional streamflow distribution. Common loss functions and performance metrics employed in hydrology, including the squared error loss and the $\mathrm{NSE}$, are illustrated in Section~\ref{sec:squared_error_and_nse}. Subsequently, Section~\ref{sec:kling_gupta_loss_sec} introduces the Kling-Gupta loss, a negatively oriented transformation of the $\mathrm{KGE}$, and establishes its elementary properties. The section concludes with the mathematical formulation of the linear model (Section~\ref{sec:linear-model}) and OLS linear regression (Section~\ref{sec:estimating-linear-model-se}), both of which underpin the theoretical developments of Section~\ref{sec:kling_gupta_regressions}.
\subsection{Strictly consistent loss functions}
\label{sec:strictly_consistent_loss_functions}
A favorable property of any loss function is its (strict) consistency. To examine this concept, it is necessary to first establish the definition of a statistical functional. Let $\underline{y}$ be a random variable and let a realization of this variable be $y$. The notation $\underline{y} \sim F_{\underline{y}}$ indicates that the random variable $\underline{y}$ follows the cumulative distribution function (CDF) $F_{\underline{y}}$, which is defined as:
\begin{equation}
F_{\underline{y}}(y) := P(\underline{y} \leq y)
\end{equation}
A one-dimensional statistical functional $T$ (often referred to simply as a functional) is a mapping \citep{gneiting2011}:
\begin{equation}
T: \mathcal{F} \rightarrow \mathcal{P}(D), F_{\underline{y}} \mapsto T(F_{\underline{y}}) \subseteq D
\end{equation}
Here, $\mathcal{F}$ represents a specified class of probability distributions. This mapping assigns every distribution $F_{\underline{y}} \in \mathcal{F}$ to a corresponding subset $T(F_{\underline{y}})$ within a set $D \subseteq \R$, thereby constituting an element of the power set $\mathcal{P}(D)$.

A given loss function $L$ is said to be $\mathcal{F}$-consistent for a specified functional $T$ if its expectation satisfies the following inequality for all probability distributions in the class $\mathcal{F}$ \citep{gneiting2011, murphy1985}:
\begin{equation}
\label{eq:consistency_definition}
\mathbb{E}_{F_{\underline{y}}}[L(t, \underline{y})] \leq \mathbb{E}_{F_{\underline{y}}}[L(z, \underline{y})] \forall F_{\underline{y}} \in \mathcal{F}, t \in T(F_{\underline{y}}), z \in D
\end{equation}
Furthermore, $L$ is characterized as strictly $\mathcal{F}$-consistent if it satisfies $\mathcal{F}$-consistency and equality in eq.~\eqref{eq:consistency_definition} holds strictly when the prediction $z \in T(F_{\underline{y}})$.

To illustrate the practical significance of strict consistency, suppose a modeler is tasked with forecasting the expected value (mean) of the random variable $\underline{y}$. The squared error loss is strictly consistent for the mean \citep{gneiting2011}. Therefore, a modeler who correctly predicts the true mean will, in expectation, achieve the minimum of this specific loss function.

Because the expected loss expressed in eq.~\eqref{eq:consistency_definition} is a population-level characteristic, practical applications that assess and compare predictive performance must use its empirical counterpart. This is formulated as the sample average loss:
\begin{equation}
\label{eq:sample_average_loss}
\overline{L}(\vect{z}_n, \vect{y}_n) := (1/n) \sum_{i = 1}^n L(z_i, y_i)
\end{equation}
where predictions with a lower value of $\overline{L}$ are ranked as better.
\subsection{Extremum and \texorpdfstring{$M$}{M}-estimators}
\label{sec:extremum_and_m_estimators}
The theoretical framework for extremum estimators was established by \citep{amemiya1973, amemiya1985} and \citep{newey1994}. Let a loss function $L(\vect{\theta}, \vect{y}_n)$ be designed to measure the discrepancy between a parameter vector $\vect{\theta}$ and a realized vector $\vect{y}_n$ generated from a random variable $\underline{y}$ that follows the CDF $F_{\underline{y}}$. We assume that the true parameter value $\vect{\theta}_0$ of $F_{\underline{y}}$ lies within the parameter space $\vect{\Theta}$ and corresponds to a specific statistical functional of the distribution $F_{\underline{y}}$ (such as the expectation, $\mathbb{E}_{F_{\underline{y}}}[\underline{y}]$). An extremum estimator $\widehat{\vect{\theta}}(\underline{\vect{y}}_n)$ for the parameter $\vect{\theta}_0$ is defined as the solution to the following optimization problem:
\begin{equation}
\label{eq:extremum_estimator}
\widehat{\vect{\theta}}(\underline{\vect{y}}_n) := \underset{\vect{\theta} \in \vect{\Theta}}{\operatorname*{arg\,min}} L(\vect{\theta}, \underline{\vect{y}}_n)
\end{equation}
The statistical behavior and asymptotic properties of the estimator $\widehat{\vect{\theta}}(\underline{\vect{y}}_n)$ are inherently governed by the sample size $n$, which appears in the empirical observations $\vect{y}_n$, as well as by the specific form of the loss function $L$. The formulation in eq.~\eqref{eq:extremum_estimator} refers to models without predictors, e.g., parameter estimation for probability distributions. Semiparametric regression models will be introduced later in Section~\ref{sec:hydrologic_models}, while extremum estimators for such models will be formulated in Section~\ref{sec:kling_gupta_loss_sec}.

While the loss function $L$ can assume a wide variety of forms, settings characterized by additivity have been the subject of attention in the statistical literature. Estimators arising from such additive loss functions are classified as $M$-estimators, a conceptual class introduced by \citep{huber1964, huber1967}. In the context of a scalar parameter $\theta_0$, an $M$-estimator takes the following form:
\begin{equation}
\label{eq:scalar_M_estimator}
\widehat{\theta}(\underline{\vect{y}}_n) := \underset{\theta \in \Theta}{\operatorname*{arg\,min}} (1/n) \sum_{i = 1}^n L(\theta, \underline{y}_i)
\end{equation}
An illustration of this class is the estimator arising from the squared error loss, which corresponds to the consistent $M$-estimator of the population expectation $\mathbb{E}_{F_{\underline{y}}}[\underline{y}]$. This estimator is expressed as:
\begin{equation}
\widehat{\mu}(\underline{\vect{y}}_n) := \underset{\theta \in \Theta}{\operatorname*{arg\,min}} (1/n) \sum_{i = 1}^n (\theta - \underline{y}_i)^2 = \underset{\theta \in \Theta}{\operatorname*{arg\,min}} \mathrm{MSE}(\theta \1_n, \underline{\vect{y}}_n), \Theta \subseteq \R
\end{equation}
A favorable property for any statistical estimator is consistency. An estimator $\widehat{\vect{\theta}}(\underline{\vect{y}}_n)$ of true underlying parameter $\vect{\theta}_0$ is consistent if, as the sample size increases, it converges in probability toward $\vect{\theta}_0$, meaning $\widehat{\vect{\theta}}(\underline{\vect{y}}_n) \overset{P}{\longrightarrow} \vect{\theta}_0$. As a practical example, $\widehat{\mu}(\underline{\vect{y}}_n)$ is a consistent estimator for the expectation $\mathbb{E}_{F_{\underline{y}}}[\underline{y}]$ as demonstrated by \citep{dimitriadis2024} among others.
\subsection{Estimation and evaluation: A unified framework}
\label{sec:estimation_and_evaluation}
Although model estimation (training) and prediction evaluation (testing) are often treated as distinct procedural steps, they are unified within the framework of $M$-estimators. Specifically, \citep{dimitriadis2024} established a mathematical equivalence: an $M$-estimator of a true parameter $\theta_0$ is consistent if and only if the estimator's generating loss function (as formulated in eq.~\eqref{eq:scalar_M_estimator}) is strictly consistent for the statistical functional $T_0$ corresponding to $\theta_0$.

It is important to distinguish between the statistical consistency of an estimator and the strict consistency of a loss function, as they apply to separate procedures. Nevertheless, their equivalence implies that minimizing a strictly consistent loss function during the training of a semiparametric regression model incentivizes the model to predict the true underlying functional. This outcome arises because the expected loss $\mathbb{E}_{F_{\underline{y}}}[L(z, \underline{y})]$ attains its global minimum when the prediction $z$ aligns with the functional's true value $T_0$.

Given that a single statistical functional can be associated with an infinite family of strictly consistent loss functions, the specific metric chosen for evaluation must be disclosed to the modeler a priori. As demonstrated by \citep{patton2020}, model rankings are highly sensitive to this selection in practical applications, particularly when the data-generating process is unknown. Moreover, different loss functions belonging to the same family (i.e., all strictly consistent for a given functional) can lead to entirely divergent hierarchies of model performance.
\subsection{Hydrologic models as semiparametric regression models}
\label{sec:hydrologic_models}
Let the random variable $\underline{y}$ represent the streamflow process. We define a hydrologic model $\mathcal{HM}$, which is characterized by a parameter vector $\vect{\theta}$ belonging to a specified parameter space $\vect{\Theta}$. The model processes a vector of predictor variables $\vect{x}_p$, such as precipitation, temperature or evapotranspiration, to issue a specific prediction $z$. Importantly, $z$ is not a direct physical realization of streamflow; rather, it is a specific statistical functional computed from the streamflow conditional probability distribution of $\underline{y}|\underline{\vect{x}}_p$ given the meteorological forcing data. This predictive relationship is formalized as:
\begin{equation}
\label{eq:hydrologic_model}
z = \mathcal{HM}(\vect{x}_p; \vect{\theta})
\end{equation}
Such structural representations of hydrologic systems are well-documented in the literature (e.g., \citep{montanari2012} present a distributional regression equivalent of model \eqref{eq:hydrologic_model}). \citep{vrugt2024} specifically refers to such predictions $z$ of hydrologic models as ``hydrograph functionals''. In this context, the standard modeling workflow comprises two distinct phases: model estimation (calibration) and model testing (validation), as detailed by \citep{klemes1986} and \citep{beven2025}.

For a hydrologic model that is correctly specified and point-identified, a unique true parameter vector $\vect{\theta}_0 \in \vect{\Theta}$ exists such that the model's prediction perfectly matches the target functional:
\begin{equation}
T(F_{\underline{y}|\underline{\vect{x}}_p}) = \mathcal{HM}(\underline{\vect{x}}_p; \vect{\theta}_0)
\end{equation}
While this theoretical framework was explicitly formalized by \citep{dimitriadis2024} for general semiparametric regression analysis, it applies directly to hydrologic modeling architectures.

In practice, even if the structure of $\mathcal{HM}$ is assumed to represent the true data-generating process, the true parameter vector $\vect{\theta}_0$ remains unknown and must be estimated from observations. The standard estimation approach involves selecting a specific loss function $L$, to construct an $M$-estimator for $\vect{\theta}_0$. Given a sequence of predictor data $\mat{X}_{n \times p}$ (the subscript $n \times p$ indicates the matrix dimensions, as defined in eq.~\eqref{eq:matrix}), where each row of $\mat{X}_{n \times p}$ corresponds to a single observation of the predictor vector $\vect{x}_p$, and corresponding observed streamflow targets $\vect{y}_n$, the estimator is found by solving the following optimization problem:
\begin{equation}
\widehat{\vect{\theta}}(\underline{\mat{X}}_{n \times p}, \underline{\vect{y}}_n) := \underset{\vect{\theta} \in \vect{\Theta}}{\operatorname*{arg\,min}} (1/n) \sum_{i = 1}^n L(\mathcal{HM}(\underline{\mat{X}}_{i, \bullet}^\mathsf{T}; \vect{\theta}), \underline{y}_i)
\end{equation}
where $\mat{X}_{i, \bullet}$ is the $i$\textsuperscript{th} row of $\mat{X}_{n \times p}$ (see eq.~\eqref{eq:matrixrow}). The choice of the loss function $L$ determines the specific statistical functional that the model is trained to predict \citep{dimitriadis2024}, as explained in Section~\ref{sec:estimation_and_evaluation}. For instance, minimizing the squared error loss incentivizes the model to predict the conditional expectation.

Because the true underlying data-generating mechanism of a catchment is never known, modelers typically evaluate a range of candidate models. The testing phase involves ranking these alternative models against a test set to identify the best predictive structure \citep{klemes1986}. Loss functions are equally central to this evaluation phase; models are compared based on their average empirical loss, computed as in eq.~\eqref{eq:sample_average_loss}. The candidate model that achieves the lowest average loss on the test set is then selected as the best model for predicting that specific functional \citep{gneiting2011}.
\subsection{Squared error loss and Nash-Sutcliffe efficiency}
\label{sec:squared_error_and_nse}
The squared error loss, defined in eq.~\eqref{eq:squared_error_loss}, is an important loss function and belongs to the broader family of Bregman loss functions \citep{banerjee2005, gneiting2011, patton2011, reichelstein1984, saerens2000, savage1971}. As established in Section~\ref{sec:strictly_consistent_loss_functions}, this loss function is strictly consistent for the mean functional. Substituting the squared error loss into eq.~\eqref{eq:sample_average_loss} gives the $\mathrm{MSE}$, defined in eq.~\eqref{eq:mean_squared_error}. Therefore, minimizing the $\mathrm{MSE}$ during model estimation functions as an $M$-estimator that incentivizes a correctly specified semiparametric regression model to predict the conditional expectation.

To facilitate interpretable comparisons across models, the $\mathrm{NSE}$, defined in eq.~\eqref{eq:Nash_Sutcliffe_Efficiency}, is widely adopted as a skill score. The $\mathrm{NSE}$ normalizes the $\mathrm{MSE}$ of a model's predictions by the $\mathrm{MSE}$ of the mean climatology $\mu(\vect{y}_n)$, thereby transforming the squared error loss into a homogeneous efficiency metric of degree zero that maps to the interval $(-\infty, 1]$. Within this framework, an $\mathrm{NSE}$ of unity indicates perfect predictive performance, a value of zero signifies that the model performs equivalently to the mean climatology and negative values indicate that the mean climatology outperforms the model. As a skill score, the $\mathrm{NSE}$ thus quantifies predictive skill relative to the mean climatology benchmark \citep{gneiting2023, murphy1988}. For evaluations involving a single time series, the $\mathrm{NSE}$ retains the model rankings established by the $\mathrm{MSE}$ but is often preferred for its interpretability \citep{tyralis2026nse}. $\mathrm{MSE}$ and $\mathrm{NSE}$ lead to identical parameter estimates when fitting a semiparametric regression model to a single time series.
\subsection{Kling-Gupta loss}
\label{sec:kling_gupta_loss_sec}
Unlike the $\mathrm{NSE}$, the $\mathrm{KGE}$ is not a skill score, nor does it function as an $M$-estimator for single time series problems, because, unlike the $\mathrm{MSE}$, it lacks an additive structure. As indicated in Section~\ref{sec:introduction}, we adopt a bijective transformation of the $\mathrm{KGE}$, which we refer to as the Kling-Gupta loss, to facilitate our analysis. This loss function is defined as:
\begin{equation}
\label{eq:Kling_Gupta_loss}
\begin{aligned}
L_\mathrm{KG}(\vect{z}_n, \vect{y}_n) &:= (1 - \frac{\mu(\vect{z}_n)}{\mu(\vect{y}_n)})^2 + (1 - \frac{\sigma(\vect{z}_n)}{\sigma(\vect{y}_n)})^2 + (1 - \rho(\vect{z}_n, \vect{y}_n))^2,\\
\mu(\vect{y}_n), \sigma(\vect{y}_n), \sigma(\vect{z}_n) &\in \R\backslash\{0\}
\end{aligned}
\end{equation}
For single time series problems, working with $\mathrm{KGE}$ or $L_\mathrm{KG}$ is equivalent, as the two are related by the bijective transformation:
\begin{equation}
\label{eq:KGE_LKG_relation}
L_\mathrm{KG}(\vect{z}_n, \vect{y}_n) = (1 - \mathrm{KGE}(\vect{z}_n, \vect{y}_n))^2, \mathrm{KGE}(\vect{z}_n, \vect{y}_n) = 1 - \sqrt{L_\mathrm{KG}(\vect{z}_n, \vect{y}_n)}
\end{equation}
Eliminating the square root from the $\mathrm{KGE}$ expression facilitates analytical computations (e.g., differentiation) without compromising interpretability. $L_\mathrm{KG}$ possesses the structure of a loss function, meaning it is negatively oriented and bounded below by zero. A value of $L_\mathrm{KG} = 0$ corresponds to $\mathrm{KGE} = 1$. Unconditional forecasts from a probability distribution with expectation $\mathbb{E}_{F_{\underline{z}}}[\underline{z}] = \mu(\vect{y}_n)$ and variance $\mathrm{Var}_{F_{\underline{z}}}[\underline{z}] = \sigma^2(\vect{y}_n)$ (defined in eq.~\eqref{eq:variance}) give, in expectation, $L_\mathrm{KG} = 1$ and $\mathrm{KGE} = 0$.

As established in \hyperref[proof:b1]{Proof~B.1}, $L_\mathrm{KG}$ attains its unique minimum value of zero if and only if $\vect{z}_n = \vect{y}_n$. To the best of our knowledge, this result, specifically the demonstration of equivalence, has not previously appeared in the $\mathrm{KGE}$-related literature. Moreover, $L_\mathrm{KG}$ exhibits homogeneity of degree zero, as expressed by:
\begin{equation}
L_\mathrm{KG}(c \vect{z}_n, c \vect{y}_n) = L_\mathrm{KG}(\vect{z}_n, \vect{y}_n), c \neq 0
\end{equation}
The formulation in eq.~\eqref{eq:extremum_estimator}, with $L$ replaced by $L_\mathrm{KG}$, corresponds to extremum parameter estimation problems with $\mathrm{KGE}$ that lack predictors and involve a scalar parameter $\theta$ (e.g., estimating a parameter of a probability distribution). In such cases, the optimization would take the form of eq.~\eqref{eq:extremum_estimator}, with $\vect{z}_n = \theta \1_n$. However, the Kling-Gupta loss is undefined when $\vect{z}_n = \theta \1_n$, because this implies $\sigma(\vect{z}_n) = 0$, which violates the domain restriction specified in eq.~\eqref{eq:Kling_Gupta_loss}. Therefore, the estimation problem cannot be posed within this framework.

Nevertheless, $L_\mathrm{KG}$ can function as an extremum estimator, according to eq.~\eqref{eq:extremum_estimator}, for semiparametric regression models. In this context, for a semiparametric model of the form introduced in Section~\ref{sec:hydrologic_models}, we refer to the extremum estimator as \textit{Kling-Gupta extremum estimator}, which takes the following form:
\begin{equation}
\widehat{\vect{\theta}}(\underline{\mat{X}}_{n \times p}, \underline{\vect{y}}_n) := \underset{\vect{\theta} \in \vect{\Theta}}{\operatorname*{arg\,min}} L_\mathrm{KG}(\vect{\mathcal{HM}}_n(\underline{\mat{X}}_{n \times p}; \vect{\theta}), \underline{\vect{y}}_n)
\end{equation}
where $\vect{\mathcal{HM}}_n$ represents the $n$ predictions of the hydrologic model $\mathcal{HM}$ specified in eq.~\eqref{eq:hydrologic_model} and $\mat{X}_{n \times p}$ is the predictor matrix, each row of which corresponds to a single observation of the predictor vector $\vect{x}_p$:
\begin{equation}
\vect{\mathcal{HM}}_n(\mat{X}_{n \times p}; \vect{\theta}) = (\mathcal{HM}(\mat{X}_{1, \bullet}^\mathsf{T}; \vect{\theta}), \ldots, \mathcal{HM}(\mat{X}_{n, \bullet}^\mathsf{T}; \vect{\theta}))^\mathsf{T}
\end{equation}
\subsection{The linear model}
\label{sec:linear-model}
This section introduces the general formulation of a linear model, which is independent of the specific estimator used for the model parameters. Sections~\ref{sec:estimating-linear-model-se} and \ref{sec:kling_gupta_regressions} then show how these baseline expressions take distinct forms under different estimation methods.
\subsubsection{Linear model with multiple predictors}
\label{sec:linear_model_multiple_predictors}
Let $\underline{y}$ be a scalar response random variable and let $\underline{\vect{x}}_p$ represent a $p$-dimensional predictor vector. A linear model designed to predict a specific statistical functional of the response variable $\underline{y}$ is defined as follows:
\begin{equation}
\label{eq:linear_model}
z = \vect{a}_p^\mathsf{T} \vect{x}_p + b_p
\end{equation}
where $\vect{a}_p \in \R^p$ is the vector of regression coefficients (slopes) and $b_p \in \R$ is the intercept.

Let
\begin{equation}
\label{eq:multi_linear_model_parameter}
\vect{\theta}_{1 \times (p + 1)} = (a_1, \ldots, a_p, b_p) \in \R^{1 \times (p + 1)}
\end{equation}
be the parameter vector of the coefficients and the intercept. Given a sample of $n$ observations, let $\vect{y}_n$ be the vector of realizations of the response variable and let $\mat{X}_{n \times p}$ be the predictor matrix. Each row of $\mat{X}_{n \times p}$ corresponds to a single observation of the predictor vector $\vect{x}_p$. We assume the sample size is sufficiently large such that $n \geq p$. We further assume that $\mat{X}_{n \times p}$ has full column rank and that the all-ones vector $\1_n$ does not lie in its column space.

Because the linear model in eq.~\eqref{eq:linear_model} issues scalar predictions, the complete vector of predictions for all $n$ observations can be expressed compactly in matrix form as:
\begin{equation}
\label{eq:lm_training_predictions}
\vect{z}_n = \widetilde{\mat{X}}_{n \times (p + 1)} \vect{\theta}_{1 \times (p + 1)}^\mathsf{T} = \mat{X}_{n \times p} \vect{a}_p + b_p \1_n
\end{equation}
where $\vect{z}_n$ is the $n$-dimensional vector of predictions and the augmented predictor matrix $\widetilde{\mat{X}}_{n \times (p + 1)}$ is constructed by appending a column of ones to $\mat{X}_{n \times p}$:
\begin{equation}
\widetilde{\mat{X}}_{n \times (p + 1)} := \begin{bmatrix}
\mat{X}_{n \times p} & \1_n
\end{bmatrix}
\end{equation}
To facilitate variance and covariance computations, we define the centered (mean-subtracted) predictor matrix as (eq.~\eqref{eq:centmat}):
\begin{equation}
\mat{X}_{n \times p, \mathrm{c}} = \mat{P}_{n \times n}\mat{X}_{n \times p}
\end{equation}
where $\mat{P}_{n \times n}$ is the centering matrix defined in eq.~\eqref{eq:centeringmat}. Applying this centering operation to the predictions leads to the centered prediction vector $\vect{z}_{n \mathrm{c}}$ (defined in eq.~\eqref{eq:centvec}):
\begin{equation}
\label{eq:zcent}
\vect{z}_{n \mathrm{c}} = \mat{P}_{n \times n} \vect{z}_n = \mat{P}_{n \times n}(\mat{X}_{n \times p} \vect{a}_p + b_p \1_n) = \mat{X}_{n \times p, \mathrm{c}} \vect{a}_p
\end{equation}
The sample mean of the predictions $\mu(\vect{z}_n)$ is given by:
\begin{equation}
\label{eq:mean_predictions_multi}
\mu(\vect{z}_n) = (\vect{\mu}(\mat{X}_{n \times p}))^\mathsf{T} \vect{a}_p + b_p
\end{equation}
where $\vect{\mu}(\mat{X}_{n \times p})$ is the vector of component-wise sample means for the columns of $\mat{X}_{n \times p}$, defined as in eq.~\eqref{eq:matrixmean}:
\begin{equation}
\vect{\mu}(\mat{X}_{n \times p}) = (\mu(\mat{X}_{\bullet ,1}), \ldots, \mu(\mat{X}_{\bullet ,p}))^\mathsf{T}
\end{equation}
By definition, the sample mean of the centered prediction vector $\vect{z}_{n \mathrm{c}}$ is strictly zero (see eq.~\eqref{eq:meancentvec}):
\begin{equation}
\mu(\vect{z}_{n \mathrm{c}}) = 0
\end{equation}
The sample standard deviation of the predictions follows from eqs.~\eqref{eq:zcent} and \eqref{eq:samplestd}:
\begin{equation}
\label{eq:sigma_predictions}
\sigma(\vect{z}_n) = \sqrt{(\vect{z}_{n \mathrm{c}}^\mathsf{T} \vect{z}_{n \mathrm{c}})/n} = \sqrt{(\vect{a}_p^\mathsf{T} \mat{X}_{n \times p, \mathrm{c}}^\mathsf{T} \mat{X}_{n \times p, \mathrm{c}} \vect{a}_p)/n} = \sqrt{\vect{a}_p^\mathsf{T} \mat{S}_{p \times p} \vect{a}_p}
\end{equation}
where $\mat{S}_{p \times p}$ represents the sample covariance matrix of the predictor variables:
\begin{equation}
\label{eq:covariance_matrix}
\mat{S}_{p \times p} = (\mat{X}_{n \times p, \mathrm{c}}^\mathsf{T} \mat{X}_{n \times p, \mathrm{c}})/n
\end{equation}
The sample covariance matrix $\mat{S}_{p \times p}$ is symmetric and square. Given our prior assumption that $\mat{X}_{n \times p}$ has full column rank and that the all-ones vector $\1_n$ does not lie in its column space, Property~4 in \citep{amrhein1996} ensures that the centered matrix $\mat{X}_{n \times p, \mathrm{c}}$ also maintains full column rank. Therefore, the cross-product matrix $\mat{X}_{n \times p, \mathrm{c}}^\mathsf{T}\mat{X}_{n \times p, \mathrm{c}}$ is both full rank and positive definite \cite[p.~138]{gentle2024}. This positive definiteness ensures the matrix $\mat{S}_{p \times p}$ is nonsingular \cite[p.~122]{gentle2024} and therefore invertible \cite[p.~129]{gentle2024}.

Finally, the Pearson sample correlation between the model predictions and the realizations is calculated using eqs.~\eqref{eq:zcent} and \eqref{eq:pearson}. The vector inequality in eq.~\eqref{eq:correlation_predictions_multi} follows the convention of eq.~\eqref{eq:vecne}:
\begin{equation}
\label{eq:correlation_predictions_multi}
\begin{aligned}
\rho(\vect{z}_n, \vect{y}_n) &= \frac{\vect{z}_{n \mathrm{c}}^\mathsf{T} \vect{y}_{n \mathrm{c}}}{n \sigma(\vect{z}_n) \sigma(\vect{y}_n)} = \frac{\vect{a}_p^\mathsf{T} \mat{X}_{n \times p, \mathrm{c}}^\mathsf{T} \vect{y}_{n \mathrm{c}}}{n \sigma(\vect{z}_n) \sigma(\vect{y}_n)} = \frac{\vect{a}_p^\mathsf{T} \vect{s}_p}{\sigma(\vect{z}_n) \sigma(\vect{y}_n)}\\
&= \frac{\vect{a}_p^\mathsf{T} \vect{s}_p}{\sqrt{\vect{a}_p^\mathsf{T} \mat{S}_{p \times p} \vect{a}_p} \sigma(\vect{y}_n)}, \vect{a}_p \neq \0_p
\end{aligned}
\end{equation}
where $\vect{s}_p$ represents the sample cross-covariance vector between the centered predictors and the centered response variable and $\0_p$ is the zero vector, defined in eq.~\eqref{eq:zerovec}:
\begin{equation}
\label{eq:cross_covariance_vector}
\vect{s}_p = (\mat{X}_{n \times p, \mathrm{c}}^\mathsf{T} \vect{y}_{n \mathrm{c}})/n = (\mat{X}_{n \times p}^\mathsf{T} \mat{P}_{n \times n} \vect{y}_n)/n
\end{equation}
\subsubsection{Linear model with a single predictor}
\label{sec:linear-model-single-predictor}
This section addresses the simplified special case involving only a single predictor variable, that arises from eq.~\eqref{eq:linear_model} for $p = 1$. Let $\vect{x}_n$ represent the vector of predictor values, corresponding directly to the single column of the predictor matrix $\mat{X}_{n \times 1}$ discussed in Section~\ref{sec:linear_model_multiple_predictors}. The linear model formulation simplifies to:
\begin{equation}
\label{eq:single_pred_model}
z = ax + b
\end{equation}
where $a \in \R$ is the slope (regression) coefficient and $b \in \R$ is the intercept.

Let
\begin{equation}
\vect{\theta}_{1 \times 2} = (a, b) \in \R^{1 \times 2}
\end{equation}
be the parameter vector of the slope and the intercept. Accordingly, the complete vector of predictions for the sample can be expressed concisely as:
\begin{equation}
\label{eq:single_pred_model_vec}
\vect{z}_n = \widetilde{\mat{X}}_{n \times 2} \vect{\theta}_{1 \times 2}^\mathsf{T} = a \vect{x}_n + b \1_n
\end{equation}
The sample statistics of these predictions scale with the chosen parameters:
\begin{equation}
\mu(\vect{z}_n) = a \mu(\vect{x}_n) + b
\end{equation}
\begin{equation}
\sigma(\vect{z}_n) = |a|\sigma(\vect{x}_n)
\end{equation}
\begin{equation}
\rho(\vect{z}_n, \vect{y}_n) = \operatorname{sign}(a)\rho(\vect{x}_n, \vect{y}_n), a \neq 0
\end{equation}
where the $\operatorname{sign}$ function is defined in eq.~\eqref{eq:sign}.
\subsection{Estimating the linear model with the squared error loss}
\label{sec:estimating-linear-model-se}
We examine the consequences of training the linear model introduced in Section~\ref{sec:linear-model} using the squared error loss, because these linear models will underpin much of the subsequent theory development. This estimation procedure is referred to in the literature as OLS linear regression, and its properties are well established \citep{gentle2024}. As discussed in Section~\ref{sec:estimation_and_evaluation}, assuming that the data-generating process follows a linear model, training the model with the squared error loss ensures that it predicts the conditional mean of the dependent variable.
\subsubsection{Ordinary least squares linear regression}
\label{sec:ols-linear-regression}
We define the OLS estimate of the linear model in eq.~\eqref{eq:linear_model} as the parameter vector that minimizes the $\mathrm{MSE}$:
\begin{equation}
\widehat{\vect{\theta}}_{1 \times (p + 1), \mathrm{OLS}}(\mat{X}_{n \times p}, \vect{y}_n) := \underset{\vect{\theta}_{1 \times (p + 1)} \in \vect{\Theta}}{\operatorname*{arg\,min}} \mathrm{MSE}(\widetilde{\mat{X}}_{n \times (p + 1)}\vect{\theta}_{1 \times (p + 1)}^\mathsf{T}, \vect{y}_n)
\end{equation}
where $\vect{\theta}_{1 \times (p + 1)}$ refers to the parameters of the linear model as specified in eq.~\eqref{eq:multi_linear_model_parameter}. The estimate of the regression coefficients $\widehat{\vect{a}}_{p, \mathrm{OLS}}$ is given by \cite[p.~441]{gentle2024}:
\begin{equation}
\label{eq:ols_slope_multi}
\widehat{\vect{a}}_{p, \mathrm{OLS}} = (\mat{X}_{n \times p, \mathrm{c}}^\mathsf{T} \mat{X}_{n \times p, \mathrm{c}})^{-1}\mat{X}_{n \times p,\mathrm{c}}^\mathsf{T} \vect{y}_{n \mathrm{c}} = \mat{S}_{p \times p}^{-1} \vect{s}_p
\end{equation}
and the corresponding estimate of the intercept parameter is \cite[p.~441]{gentle2024}:
\begin{equation}
\label{eq:ols_intercept_multi}
\widehat{b}_{p, \mathrm{OLS}} = \mu(\vect{y}_n) - (\vect{\mu}(\mat{X}_{n \times p}))^\mathsf{T} \widehat{\vect{a}}_{p, \mathrm{OLS}}
\end{equation}
The predictions of the OLS linear regression model on the training set take the form:
\begin{equation}
\label{eq:ols_predictions_multi}
\vect{z}_{n, \mathrm{OLS}} = \mat{X}_{n \times p} \widehat{\vect{a}}_{p, \mathrm{OLS}} + \widehat{b}_{p, \mathrm{OLS}} \1_n
\end{equation}
Given a test predictor matrix $\mat{X}_{k \times p, \mathrm{test}}$, the predictions from the estimated OLS linear regression model on the test set are given by:
\begin{equation}
\label{eq:ols_test_predictions_multi}
\vect{z}_{k, \mathrm{OLS}, \mathrm{test}} = \mat{X}_{k \times p, \mathrm{test}} \widehat{\vect{a}}_{p, \mathrm{OLS}} + \widehat{b}_{p, \mathrm{OLS}} \1_k
\end{equation}
From the intercept eq.~\eqref{eq:ols_intercept_multi} and the prediction eq.~\eqref{eq:ols_predictions_multi}, it follows directly that the sample mean of the predictions on the training set is:
\begin{equation}
\label{eq:ols_mean_property}
\mu_\mathrm{OLS} := \mu(\vect{z}_{n, \mathrm{OLS}}) = \mu(\vect{y}_n)
\end{equation}
We write the standard deviation of the predictions as:
\begin{equation}
\label{eq:ols_sigma_def}
\sigma_\mathrm{OLS} := \sigma(\vect{z}_{n, \mathrm{OLS}})
\end{equation}
Because the OLS predictions have the same mean as the observations (see eq.~\eqref{eq:ols_mean_property}), substituting the intercept estimate from eq.~\eqref{eq:ols_intercept_multi} into eq.~\eqref{eq:ols_predictions_multi} gives the centered predictions:
\begin{equation}
\label{eq:ols_centered_predictions}
\vect{z}_{n \mathrm{c}, \mathrm{OLS}} = \mat{X}_{n \times p, \mathrm{c}} \widehat{\vect{a}}_{p, \mathrm{OLS}}
\end{equation}
Then, from eqs.~\eqref{eq:sigma_predictions} and \eqref{eq:ols_sigma_def}, $\sigma_\mathrm{OLS}$ can be expanded as:
\begin{equation}
\label{eq:ols_sigma_expansion}
\begin{aligned}
\sigma_{\mathrm{OLS}} &= \sqrt{(\vect{z}_{n \mathrm{c}, \mathrm{OLS}}^\mathsf{T} \vect{z}_{n \mathrm{c}, \mathrm{OLS}})/n} = \sqrt{(1/n) \widehat{\vect{a}}_{p, \mathrm{OLS}}^\mathsf{T} \mat{X}_{n \times p, \mathrm{c}}^\mathsf{T} \mat{X}_{n \times p, \mathrm{c}}\widehat{\vect{a}}_{p, \mathrm{OLS}}} \\ &= \sqrt{\widehat{\vect{a}}_{p, \mathrm{OLS}}^\mathsf{T} \mat{S}_{p \times p} \widehat{\vect{a}}_{p, \mathrm{OLS}}}
\end{aligned}
\end{equation}
For the sample cross-covariance, eq.~\eqref{eq:cross_covariance_vector} implies:
\begin{equation}
\label{eq:ols_cross_covariance}
\begin{aligned}
\vect{s}_p &= (\mat{X}_{n \times p, \mathrm{c}}^\mathsf{T} \vect{y}_{n \mathrm{c}})/n = (1/n)\mat{X}_{n \times p, \mathrm{c}}^\mathsf{T} \mat{X}_{n \times p, \mathrm{c}}(\mat{X}_{n \times p, \mathrm{c}}^\mathsf{T} \mat{X}_{n \times p, \mathrm{c}})^{-1} \mat{X}_{n \times p, \mathrm{c}}^\mathsf{T} \vect{y}_{n \mathrm{c}} \\
&= \mat{S}_{p \times p} \widehat{\vect{a}}_{p, \mathrm{OLS}}
\end{aligned}
\end{equation}
We define the OLS correlation $\rho_\mathrm{OLS}$ as:
\begin{equation}
\label{eq:ols_rho_def}
\rho_\mathrm{OLS} := \rho(\vect{z}_{n, \mathrm{OLS}}, \vect{y}_n), \vect{s}_p \neq \0_p
\end{equation}
From eq.~\eqref{eq:ols_cross_covariance} the condition $\vect{s}_p = \0_p$ is equivalent to $\widehat{\vect{a}}_{p, \mathrm{OLS}} = \0_p$. In this zero cross-covariance case, $\widehat{b}_{p, \mathrm{OLS}} = \mu(\vect{y}_n)$, the OLS model predictions reduce to the constant vector $\vect{z}_{n, \mathrm{OLS}} = \mu(\vect{y}_n) \1_n$ and $\rho_\mathrm{OLS}$ becomes undefined; thus the domain of $\rho_\mathrm{OLS}$ is restricted to $\vect{s}_p \neq \0_p$. Using the centered vector expressions from eqs.~\eqref{eq:correlation_predictions_multi}, \eqref{eq:ols_centered_predictions}, \eqref{eq:ols_cross_covariance} and \eqref{eq:pearson}, we establish the following relationship for the correlation:
\begin{equation}
\label{eq:ols_correlation_deriv}
\begin{aligned}
\rho_\mathrm{OLS} &= \frac{\vect{z}_{n \mathrm{c}, \mathrm{OLS}}^\mathsf{T} \vect{y}_{n \mathrm{c}}}{n \sigma(\vect{z}_{n, \mathrm{OLS}}) \sigma(\vect{y}_n)} = \frac{\widehat{\vect{a}}_{p, \mathrm{OLS}}^\mathsf{T} \mat{X}_{n \times p,\mathrm{c}}^\mathsf{T} \vect{y}_{n \mathrm{c}}}{n \sigma(\vect{z}_{n, \mathrm{OLS}}) \sigma(\vect{y}_n)} = \frac{\widehat{\vect{a}}_{p, \mathrm{OLS}}^\mathsf{T} \vect{s}_p}{\sigma_\mathrm{OLS} \sigma(\vect{y}_n)}\\
&= \frac{\widehat{\vect{a}}_{p, \mathrm{OLS}}^\mathsf{T} \mat{S}_{p \times p} \widehat{\vect{a}}_{p, \mathrm{OLS}}}{\sigma_\mathrm{OLS} \sigma(\vect{y}_n)} = \frac{\sigma_\mathrm{OLS}^2}{\sigma_\mathrm{OLS} \sigma(\vect{y}_n)} = \frac{\sigma_\mathrm{OLS}}{\sigma(\vect{y}_n)} > 0, \vect{s}_p \neq \0_p
\end{aligned}
\end{equation}
The correlation is positive, within its domain, $\widehat{\vect{a}}_{p, \mathrm{OLS}} \neq \0_p$, and therefore $\sigma_\mathrm{OLS} > 0$ by eq.~\eqref{eq:ols_sigma_expansion}. Eq.~\eqref{eq:ols_correlation_deriv} immediately gives:
\begin{equation}
\label{eq:ols_variance_reduction}
\sigma_\mathrm{OLS} \leq \sigma(\vect{y}_n)
\end{equation}
In summary, the key statistical properties of the OLS predictions on the training set are:
\begin{equation}
\mu(\vect{z}_{n, \mathrm{OLS}}) = \mu(\vect{y}_n)
\end{equation}
\begin{equation}
\sigma(\vect{z}_{n, \mathrm{OLS}}) = \sigma_\mathrm{OLS} = \sqrt{\frac{1}{n} \vect{y}_n^\mathsf{T} \mat{P}_{n \times n} \mat{X}_{n \times p} (\mat{X}_{n \times p}^\mathsf{T} \mat{P}_{n \times n} \mat{X}_{n \times p})^{-1} \mat{X}_{n \times p}^\mathsf{T} \mat{P}_{n \times n} \vect{y}_n}
\end{equation}
\begin{equation}
\rho(\vect{z}_{n, \mathrm{OLS}}, \vect{y}_n) = \rho_\mathrm{OLS} = \frac{\sigma_\mathrm{OLS}}{\sigma(\vect{y}_n)}, \vect{s}_p \neq \0_p
\end{equation}
\subsubsection{Ordinary least squares linear regression with a single predictor}
\label{sec:ols-linear-regression-single-predictor}
We define the OLS estimator of the linear model in eq.~\eqref{eq:single_pred_model} as the parameter vector that minimizes the $\mathrm{MSE}$:
\begin{equation}
\widehat{\vect{\theta}}_{1 \times 2, \mathrm{OLS}}(\vect{x}_n, \vect{y}_n) := \underset{(a, b) \in \R^{1 \times 2}}{\operatorname*{arg\,min}} \mathrm{MSE}(a \vect{x}_n + b \1_n, \vect{y}_n)
\end{equation}
This formulation is a special case of the model estimated in Section~\ref{sec:ols-linear-regression}, corresponding to $p = 1$. Therefore, the equations that follow emerge directly from their counterparts in Section~\ref{sec:ols-linear-regression} with $p = 1$. The predictions of the OLS linear regression model on the training set take the form:
\begin{equation}
\label{eq:ols_predictions_single}
\vect{z}_{n, \mathrm{OLS}} = \widehat{a}_{\mathrm{OLS}} \vect{x}_n + \widehat{b}_{\mathrm{OLS}} \1_n
\end{equation}
Given a test predictor vector $\vect{x}_k$, the predictions from the estimated OLS linear regression model on the test set are given by:
\begin{equation}
\label{eq:ols_test_predictions_single}
\vect{z}_{k, \mathrm{OLS}, \mathrm{test}} = \widehat{a}_{\mathrm{OLS}} \vect{x}_k + \widehat{b}_{\mathrm{OLS}} \1_k
\end{equation}
The estimate of the slope parameter is:
\begin{equation}
\label{eq:ols_slope_single}
\widehat{a}_\mathrm{OLS} = \rho(\vect{x}_n, \vect{y}_n) \frac{\sigma(\vect{y}_n)}{\sigma(\vect{x}_n)}
\end{equation}
The estimate of the intercept parameter is:
\begin{equation}
\label{eq:ols_intercept_single}
\widehat{b}_\mathrm{OLS} = \mu(\vect{y}_n) - \mu(\vect{x}_n) \widehat{a}_\mathrm{OLS}
\end{equation}
In summary, the key statistical properties of the OLS predictions on the training set simplify to:
\begin{equation}
\mu(\vect{z}_{n, \mathrm{OLS}}) = \mu(\vect{y}_n)
\end{equation}
\begin{equation}
\sigma(\vect{z}_{n, \mathrm{OLS}}) = \sigma_\mathrm{OLS} = |\rho(\vect{x}_n, \vect{y}_n)| \sigma(\vect{y}_n)
\end{equation}
\begin{equation}
\label{eq:single_predictor_predict_obs_cor}
\rho(\vect{z}_{n, \mathrm{OLS}}, \vect{y}_n) = \rho_\mathrm{OLS} = |\rho(\vect{x}_n, \vect{y}_n)|, \rho(\vect{x}_n, \vect{y}_n) \neq 0
\end{equation}
\section{Kling-Gupta regression}
\label{sec:kling_gupta_regressions}
This section details the primary theoretical developments of the manuscript. Building on the Kling-Gupta loss established in eq.~\eqref{eq:Kling_Gupta_loss}, we specify the extremum estimation framework for linear models (Section~\ref{sec:linear-model}) and identify closed-form expressions for the associated parameter estimates. We define this estimation procedure as \textit{Kling-Gupta linear regression} (Section~\ref{sec:kling-gupta-linear-regression}). The estimator functions as a variance-inflated variant of OLS, ensuring predictions replicate the sample variance of the observations while maintaining mean and correlation properties. Section~\ref{sec:asymptotic_properties} establishes the asymptotic limits of the estimators, while Section~\ref{sec:comparative_training_data} quantifies the trade-offs between $\mathrm{MSE}$, $\mathrm{NSE}$, $L_\mathrm{KG}$, and $\mathrm{KGE}$ within the training set. These findings extend to infinite-sample settings in Sections~\ref{sec:comparative_infinite_training_data} and~\ref{sec:comparative_infinite_training_test_data}, demonstrating that identical performance relationships persist asymptotically on independent test sets. Finally, we investigate constrained estimation: Section~\ref{sec:fixed_regression_coefficient} analyzes fixed coefficient vectors, whereas Section~\ref{sec:fixed_intercept_kg_linear_regression} investigates the single-predictor model with a fixed intercept, characterizing conditions under which a global minimum for the slope parameter does not exist.
\subsection{Kling-Gupta linear regression}
\label{sec:kling-gupta-linear-regression}
Section~\ref{sec:kling_gupta_loss_sec} established the Kling-Gupta loss as an extremum estimator; we now apply this framework to the linear models from Section~\ref{sec:linear-model}. A comparison between these estimates and those from the OLS linear regression in Section~\ref{sec:estimating-linear-model-se} clarifies how minimizing the Kling-Gupta loss acts on the linear model during estimation. This approach defines the Kling-Gupta linear regression introduced in Section~\ref{sec:introduction}.
\subsubsection{The Kling-Gupta linear regression with multiple predictors}
\label{sec:kling_gupta_linear_regression_multiple_predictors}
The analytical results that follow are established in \hyperref[proof:b2]{Proof~B.2}. We determine the parameter estimates for the linear model $z = \vect{a}_p^\mathsf{T} \vect{x}_p + b_p$ specified in eq.~\eqref{eq:linear_model} by minimizing the Kling-Gupta loss function $L_\mathrm{KG}$ defined in eq.~\eqref{eq:Kling_Gupta_loss}, given observations $\vect{y}_n$ of the response variable and $\mat{X}_{n \times p}$ of the predictor variables, while following the notational conventions of Section~\ref{sec:linear_model_multiple_predictors}. The sample size is assumed to be sufficiently large such that $n \geq p$. We further assume that $\mat{X}_{n \times p}$ has full column rank and that the all-ones vector $\1_n$ does not lie in its column space. Kling-Gupta linear regression minimizes $L_\mathrm{KG}$ over the parameter vector $\vect{\theta}_{1 \times (p + 1)}$ defined in eq.~\eqref{eq:multi_linear_model_parameter}:
\begin{equation}
\label{eq:kg_estimator_multi_def}
\widehat{\vect{\theta}}_{1 \times (p + 1), \mathrm{KG}}(\mat{X}_{n \times p}, \vect{y}_n) := \underset{\vect{\theta}_{1 \times (p + 1)} \in \vect{\Theta}}{\operatorname*{arg\,min}} L_\mathrm{KG}(\widetilde{\mat{X}}_{n \times (p + 1)} \vect{\theta}_{1 \times (p + 1)}^\mathsf{T}, \vect{y}_n)
\end{equation}
Prior to analyzing the parameter estimates, we contextualize them relative to the OLS estimates discussed in Section~\ref{sec:ols-linear-regression}. Recall that OLS identifies the unique linear model that minimizes the $\mathrm{MSE}$, with the coefficient vector $\widehat{\vect{a}}_{p, \mathrm{OLS}}$ and intercept $\widehat{b}_{p, \mathrm{OLS}}$ estimates given by eqs.~\eqref{eq:ols_slope_multi} and \eqref{eq:ols_intercept_multi}, respectively. Importantly, as established in eq.~\eqref{eq:ols_variance_reduction}, OLS predictions on the training set exhibit a sample standard deviation $\sigma_\mathrm{OLS}$ that is less than or equal to that of the observations, $\sigma(\vect{y}_n)$, as a result of minimizing squared vertical distances. The OLS intercept, however, perfectly aligns the sample means, ensuring $\mu(\vect{z}_{n, \mathrm{OLS}}) = \mu(\vect{y}_n)$.

In cases where the sample cross-covariance vector is non-zero ($\vect{s}_p \neq \0_p$), the Kling-Gupta coefficient parameter estimate is an explicit scaling of the OLS coefficients:
\begin{equation}
\label{eq:kg_slope_multi}
\widehat{\vect{a}}_{p, \mathrm{KG}} = \frac{\sigma(\vect{y}_n)}{\sigma_\mathrm{OLS}} \widehat{\vect{a}}_{p, \mathrm{OLS}}
\end{equation}
Eq.~\eqref{eq:kg_slope_multi} implies that the Kling-Gupta extremum estimator maintains the direction of the regression vector determined by OLS. It amplifies the magnitude of this vector by a factor $\frac{\sigma(\vect{y}_n)}{\sigma_\mathrm{OLS}} \geq 1$, which is the specific multiplier that eliminates the variance reduction inherent to OLS (see eq.~\eqref{eq:ols_variance_reduction}). Accordingly, Kling-Gupta linear regression can be viewed as OLS linear regression followed by a subsequent variance-inflation correction.

When the cross-covariance vanishes ($\vect{s}_p = \0_p$), the OLS estimator produces a constant prediction equal to $\mu(\vect{y}_n)$. The Kling-Gupta estimator, however, cannot admit a zero-variance prediction due to the domain restriction in eq.~\eqref{eq:Kling_Gupta_loss}. The minimization problem therefore does not result in a unique point estimate; rather, the solution lies on the ellipsoid defined by the predictor covariance structure:
\begin{equation}
\label{eq:kg_ellipsoid_condition}
\widehat{\vect{a}}_{p, \mathrm{KG}}^\mathsf{T} \mat{S}_{p \times p} \widehat{\vect{a}}_{p, \mathrm{KG}} = \sigma^2(\vect{y}_n)
\end{equation}
Any coefficient vector satisfying eq.~\eqref{eq:kg_ellipsoid_condition} forces the estimated Kling-Gupta linear model to issue predictions whose sample variance equals that of the observations, thereby setting the variance penalty term of $L_\mathrm{KG}$ to zero. This behavior highlights a characteristic of optimizing the Kling-Gupta loss when the correlation between predictions and observations is zero: the estimator favors a model that generates high-variance fluctuations over the constant mean climatology.

Despite the pronounced differences in slope estimation, the mechanism for enforcing mean-unbiasedness remains structurally identical to that of OLS. In all cases, the Kling-Gupta linear regression model determines its intercept estimate by aligning the sample means:
\begin{equation}
\label{eq:kg_intercept_multi}
\widehat{b}_{p, \mathrm{KG}} = \mu(\vect{y}_n) - (\vect{\mu}(\mat{X}_{n \times p}))^\mathsf{T} \widehat{\vect{a}}_{p, \mathrm{KG}}
\end{equation}
Similar to the OLS estimate (eq.~\eqref{eq:ols_intercept_multi}), the intercept functions as a centering term that eliminates any bias introduced by the slope coefficients.

To facilitate a unified treatment of the correlation between Kling-Gupta predictions and observations independent of the value of $\vect{s}_p$, we establish the following convention:
\begin{equation}
\label{eq:rho_ols_star_def}
\rho_{\mathrm{OLS}_{(0)}} =
\begin{cases}
\rho_\mathrm{OLS}, & \text{if } \vect{s}_p \neq \0_p \\
0, & \text{if } \vect{s}_p = \0_p
\end{cases}
\end{equation}
Within the training set, the predictions of the Kling-Gupta estimated model are expressed as:
\begin{equation}
\label{eq:kg_predictions_multi}
\vect{z}_{n, \mathrm{KG}} = \mat{X}_{n \times p} \widehat{\vect{a}}_{p, \mathrm{KG}} + \widehat{b}_{p, \mathrm{KG}} \1_n
\end{equation}
Given a test predictor matrix $\mat{X}_{k \times p, \mathrm{test}}$, the predictions from the estimated Kling-Gupta linear regression model on the test set are given by:
\begin{equation}
\label{eq:kg_test_predictions_multi}
\vect{z}_{k, \mathrm{KG}, \mathrm{test}} = \mat{X}_{k \times p, \mathrm{test}} \widehat{\vect{a}}_{p, \mathrm{KG}} + \widehat{b}_{p, \mathrm{KG}} \1_k
\end{equation}
The statistical properties of the model's predictions on the training set are characterized as follows:
\begin{equation}
\label{eq:KG_mu_predictions}
\mu(\vect{z}_{n, \mathrm{KG}}) = \mu(\vect{y}_n)
\end{equation}
\begin{equation}
\label{eq:KG_sigma_predictions}
\sigma(\vect{z}_{n, \mathrm{KG}}) = \sigma(\vect{y}_n)
\end{equation}
\begin{equation}
\label{eq:KG_rho_predictions}
\rho(\vect{z}_{n, \mathrm{KG}}, \vect{y}_n) = \rho_{\mathrm{OLS}_{(0)}}
\end{equation}
These equations formalize the trade-off of Kling-Gupta linear regression. The estimator matches the observed variance at the expense of the least-squares optimality of OLS, which involves predictive variance reduction. In both models, the predictive means and correlations remain identical, and the predictive mean equals the observed mean.
\subsubsection{The Kling-Gupta linear regression with a single predictor}
\label{sec:kling_gupta_linear_regression_single_predictor}
The Kling-Gupta extremum estimator for the linear model $z = a x + b$ in eq.~\eqref{eq:single_pred_model} is defined as the parameter vector that minimizes $L_\mathrm{KG}$:
\begin{equation}
\widehat{\vect{\theta}}_{1 \times 2, \mathrm{KG}}(\vect{x}_n, \vect{y}_n) := \underset{(a, b) \in \R^{1 \times 2}}{\operatorname*{arg\,min}} L_\mathrm{KG}(a \vect{x}_n + b \1_n, \vect{y}_n)
\end{equation}
The single-predictor case follows immediately from the general multiple-predictor framework of Section~\ref{sec:kling_gupta_linear_regression_multiple_predictors} by setting $p = 1$. While the results are therefore special cases of those already established, examining the simpler model offers substantial interpretive value. In instances where $\rho(\vect{x}_n, \vect{y}_n) \neq 0$, the coefficient parameter estimate is formulated as:
\begin{equation}
\label{eq:kg_slope_single}
\widehat{a}_\mathrm{KG} = \frac{\sigma(\vect{y}_n)}{\sigma_\mathrm{OLS}} \widehat{a}_{\mathrm{OLS}} = \operatorname{sign}(\rho(\vect{x}_n, \vect{y}_n)) \frac{\sigma(\vect{y}_n)}{\sigma(\vect{x}_n)}
\end{equation}
Therefore, the Kling-Gupta slope inherits the sign of the OLS slope but amplifies its magnitude by the factor $\frac{\sigma(\vect{y}_n)}{\sigma_\mathrm{OLS}}$.

When the predictor and response are uncorrelated ($\rho(\vect{x}_n, \vect{y}_n) = 0$), the OLS estimator issues a constant prediction equal to $\mu(\vect{y}_n)$. The Kling-Gupta extremum estimator, however, cannot issue a constant prediction without violating the domain condition $\sigma(\vect{z}_n) \neq 0$. The optimization therefore determines two symmetric estimates:
\begin{equation}
\widehat{a}_\mathrm{KG} = \pm \frac{\sigma(\vect{y}_n)}{\sigma(\vect{x}_n)}
\end{equation}
This sign ambiguity arises from the structure of the Kling-Gupta loss. The correlation term is determined by $\operatorname{sign}(a) \rho(\vect{x}_n, \vect{y}_n)$, while the variability term depends only on $|a|$. When $\rho(\vect{x}_n, \vect{y}_n) = 0$, the correlation term simplifies to unity regardless of the sign of $a$, and the loss function treats positive and negative slopes as equally acceptable, as long as they align with the observed variance.

In all cases, the intercept parameter estimate is determined to align the sample means, consistent with the OLS estimate:
\begin{equation}
\label{eq:kg_intercept_single}
\widehat{b}_\mathrm{KG} = \mu(\vect{y}_n) - \mu(\vect{x}_n) \widehat{a}_\mathrm{KG}
\end{equation}
Within the training set, the predictions of the Kling-Gupta estimated model are expressed as:
\begin{equation}
\label{eq:kg_predictions_single}
\vect{z}_{n, \mathrm{KG}} = \widehat{a}_\mathrm{KG} \vect{x}_n + \widehat{b}_\mathrm{KG} \1_n
\end{equation}
Given a test predictor vector $\vect{x}_k$, the predictions from the estimated Kling-Gupta linear regression model on the test set are given by:
\begin{equation}
\label{eq:kg_test_predictions_single}
\vect{z}_{k, \mathrm{KG}, \mathrm{test}} = \widehat{a}_{\mathrm{KG}} \vect{x}_k + \widehat{b}_{\mathrm{KG}} \1_k
\end{equation}
The statistical properties of the model's predictions on the training set are characterized as follows:
\begin{equation}
\mu(\vect{z}_{n, \mathrm{KG}}) = \mu(\vect{y}_n)
\end{equation}
\begin{equation}
\sigma(\vect{z}_{n, \mathrm{KG}}) = \sigma(\vect{y}_n)
\end{equation}
\begin{equation}
\label{eq:rho_single_ols_predictions}
\rho(\vect{z}_{n, \mathrm{KG}}, \vect{y}_n) = |\rho(\vect{x}_n, \vect{y}_n)|
\end{equation}
The single-predictor case constitutes an intuitive framework for understanding the more general multiple-predictor estimator. Kling-Gupta linear regression acts on the OLS estimate by applying a scale correction that inflates the predictions until their variance equals that of the observations, while the correlation remains unchanged.
\subsubsection{Numerical illustration}
\label{sec:numerical_illustration}
To numerically and visually illustrate the theoretical properties and performance trade-offs of the OLS and Kling-Gupta estimators, we perform a simulation experiment. We generate $n = 10^6$ training samples and $k = 10^6$ test samples, where the predictor variable $\underline{x} \sim N(1, 1)$ (with probability density function defined in eq.~\eqref{eq:normal_pdf}) and the error term $\underline{\epsilon} \sim N(0, 2)$. The response variable $\underline{y}$ is constructed via the linear relation $\underline{y} = a_0 \underline{x} + b_0 + \underline{\epsilon}$, with true parameter values $a_0 = 0.60$ and $b_0 = 1.00$. Under this data-generating setup, the population characteristics of the response variable are $\mathbb{E}_{F_{\underline{y}}}[\underline{y}] = 1.60$, $\mathrm{Var}_{F_{\underline{y}}}[\underline{y}] = 4.36$, and $\mathrm{Corr}_{F_{\underline{x}, \underline{y}}}(\underline{x}, \underline{y}) \approx 0.2873$. We fit the single-predictor linear model $z = a x + b$ (defined in eq.~\eqref{eq:single_pred_model}) to the training set. The parameters $a$ and $b$ are estimated from the training data using two methods: OLS linear regression (with estimates given by eqs.~\eqref{eq:ols_slope_single} and \eqref{eq:ols_intercept_single}) and Kling-Gupta linear regression (with estimates given by eqs.~\eqref{eq:kg_slope_single} and \eqref{eq:kg_intercept_single}).

\begin{figure}[htbp]
\centering
\includegraphics[scale=1.00]{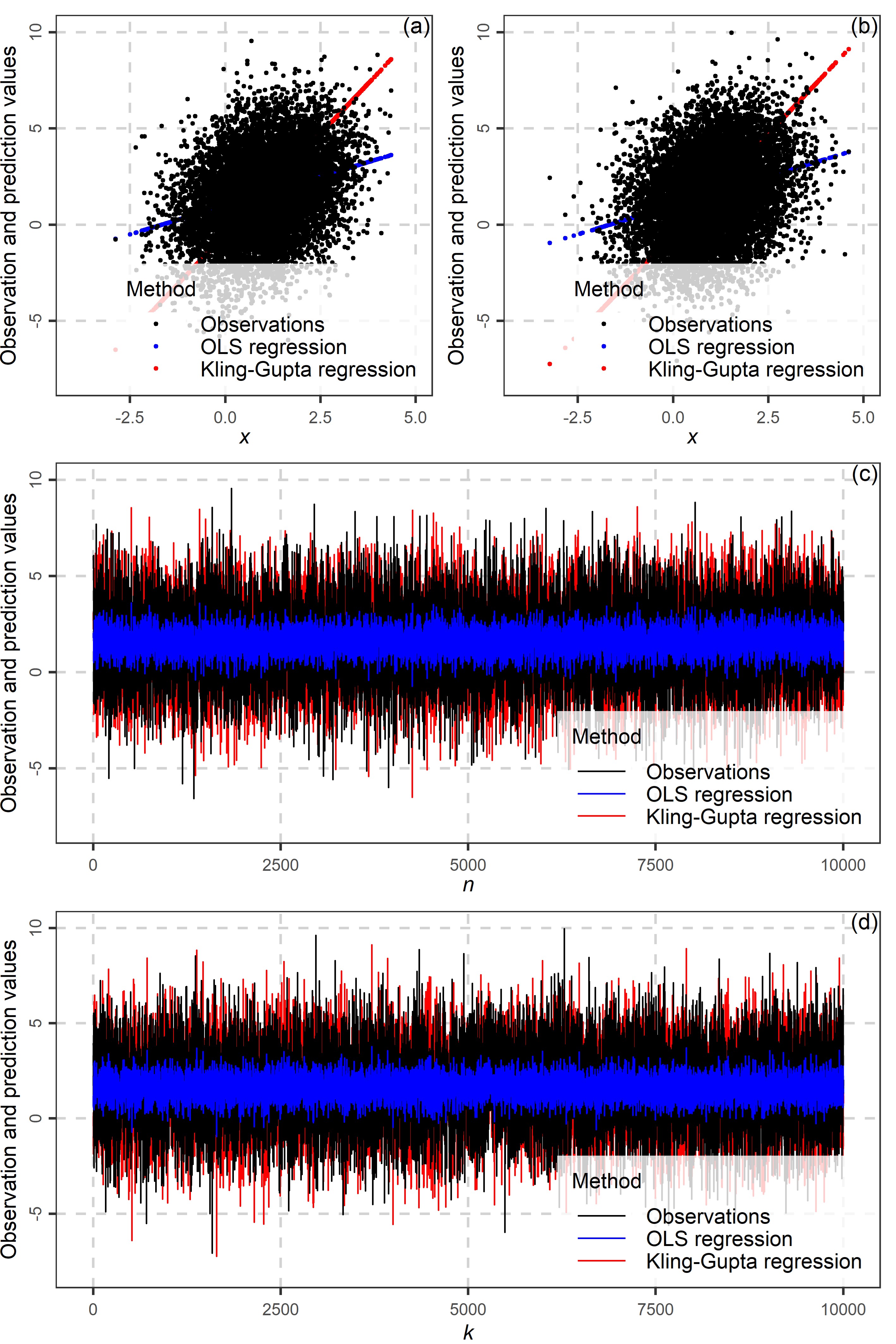}
\caption{Visual comparison of OLS (in blue) and Kling-Gupta linear regression (in red) for the simulation data in Section~\ref{sec:numerical_illustration}: (a) scatterplot of observations and predictions against the predictor $x$ for the first $10^5$ training samples, (b) scatterplot of observations and predictions against $x$ for the first $10^5$ test samples, (c) time series of observations and predictions for the first $10^5$ training samples, (d) time series of observations and predictions for the first $10^5$ test samples.}
\label{fig:figure1}
\end{figure}

Figure~\ref{fig:figure1} visually compares the predictions from both estimated models across subsets of the training set and the test set. OLS predictions are computed using eqs.~\eqref{eq:ols_predictions_single} (training) and \eqref{eq:ols_test_predictions_single} (test), while Kling-Gupta predictions are computed using eqs.~\eqref{eq:kg_predictions_single} (training) and \eqref{eq:kg_test_predictions_single} (test). For visual clarity, predictions on the complete dataset are not shown. Table~\ref{tab:numerical_results} summarizes the numerical results of the simulation, including sample means, variances, and correlations evaluated on the training and test sets, together with their population counterparts.
\begin{table}[htbp]
\centering
\setlength{\tabcolsep}{3pt}
\caption{Statistical properties of observations, OLS predictions, and Kling-Gupta predictions for the simulation example of Section~\ref{sec:numerical_illustration}. Values are reported for the population, the training set, and the test set.}
\label{tab:numerical_results}
\begin{tabular}{lrrr}
\toprule
Statistic & Population & Training set & Test set \\
\midrule
Mean of observations & 1.6000 & 1.6028 & 1.6028 \\
Mean of OLS predictions & 1.6000 & 1.6028 & 1.6018 \\
Mean of Kling-Gupta predictions & 1.6000 & 1.6028 & 1.5993 \\
\midrule
Variance of observations & 4.3600 & 4.3591 & 4.3682 \\
Variance of OLS predictions & 0.3600 & 0.3619 & 0.3608 \\
Variance of Kling-Gupta predictions & 4.3600 & 4.3591 & 4.3460 \\
\midrule
Correlation (OLS predictions, observations) & 0.2873 & 0.2881 & 0.2874 \\
Correlation (Kling-Gupta predictions, observations) & 0.2873 & 0.2881 & 0.2874 \\
\bottomrule
\end{tabular}
\end{table}

Panels (a) and (b) in Figure~\ref{fig:figure1} show scatterplots of the observations together with the estimated OLS and Kling-Gupta regression lines plotted against the predictor $x$ for the training and test sets, respectively. These plots illustrate how the variance of the predictions depends on the estimation method. The OLS regression exhibits a flatter slope, so the variability of the predictions is reduced, which is consistent with minimizing the $\mathrm{MSE}$. In contrast, the Kling-Gupta regression line exhibits a much steeper slope. This demonstrates that the estimator explicitly inflates the regression coefficients by the scaling factor $\sigma(\vect{y}_n)/\sigma_{\mathrm{OLS}}$ to ensure that the predictions reproduce the full marginal sample variance of the observations.

Panels (c) and (d) of Figure~\ref{fig:figure1} show the sequentially ordered observations together with the predictions from both models for the training and test periods, resembling a hydrograph comparison. These time series plots indicate that the sample means of the observations and the OLS predictions are approximately equal, whereas the variance of the OLS predictions is substantially lower than that of the observations, as expected from eq.~\eqref{eq:ols_correlation_deriv} for the training set, which implies $\sigma_\mathrm{OLS} \leq \sigma(\vect{y}_n)$. In contrast, the Kling-Gupta predictions exhibit a mean that aligns with the observations, while their variance also matches that of the observations, as expected from eq.~\eqref{eq:KG_sigma_predictions} for the training set. These visual findings are confirmed by the numerical values in Table~\ref{tab:numerical_results} for both the training and test sets. As demonstrated in Sections~\ref{sec:ols-linear-regression} and \ref{sec:kling-gupta-linear-regression}, such results are expected for the training set. For the test set, the results agree with the asymptotic behavior as $n \to \infty$ and $k \to \infty$, as shown in subsequent sections. Furthermore, the correlations between predictions and observations are equal for the OLS and Kling-Gupta models, as expected from eqs.~\eqref{eq:single_predictor_predict_obs_cor} and \eqref{eq:rho_single_ols_predictions}.
\subsection{Asymptotic properties of parameter estimators}
\label{sec:asymptotic_properties}
The estimators established in Section~\ref{sec:kling-gupta-linear-regression} give explicit expressions for the Kling-Gupta linear regression parameters. Although these formulas are analytic for any finite sample size $n$, they offer limited insight into the behavior of the estimator as the sample size increases. Characterizing the asymptotic properties, particularly the limits toward which the estimators converge as $n \to \infty$, is essential to identify the population functional targeted by Kling-Gupta regression, to assess its statistical consistency, and to facilitate large-sample comparisons with OLS. Therefore, we determine the probability limits of the OLS and Kling-Gupta coefficient estimators under standard regularity conditions, assuming the data consist of independent and identically distributed (i.i.d.) samples from a joint distribution $F_{\underline{\vect{x}}_p, \underline{y}}$. These developments leverage the convergence properties of general statistics established in eqs.~\eqref{eq:sample_mean_convergence}, \eqref{eq:sample_standard_deviation_convergence}, \eqref{eq:sample_Pearson_correlation_convergence}, \eqref{eq:matrix_sample_mean_convergence}, \eqref{eq:sample_covariance_matrix_convergence}, and \eqref{eq:sample_cross_covariance_vector_convergence}.

To separate realized estimates from estimators, which are random variables, we underline the notation for estimators to signify their nature as random variables. The sample quantity $\sigma_\mathrm{OLS}$ from eq.~\eqref{eq:ols_sigma_expansion} converges almost surely to:
\begin{equation}
\label{eq:sigma_OLS_asymptotic}
\underline{\sigma}_\mathrm{OLS} \overset{\mathrm{a.s.}}{\longrightarrow} \sigma_{\mathrm{OLS}^*} := \sqrt{\mathrm{Cov}_{F_{\underline{\vect{x}}_p, \underline{y}}}(\underline{\vect{x}}_p, \underline{y})^\mathsf{T} (\mathrm{Var}_{F_{\underline{\vect{x}}_p}}[\underline{\vect{x}}_p])^{-1} \mathrm{Cov}_{F_{\underline{\vect{x}}_p, \underline{y}}}(\underline{\vect{x}}_p, \underline{y})}
\end{equation}
In this expression, $\mathrm{Cov}_{F_{\underline{\vect{x}}_p, \underline{y}}}(\underline{\vect{x}}_p, \underline{y})$ is specified in eq.~\eqref{eq:covariances_vector}, whereas $\mathrm{Var}_{F_{\underline{\vect{x}}_p}}[\underline{\vect{x}}_p]$ is defined in eq.~\eqref{eq:variance_covariance_matrix}. For the Kling-Gupta linear regression estimator, we separate the case where the predictor-response covariance vector is non-zero ($\mathrm{Cov}_{F_{\underline{\vect{x}}_p, \underline{y}}}(\underline{\vect{x}}_p, \underline{y}) \neq \0_p$) from the zero cross-covariance scenario ($\mathrm{Cov}_{F_{\underline{\vect{x}}_p, \underline{y}}}(\underline{\vect{x}}_p, \underline{y}) = \0_p$). Under the i.i.d. assumption, the sample predictive correlation specified in eq.~\eqref{eq:ols_correlation_deriv} converges as $n \to \infty$:
\begin{equation}
\label{eq:rho_OLS_asymptotic}
\underline{\rho}_\mathrm{OLS} \overset{\mathrm{a.s.}}{\longrightarrow} \rho_{\mathrm{OLS}^{*}} = \frac{\sigma_{\mathrm{OLS}^*}}{\sqrt{\mathrm{Var}_{F_{\underline{y}}}[\underline{y}]}}, \mathrm{Cov}_{F_{\underline{\vect{x}}_p, \underline{y}}}(\underline{\vect{x}}_p, \underline{y}) \neq \0_p
\end{equation}
We generalize $\underline{\rho}_\mathrm{OLS}$ to the case where $\mathrm{Cov}_{F_{\underline{\vect{x}}_p, \underline{y}}}(\underline{\vect{x}}_p, \underline{y}) = \0_p$:
\begin{equation}
\label{eq:rho_OLS_0_asymptotic}
\underline{\rho}_{\mathrm{OLS}_{(0)}} \overset{\mathrm{a.s.}}{\longrightarrow} \rho_{\mathrm{OLS}_{(0)}^*} =
\begin{cases}
\rho_{\mathrm{OLS}^*}, & \text{if } \mathrm{Cov}_{F_{\underline{\vect{x}}_p, \underline{y}}}(\underline{\vect{x}}_p, \underline{y}) \neq \0_p \\
0, & \text{if } \mathrm{Cov}_{F_{\underline{\vect{x}}_p, \underline{y}}}(\underline{\vect{x}}_p, \underline{y}) = \0_p
\end{cases}
\end{equation}
along with the analogous single-predictor case, which follows from eq.~\eqref{eq:sample_Pearson_correlation_convergence}:
\begin{equation}
\rho(\underline{\vect{x}}_n, \underline{\vect{y}}_n) \overset{\mathrm{a.s.}}{\longrightarrow} \mathrm{Corr}_{F_{\underline{x}, \underline{y}}}(\underline{x}, \underline{y})
\end{equation}
\subsubsection{Linear models with multiple predictors}
\label{sec:linear_models_with_multiple_predictors}
We first summarize the asymptotic limits for the OLS estimator, which are documented in the literature and function as a benchmark. Applying the continuous mapping theorem to eqs.~\eqref{eq:ols_slope_multi} and \eqref{eq:ols_intercept_multi} implies:
\begin{equation}
\label{eq:OLS_asymptotic_coefficient_estimate}
\widehat{\underline{\vect{a}}}_{p, \mathrm{OLS}} \overset{\mathrm{a.s.}}{\longrightarrow} \vect{a}_{p, \mathrm{OLS}^*} := (\mathrm{Var}_{F_{\underline{\vect{x}}_p}}[\underline{\vect{x}}_p])^{-1} \mathrm{Cov}_{F_{\underline{\vect{x}}_p, \underline{y}}}(\underline{\vect{x}}_p, \underline{y})
\end{equation}
\begin{equation}
\label{eq:OLS_asymptotic_intercept_estimate}
\widehat{\underline{b}}_{p, \mathrm{OLS}} \overset{\mathrm{a.s.}}{\longrightarrow} b_{p, \mathrm{OLS}^*} := \mathbb{E}_{F_{\underline{y}}}[\underline{y}] - (\mathbb{E}_{F_{\underline{\vect{x}}_p}}[\underline{\vect{x}}_p])^\mathsf{T} \vect{a}_{p, \mathrm{OLS}^*}
\end{equation}
where the component-wise expectation $\mathbb{E}_{F_{\underline{\vect{x}}_p}}[\underline{\vect{x}}_p]$ is defined in eq.~\eqref{eq:component_wise_mean}.

In the non-zero cross-covariance case, eq.~\eqref{eq:kg_slope_multi} characterizes the Kling-Gupta estimator as a scaled version of the OLS coefficients. By applying the continuous mapping theorem, the asymptotic limit of the Kling-Gupta slope vector is therefore established as:
\begin{equation}
\label{eq:KG_asymptotic_coefficient_estimate}
\widehat{\underline{\vect{a}}}_{p, \mathrm{KG}} \overset{\mathrm{a.s.}}{\longrightarrow} \vect{a}_{p, \mathrm{KG}^*} := \frac{\sqrt{\mathrm{Var}_{F_{\underline{y}}}[\underline{y}]}}{\sigma_{\mathrm{OLS}^*}} \vect{a}_{p, \mathrm{OLS}^*}
\end{equation}
The intercept estimator established in eq.~\eqref{eq:kg_intercept_multi} converges analogously:
\begin{equation}
\label{eq:KG_asymptotic_intercept_estimate}
\widehat{\underline{b}}_{p, \mathrm{KG}} \overset{\mathrm{a.s.}}{\longrightarrow} b_{p, \mathrm{KG}^*} := \mathbb{E}_{F_{\underline{y}}}[\underline{y}] - (\mathbb{E}_{F_{\underline{\vect{x}}_p}}[\underline{\vect{x}}_p])^\mathsf{T} \vect{a}_{p, \mathrm{KG}^*}
\end{equation}
In the zero cross-covariance case ($\mathrm{Cov}_{F_{\underline{\vect{x}}_p, \underline{y}}}(\underline{\vect{x}}_p, \underline{y}) = \0_p$), the sample estimate $\widehat{\vect{a}}_{p, \mathrm{KG}}$ is not unique; any admissible estimator satisfies the ellipsoid condition of eq.~\eqref{eq:kg_ellipsoid_condition}. As $n \to \infty$, any accumulation point $\vect{a}_{p, \mathrm{KG}^*}$ of a sequence $\{\widehat{\vect{a}}_{p, \mathrm{KG}}\}$ must lie on the asymptotic ellipsoid:
\begin{equation}
\label{eq:a_KG_asymptotic_ellipsoid}
\vect{a}_{p, \mathrm{KG}^*}^\mathsf{T} \mathrm{Var}_{F_{\underline{\vect{x}}_p}}[\underline{\vect{x}}_p] \vect{a}_{p, \mathrm{KG}^*} = \mathrm{Var}_{F_{\underline{y}}}[\underline{y}]
\end{equation}
Therefore, the corresponding intercept estimates accumulate according to:
\begin{equation}
\label{eq:b_KG_asymptotic_intercept_estimate}
b_{p, \mathrm{KG}^*} = \mathbb{E}_{F_{\underline{y}}}[\underline{y}] - (\mathbb{E}_{F_{\underline{\vect{x}}_p}}[\underline{\vect{x}}_p])^\mathsf{T} \vect{a}_{p, \mathrm{KG}^*}
\end{equation}
\subsubsection{Linear models with a single predictor}
\label{sec:linear_models_with_single_predictor}
In the single-predictor case ($p = 1$), the limits for the OLS parameter estimators reduce to familiar expressions. According to eq.~\eqref{eq:ols_slope_single}:
\begin{equation}
\widehat{\underline{a}}_\mathrm{OLS} \overset{\mathrm{a.s.}}{\longrightarrow} a_{\mathrm{OLS}^*} := (\mathrm{Var}_{F_{\underline{x}}}[\underline{x}])^{-1} \mathrm{Cov}_{F_{\underline{x}, \underline{y}}}(\underline{x}, \underline{y})
\end{equation}
in which $\mathrm{Cov}_{F_{\underline{x}, \underline{y}}}(\underline{x}, \underline{y})$ is defined in eq.~\eqref{eq:covariance}. Eq.~\eqref{eq:ols_intercept_single} then implies:
\begin{equation}
\widehat{\underline{b}}_\mathrm{OLS} \overset{\mathrm{a.s.}}{\longrightarrow} b_{\mathrm{OLS}^*} := \mathbb{E}_{F_{\underline{y}}}[\underline{y}] - \mathbb{E}_{F_{\underline{x}}}[\underline{x}] a_{\mathrm{OLS}^*}
\end{equation}
The limits for the Kling-Gupta linear regression parameter estimators simplify, based on eq.~\eqref{eq:kg_slope_single}, to:
\begin{equation}
\widehat{\underline{a}}_\mathrm{KG} \overset{\mathrm{a.s.}}{\longrightarrow} a_{\mathrm{KG}^*} := \operatorname{sign}(\mathrm{Cov}_{F_{\underline{x}, \underline{y}}}(\underline{x}, \underline{y}))\sqrt{\frac{\mathrm{Var}_{F_{\underline{y}}}[\underline{y}]}{\mathrm{Var}_{F_{\underline{x}}}[\underline{x}]}}
\end{equation}
Furthermore, eq.~\eqref{eq:kg_intercept_single} implies:
\begin{equation}
\widehat{\underline{b}}_\mathrm{KG} \overset{\mathrm{a.s.}}{\longrightarrow} b_{\mathrm{KG}^*} := \mathbb{E}_{F_{\underline{y}}}[\underline{y}] - \mathbb{E}_{F_{\underline{x}}}[\underline{x}] a_{\mathrm{KG}^*}
\end{equation}
In the zero cross-covariance case ($\mathrm{Cov}_{F_{\underline{x}, \underline{y}}}(\underline{x}, \underline{y}) = 0$), the estimator $\widehat{a}_\mathrm{KG}$ is not unique; any admissible estimator satisfies $|\widehat{a}_\mathrm{KG}| = \sigma(\vect{y}_n)/\sigma(\vect{x}_n)$. As $n \to \infty$, any accumulation point $a_{\mathrm{KG}^*}$ of a sequence $\{\widehat{a}_\mathrm{KG}\}$ must satisfy:
\begin{equation}
a_{\mathrm{KG}^*} = \pm \sqrt{\frac{\mathrm{Var}_{F_{\underline{y}}}[\underline{y}]}{\mathrm{Var}_{F_{\underline{x}}}[\underline{x}]}}.
\end{equation}
The corresponding intercepts accumulate at:
\begin{equation}
b_{\mathrm{KG}^*} := \mathbb{E}_{F_{\underline{y}}}[\underline{y}] - \mathbb{E}_{F_{\underline{x}}}[\underline{x}] a_{\mathrm{KG}^*}
\end{equation}
\subsection{Comparative performance on training data}
\label{sec:comparative_training_data}
The parameterizations established in Sections~\ref{sec:kling_gupta_linear_regression_multiple_predictors} and~\ref{sec:kling_gupta_linear_regression_single_predictor} characterize distinct training set performance properties. Table~\ref{tab:training-performance} compares the OLS and Kling-Gupta linear regressions with respect to $\mathrm{MSE}$, $\mathrm{NSE}$, $L_\mathrm{KG}$, and $\mathrm{KGE}$; the corresponding proofs are detailed in \hyperref[proof:b3]{Proof~B.3}. OLS regression minimizes $\mathrm{MSE}$, which implies a quadratic relationship between skill and correlation ($\mathrm{NSE} = \rho_{\mathrm{OLS}_{(0)}}^2$). Conversely, Kling-Gupta regression enforces mean and variance matching ($\mu(\vect{z}_{n, \mathrm{KG}}) = \mu(\vect{y}_n)$ and $\sigma(\vect{z}_{n, \mathrm{KG}}) = \sigma(\vect{y}_n)$). Therefore, the $\mathrm{KGE}$ of a Kling-Gupta model simplifies to $\rho_{\mathrm{OLS}_{(0)}}$, which relates predictive skill linearly to the predictor correlation.
\begin{table}[ht]
\centering
\caption{Performance of OLS and Kling-Gupta linear regression on the training set, reporting the $\mathrm{MSE}$, $\mathrm{NSE}$, $L_\mathrm{KG}$, and $\mathrm{KGE}$. The upper panel addresses the general multiple-predictor linear model, whereas the lower panel treats the single-predictor special case. All entries are formulated using observable sample statistics: the response standard deviation $\sigma(\vect{y}_n)$, the OLS correlation $\rho_{\mathrm{OLS}_{(0)}}$ (eq.~\eqref{eq:rho_ols_star_def}), and, for the single-predictor case, the predictor-response correlation $\rho(\vect{x}_n, \vect{y}_n)$. For the OLS model, $L_\mathrm{KG}$ and $\mathrm{KGE}$ are defined only when the cross-covariance vector is non-zero ($\vect{s}_p \neq \0_p$).}
\label{tab:training-performance}
{\resizebox{\textwidth}{!}{%
\begin{tabular}{lll}
\toprule
\multicolumn{3}{l}{Multiple-predictor linear model} \\
Metric & OLS linear regression (eq.~\eqref{eq:ols_predictions_multi}) & Kling-Gupta linear regression (eq.~\eqref{eq:kg_predictions_multi}) \\
\midrule
$\mathrm{MSE}(\vect{z}_n, \vect{y}_n)$ & $\sigma^2(\vect{y}_n)(1 - \rho_{\mathrm{OLS}_{(0)}}^2)$ & $2 \sigma^2(\vect{y}_n) (1 - \rho_{\mathrm{OLS}_{(0)}})$ \\
$\mathrm{NSE}(\vect{z}_n, \vect{y}_n)$ & $\rho_{\mathrm{OLS}_{(0)}}^2$ & $2 \rho_{\mathrm{OLS}_{(0)}} - 1$ \\
$L_\mathrm{KG}(\vect{z}_n, \vect{y}_n)$ & $2 (1 - \rho_{\mathrm{OLS}_{(0)}})^2, \vect{s}_p \neq \0_p$ & $(1 - \rho_{\mathrm{OLS}_{(0)}})^2$ \\
$\mathrm{KGE}(\vect{z}_n, \vect{y}_n)$ & $\sqrt{2} \rho_{\mathrm{OLS}_{(0)}} + 1 - \sqrt{2}, \vect{s}_p \neq \0_p$ & $\rho_{\mathrm{OLS}_{(0)}}$ \\
\midrule
\multicolumn{3}{l}{Single-predictor linear model} \\
Metric & OLS linear regression ($p = 1$) (eq.~\eqref{eq:ols_predictions_single}) & Kling-Gupta linear regression ($p = 1$) (eq.~\eqref{eq:kg_predictions_single}) \\
\midrule
$\mathrm{MSE}(\vect{z}_n, \vect{y}_n)$ & $\sigma^2(\vect{y}_n)(1 - \rho^2(\vect{x}_n, \vect{y}_n))$ & $2 \sigma^2(\vect{y}_n) (1 - |\rho(\vect{x}_n, \vect{y}_n)|)$ \\
$\mathrm{NSE}(\vect{z}_n, \vect{y}_n)$ & $\rho^2(\vect{x}_n, \vect{y}_n)$ & $2| \rho(\vect{x}_n, \vect{y}_n)| - 1$ \\
$L_\mathrm{KG}(\vect{z}_n, \vect{y}_n)$ & $2 (1 - |\rho(\vect{x}_n, \vect{y}_n)|)^2, \rho(\vect{x}_n, \vect{y}_n) \neq 0$ & $(1 - |\rho(\vect{x}_n, \vect{y}_n)|)^2$ \\
$\mathrm{KGE}(\vect{z}_n, \vect{y}_n)$ & $\sqrt{2} |\rho(\vect{x}_n, \vect{y}_n) | + 1 - \sqrt{2}, \rho(\vect{x}_n, \vect{y}_n) \neq 0$ & $|\rho(\vect{x}_n, \vect{y}_n)|$ \\
\bottomrule
\end{tabular}%
}
}
\end{table}

An evaluation of $\mathrm{NSE}$ values demonstrates a counterintuitive outcome. While OLS maximizes $\mathrm{NSE}$ by construction, the Kling-Gupta estimator produces an $\mathrm{NSE}$ of $2 \rho_{\mathrm{OLS}_{(0)}} - 1$, which remains lower than the OLS $\mathrm{NSE}$ of $\rho_{\mathrm{OLS}_{(0)}}^2$, because $0 \leq \rho_{\mathrm{OLS}_{(0)}} \leq 1$ from eqs.~\eqref{eq:ols_correlation_deriv} and \eqref{eq:rho_ols_star_def}. If $\rho_{\mathrm{OLS}_{(0)}}$ falls below $0.5$, the Kling-Gupta regression produces a negative $\mathrm{NSE}$. In such cases, the model performs worse than the mean climatology benchmark despite reproducing the observed mean and variance. Furthermore, $L_\mathrm{KG}$ for the OLS model equals $2 (1 - \rho_{\mathrm{OLS}_{(0)}})^2$, but the Kling-Gupta model halves this penalty to $(1 - \rho_{\mathrm{OLS}_{(0)}})^2$ by eliminating variability errors.

These results highlight a structural trade-off. Training a model with the Kling-Gupta loss leads to higher $\mathrm{KGE}$ values but markedly reduces $\mathrm{NSE}$ performance. Conversely, training with the squared error loss (OLS) optimizes $\mathrm{MSE}$ and $\mathrm{NSE}$ at the expense of the model's $\mathrm{KGE}$. This trade-off underscores the sensitivity of model rankings to the choice of evaluation metric, as discussed in Section~\ref{sec:estimation_and_evaluation}.
\subsection{Comparative performance on infinite training data}
\label{sec:comparative_infinite_training_data}
A relevant question is how the application of the limiting estimators defined in Section~\ref{sec:asymptotic_properties} influences model performance as the training sample size $n \to \infty$. Evaluating the metrics originally presented in Table~\ref{tab:training-performance} under this asymptotic limit facilitates large-sample comparisons between OLS and Kling-Gupta regressions. Table~\ref{tab:infinite-training-performance} summarizes the values calculated for this infinite-data limit.
\begin{table}[ht]
\centering
\caption{Performance of OLS and Kling-Gupta linear regression on the training set, summarizing the $\mathrm{MSE}$, $\mathrm{NSE}$, $L_\mathrm{KG}$, and $\mathrm{KGE}$ from Table~\ref{tab:training-performance} as $n \to \infty$. The upper panel covers the general multiple-predictor linear model, whereas the lower panel addresses the single-predictor special case. All entries are formulated using population statistics. For the OLS model, $L_\mathrm{KG}$ and $\mathrm{KGE}$ are defined only when the cross-covariance vector is non-zero ($\mathrm{Cov}_{F_{\underline{\vect{x}}_p, \underline{y}}}(\underline{\vect{x}}_p, \underline{y}) \neq \0_p$).}
\label{tab:infinite-training-performance}
{\resizebox{\textwidth}{!}{%
\begin{tabular}{lll}
\toprule
\multicolumn{3}{l}{Multiple-predictor linear model} \\
Metric & OLS linear regression (eq.~\eqref{eq:ols_predictions_multi}) & Kling-Gupta linear regression (eq.~\eqref{eq:kg_predictions_multi}) \\
\midrule
$\mathrm{MSE}(\vect{z}_{\infty}, \vect{y}_{\infty})$ & $\mathrm{Var}_{F_{\underline{y}}}[\underline{y}] (1 - \rho_{\mathrm{OLS}_{(0)}^*}^2)$ & $2 \mathrm{Var}_{F_{\underline{y}}}[\underline{y}](1 - \rho_{\mathrm{OLS}_{(0)}^*})$ \\
$\mathrm{NSE}(\vect{z}_{\infty}, \vect{y}_{\infty})$ & $ \rho_{\mathrm{OLS}_{(0)}^*}^2$ & $2 \rho_{\mathrm{OLS}_{(0)}^*} - 1$ \\
$L_\mathrm{KG}(\vect{z}_{\infty}, \vect{y}_{\infty})$ & $2 (1 - \rho_{\mathrm{OLS}_{(0)}^*})^2, \mathrm{Cov}_{F_{\underline{\vect{x}}_p, \underline{y}}}(\underline{\vect{x}}_p, \underline{y}) \neq \0_p$ & $(1 - \rho_{\mathrm{OLS}_{(0)}^*})^2$ \\
$\mathrm{KGE}(\vect{z}_{\infty}, \vect{y}_{\infty})$ & $\sqrt{2} \rho_{\mathrm{OLS}_{(0)}^*} + 1 - \sqrt{2}, \mathrm{Cov}_{F_{\underline{\vect{x}}_p, \underline{y}}}(\underline{\vect{x}}_p, \underline{y}) \neq \0_p$ & $\rho_{\mathrm{OLS}_{(0)}^*}$ \\
\midrule
\multicolumn{3}{l}{Single-predictor linear model} \\
Metric & OLS linear regression ($p = 1$) (eq.~\eqref{eq:ols_predictions_single}) & Kling-Gupta linear regression ($p = 1$) (eq.~\eqref{eq:kg_predictions_single}) \\
\midrule
$\mathrm{MSE}(\vect{z}_{\infty}, \vect{y}_{\infty})$ & $\mathrm{Var}_{F_{\underline{y}}}[\underline{y}] (1 - \mathrm{Corr}_{F_{\underline{x}, \underline{y}}}^2(\underline{x}, \underline{y}))$ & $2 \mathrm{Var}_{F_{\underline{y}}}[\underline{y}](1 - |\mathrm{Corr}_{F_{\underline{x}, \underline{y}}}(\underline{x}, \underline{y})|)$ \\
$\mathrm{NSE}(\vect{z}_{\infty}, \vect{y}_{\infty})$ & $\mathrm{Corr}^2_{F_{\underline{x}, \underline{y}}}(\underline{x}, \underline{y})$ & $2 |\mathrm{Corr}_{F_{\underline{x}, \underline{y}}}(\underline{x}, \underline{y})| - 1$ \\
$L_\mathrm{KG}(\vect{z}_{\infty}, \vect{y}_{\infty})$ & $2 (1 - |\mathrm{Corr}_{F_{\underline{x}, \underline{y}}}(\underline{x}, \underline{y})|)^2, \mathrm{Corr}_{F_{\underline{x}, \underline{y}}}(\underline{x}, \underline{y}) \neq 0$ & $(1 - |\mathrm{Corr}_{F_{\underline{x}, \underline{y}}}(\underline{x}, \underline{y})|)^2$ \\
$\mathrm{KGE}(\vect{z}_{\infty}, \vect{y}_{\infty})$ & $\sqrt{2} |\mathrm{Corr}_{F_{\underline{x}, \underline{y}}}(\underline{x}, \underline{y})| + 1 - \sqrt{2}, \mathrm{Corr}_{F_{\underline{x}, \underline{y}}}(\underline{x}, \underline{y}) \neq 0$ & $|\mathrm{Corr}_{F_{\underline{x}, \underline{y}}}(\underline{x}, \underline{y})|$ \\
\bottomrule
\end{tabular}%
}
}
\end{table}

The structural correspondence between Table~\ref{tab:training-performance} and Table~\ref{tab:infinite-training-performance} follows from the almost sure convergence of sample statistics to their population counterparts as $n \to \infty$, via the continuous mapping theorem. Specifically, the sample standard deviation $\sigma(\vect{y}_n)$ converges to $\sqrt{\mathrm{Var}_{F_{\underline{y}}}[\underline{y}]}$ (eq.~\eqref{eq:sample_standard_deviation_convergence}), while the sample OLS correlation $\rho_{\mathrm{OLS}_{(0)}}$ converges to its asymptotic limit $\rho_{\mathrm{OLS}_{(0)}^*}$ (eq.~\eqref{eq:rho_OLS_0_asymptotic}). Therefore, the performance trade-offs observed on the training set persist in the asymptotic limit, where population properties replace their sample-based analogs.
\subsubsection{Illustration of comparative performance on training data}
\begin{figure}[htbp]
\centering
\includegraphics{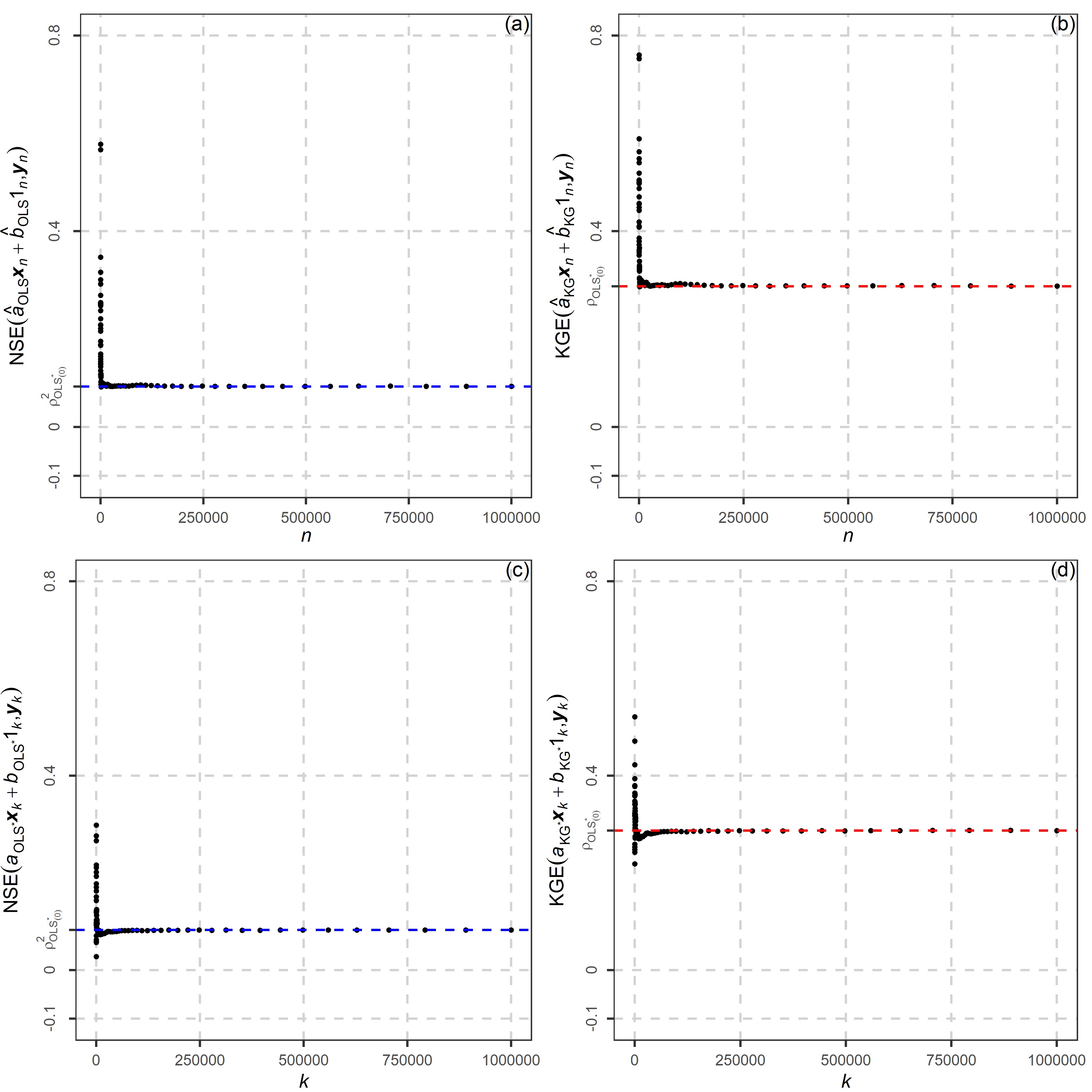}
\caption{Asymptotic convergence of performance metrics as a function of sample size for the training and test sets: (a) $\mathrm{NSE}$ of the OLS linear regression model on the training set versus sample size $n$, (b) $\mathrm{KGE}$ of the Kling-Gupta linear regression model on the training set versus sample size $n$, (c) $\mathrm{NSE}$ of the OLS linear regression model on the test set versus sample size $k$, (d) $\mathrm{KGE}$ of the Kling-Gupta linear regression model on the test set versus sample size $k$. The horizontal dashed lines indicate the population limits: $\rho_{\mathrm{OLS}_{(0)}^*}^2$ for the OLS $\mathrm{NSE}$ baseline and $\rho_{\mathrm{OLS}_{(0)}^*}$ for the Kling-Gupta $\mathrm{KGE}$ baseline.}
\label{fig:figure2}
\end{figure}

We investigate the behavior of the estimators by evaluating their empirical performance on the training data across a range of sample sizes $n$, using the predictions of the fitted models. To construct panels (a) and (b) of Figure~\ref{fig:figure2}, we fit the single-predictor model to the generated observations from Section~\ref{sec:numerical_illustration} sequentially for a sequence of selected sample sizes $n$. For each sample, we compute the $\mathrm{NSE}$ of the OLS predictions and the $\mathrm{KGE}$ of the Kling-Gupta predictions on the training set, using the expressions in Table~\ref{tab:training-performance}.

As panels (a) and (b) of Figure~\ref{fig:figure2} show, small sample sizes produce large fluctuations in the empirical metrics. As the sample size $n$ increases, these metrics stabilize and converge to their population values (Table~\ref{tab:infinite-training-performance}), confirming the asymptotic results listed therein.
\subsection{Asymptotic performance on independent test sets}
\label{sec:comparative_infinite_training_test_data}
Candidate models are evaluated on a test set to establish a performance ranking (Section~\ref{sec:strictly_consistent_loss_functions}). In finite samples, estimated model performance is in expectation worse on the test set than on the training set. Asymptotically, however, models estimated via $M$-estimation achieve equivalent performance on both sets, a consequence of the empirical risk minimization principle \citep{vapnik1998}. We show that this property extends to Kling-Gupta linear regression, even though the estimator does not belong to the $M$-estimator class. This finding implies that when a model is estimated using either $\mathrm{MSE}$ or $\mathrm{KGE}$, the expected performance determines model rankings, and those rankings depend on the loss functions chosen for both estimation and evaluation.

Given a test set of $k$ observations $\vect{y}_{k, \mathrm{test}} = (y_{1, \mathrm{test}}, \ldots, y_{k, \mathrm{test}})^\mathsf{T}$ and the predictor matrix $\mat{X}_{k \times p, \mathrm{test}}$, let predictions of the estimated models from eqs.~\eqref{eq:ols_test_predictions_multi} and \eqref{eq:kg_test_predictions_multi} be:
\begin{equation}
\label{eq:ols_n_infinity_estimated_lm}
\vect{z}_{k, \mathrm{OLS}, \mathrm{test}} = \mat{X}_{k \times p, \mathrm{test}} \vect{a}_{p, \mathrm{OLS}^*} + b_{p, \mathrm{OLS}^*} \1_k
\end{equation}
\begin{equation}
\label{eq:kg_n_infinity_estimated_lm}
\vect{z}_{k, \mathrm{KG}, \mathrm{test}} = \mat{X}_{k \times p, \mathrm{test}} \vect{a}_{p, \mathrm{KG}^*} + b_{p, \mathrm{KG}^*} \1_k
\end{equation}
where the parameter estimates $\widehat{\vect{a}}_{p, \mathrm{OLS}}$, $\widehat{b}_{p, \mathrm{OLS}}$, $\widehat{\vect{a}}_{p, \mathrm{KG}}$, and $\widehat{b}_{p, \mathrm{KG}}$ are replaced by their respective limits as $n \to \infty$. As $k \to \infty$, and assuming identical statistical properties for the training and test data, the test set loss function values (established in \hyperref[proof:b4]{Proof~B.4}) are identical to those reported in Table~\ref{tab:infinite-training-performance}.
\subsubsection{Asymptotic trade-offs and metric sensitivity}
\begin{figure}[htbp]
\centering
\includegraphics{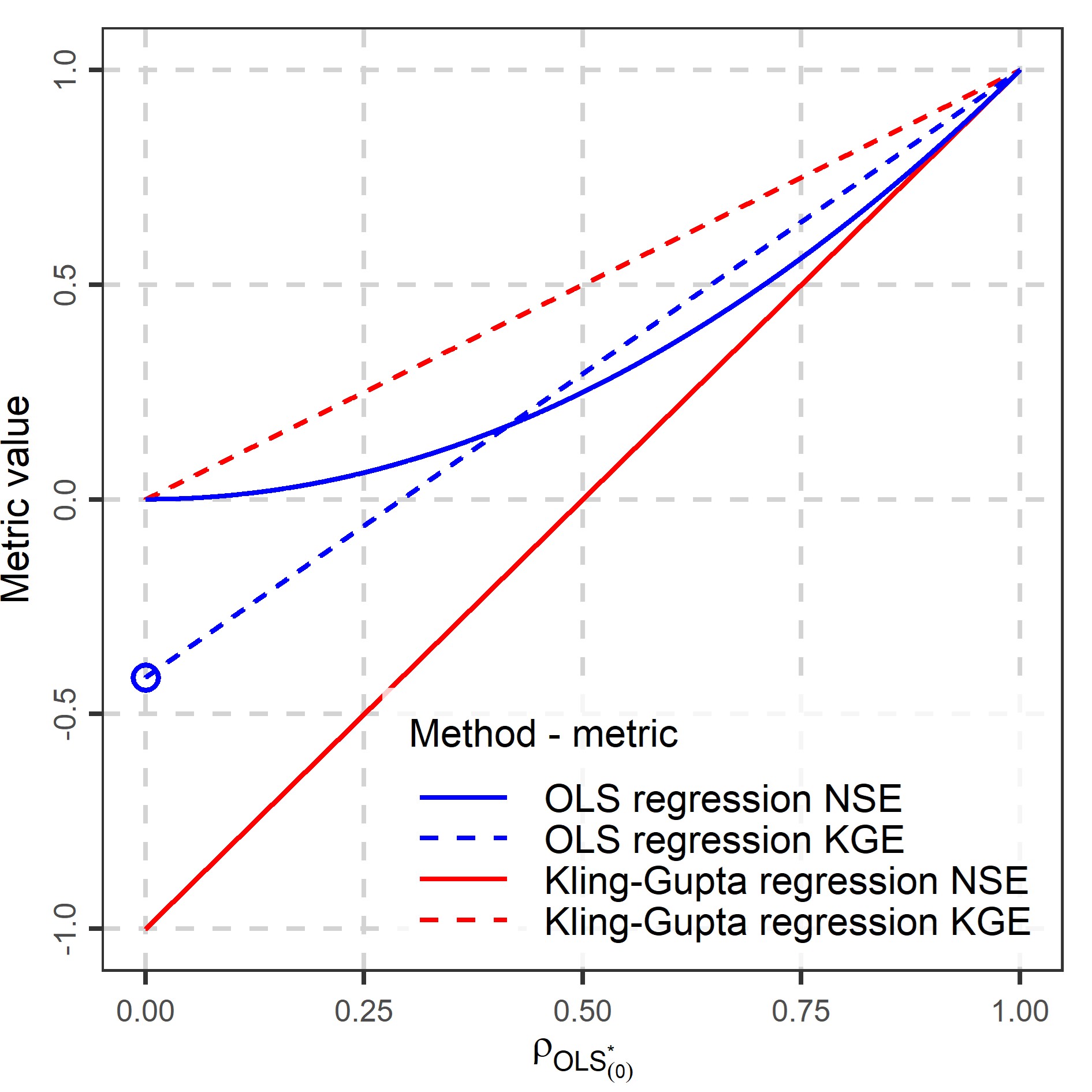}
\caption{Asymptotic performance metrics ($\mathrm{NSE}$ and $\mathrm{KGE}$) from Table~\ref{tab:infinite-training-performance} for OLS and Kling-Gupta linear regressions, plotted as a function of the OLS correlation coefficient $\rho_{\mathrm{OLS}_{(0)}^*}$. An open blue circle indicates the point where the $\mathrm{KGE}$ of the OLS model is undefined.}
\label{fig:figure3}
\end{figure}

Figure~\ref{fig:figure3} illustrates the trade-offs that govern the OLS and Kling-Gupta estimators in large samples, ($n \to \infty$), by plotting the asymptotic $\mathrm{NSE}$ and $\mathrm{KGE}$ of both regressions on the training set, thus reproducing the results of Table~\ref{tab:infinite-training-performance}. As established earlier in this section, the same pattern emerges when the estimated models are evaluated on test sets, as $n \to \infty$ and $k \to \infty$.

In the figure, blue curves represent the performance of OLS linear regression and red curves represent the performance of Kling-Gupta linear regression. Because every pair of response and predictor random variables maps to a unique value of $\rho_{\mathrm{OLS}_{(0)}^*}$, these curves describe the trade-offs between the two estimation objectives across the entire domain of the correlation. OLS regression maximizes $\mathrm{NSE}$ and outperforms Kling-Gupta regression in terms of $\mathrm{NSE}$. Conversely, Kling-Gupta regression maximizes $\mathrm{KGE}$ and outperforms OLS regression in terms of $\mathrm{KGE}$. For $\rho_{\mathrm{OLS}_{(0)}^*} < 1/2$, the $\mathrm{NSE}$ of Kling-Gupta regression becomes negative, whereas the $\mathrm{NSE}$ of OLS regression remains non-negative for all values of $\rho_{\mathrm{OLS}_{(0)}^*}$. Similarly, for $\rho_{\mathrm{OLS}_{(0)}^*} < 1 - \frac{1}{\sqrt{2}}$, the $\mathrm{KGE}$ of OLS regression becomes negative, while the $\mathrm{KGE}$ of Kling-Gupta regression stays non-negative for all $\rho_{\mathrm{OLS}_{(0)}^*}$.
\subsubsection{Asymptotic performance on test sets}
We evaluate the generalization capability of the models estimated in Section~\ref{sec:numerical_illustration} by measuring their performance on a test set as a function of the test sample size $k$. To construct panels (c) and (d) of Figure~\ref{fig:figure2}, we fix the parameter estimates from the full training sample and use them to generate predictions for a test dataset via eqs.~\eqref{eq:ols_test_predictions_single} and \eqref{eq:kg_test_predictions_single}. We then compute the test set $\mathrm{NSE}$ for the OLS model and the test set $\mathrm{KGE}$ for the Kling-Gupta model over an increasing sequence of test sample sizes $k$.

The empirical patterns in panels (c) and (d) of Figure~\ref{fig:figure2} show that for small $k$, the test metrics exhibit substantial variability and are lower than their population counterparts, but they increase as $k$ increases. As $k \to \infty$, the performance scores converge to the same theoretical population limits observed in the training phase. This confirms that the performance identities and trade-offs established in Table~\ref{tab:infinite-training-performance} extend to test data as $n, k \to \infty$, in agreement with the asymptotic results established in Sections~\ref{sec:comparative_infinite_training_data} and \ref{sec:comparative_infinite_training_test_data}.
\subsection{Linear model with fixed (known) regression coefficient}
\label{sec:fixed_regression_coefficient}
Treating the vector of regression coefficients $\vect{a}_p \in \R^{p}$ for the linear model in eq.~\eqref{eq:linear_model} as fixed a priori simplifies the estimation problem to determining the optimal intercept $b_p \in \R$ that minimizes the Kling-Gupta loss. \hyperref[proof:b5]{Proof~B.5} shows that the intercept appears exclusively in the bias term of $L_\mathrm{KG}$, while $b_p$ leaves the correlation and variability terms unaffected. Setting the bias term to zero leads to the intercept estimate:
\begin{equation}
\label{eq:intercept_estimate_fixed_a_KG}
\widehat{b}_{p, \mathrm{KG}} = \mu(\vect{y}_n) - (\vect{\mu}(\mat{X}_{n \times p}))^\mathsf{T} \vect{a}_p
\end{equation}
Minimizing the squared error loss under the same fixed coefficient constraint produces an identical intercept estimate:
\begin{equation}
\label{eq:intercept_estimate_fixed_a_OLS}
\widehat{b}_{p, \mathrm{OLS}} = \mu(\vect{y}_n) - (\vect{\mu}(\mat{X}_{n \times p}))^\mathsf{T} \vect{a}_p
\end{equation}
The equivalence of eqs.~\eqref{eq:intercept_estimate_fixed_a_KG} and \eqref{eq:intercept_estimate_fixed_a_OLS} indicates that both loss functions enforce predictions such that $\mu(\vect{z}_{n, \mathrm{OLS}}) = \mu(\vect{z}_{n, \mathrm{KG}}) = \mu(\vect{y}_n)$ via the intercept. Therefore, for a fixed $\vect{a}_p$, the selection of an estimator does not change the sample mean of the predictions.
\subsection{Fixed intercept Kling-Gupta linear regression}
\label{sec:fixed_intercept_kg_linear_regression}
Treating the intercept $b$ as fixed and known a priori simplifies the estimation task for the single-predictor linear model in eq.~\eqref{eq:single_pred_model} to identifying the slope $a \neq 0$ that minimizes the Kling-Gupta loss. For brevity, we define:
\begin{equation}
\label{eq:u_v_w_r_parameters}
u = \frac{\mu(\vect{x}_n)}{\mu(\vect{y}_n)}, v = \frac{\sigma(\vect{x}_n)}{\sigma(\vect{y}_n)} > 0, \rho = \rho(\vect{x}_n, \vect{y}_n) \text{ and } w = 1 - \frac{b}{\mu(\vect{y}_n)}
\end{equation}
The estimates of $\widehat{a}_\mathrm{KG}$ from \hyperref[proof:b6]{Proof~B.6} are summarized as:
\begin{equation}
\label{eq:slope_estimate_fixed_b_KG}
\widehat{a}_\mathrm{KG} =
\begin{cases}
\dfrac{wu + v}{u^2 + v^2}, & \text{if } \substack{(wu \geq v\ \text{and}\ \frac{(u - wv)^2}{u^2 + v^2} \leq 4\rho + w^2 + 1) \\ \text{or } (-v < wu < v\ \text{and}\ \rho > - \frac{wuv}{u^2 + v^2})} \\
\dfrac{wu - v}{u^2 + v^2}, & \text{if } \substack{(wu \leq -v\ \text{and}\ \frac{(u + wv)^2}{u^2 + v^2} \leq -4\rho + w^2 + 1) \\ \text{or } (-v < wu < v\ \text{and}\ \rho < - \frac{wuv}{u^2 + v^2})} \\
\dfrac{wu \pm v}{u^2 + v^2}, & \text{if } -v < wu < v\ \text{and}\ \rho = - \frac{wuv}{u^2 + v^2} \\
\text{no global minimum exists}, & \text{otherwise}
\end{cases}
\end{equation}
The functional form of the estimate depends on the values of $u$, $v$, $w$, and $\rho$. A global minimum for the estimation problem does not necessarily exist, as the initial three conditions in eq.~\eqref{eq:slope_estimate_fixed_b_KG} do not exhaust all possible realizations of $\vect{x}_n$ and $\vect{y}_n$. By comparison, the OLS estimate subject to the same fixed intercept is:
\begin{equation}
\widehat{a}_\mathrm{OLS} = \frac{\rho(\vect{x}_n, \vect{y}_n) \sigma(\vect{x}_n) \sigma(\vect{y}_n) + \mu(\vect{x}_n) (\mu(\vect{y}_n) - b)}{\mu^2(\vect{x}_n) + \sigma^2(\vect{x}_n)}
\end{equation}
which is unique and finite.

Analyzing the fixed-intercept case offers insight into how the Kling-Gupta loss operates independently of mean alignment, which the intercept selection in Section~\ref{sec:kling_gupta_linear_regression_single_predictor} enforces. In this setting, the fixed $b$ cannot be optimized to eliminate bias; therefore, the loss function explicitly penalizes any deviation of the predictive mean from the observed mean through the bias term $(w - au)^2$. This establishes an interaction between the correlation term, which favors a positive sign for $a$ when $\rho > 0$ and a negative sign when $\rho < 0$, and the bias term, which might incentivize slopes that shift the predictive mean toward $\mu(\vect{y}_n)$. Simultaneously, the variability term penalizes values of $|a|$ that deviate from $1/v$. The potential absence of a global minimum constitutes an inherent property of the Kling-Gupta loss under a fixed intercept, arising from the discontinuity at $a = 0$ (which is proved in \hyperref[proof:b6]{Proof~B.6}) combined with the loss function's geometry. If the loss decreases monotonically as $a \to 0$ from a specific direction, the infimum is approached but not reached, as $a = 0$ lies outside the domain (since $\sigma(\vect{z}_n) = 0$ violates the Kling-Gupta loss definition).

Comparing the Kling-Gupta and OLS estimates in this fixed-intercept setting highlights the differences in their optimization objectives. The OLS estimate is always defined and minimizes the $\mathrm{MSE}$, balancing squared bias and variance without imposing explicit penalties on the correlation or the standard deviation ratio. The Kling-Gupta estimate, when it exists, prioritizes a high absolute correlation and a standard deviation ratio near unity, often at the cost of mean bias if the fixed intercept is suboptimally selected. The piecewise conditions identifying the global minimizer, such as the comparison between $\rho$ and $-\frac{wuv}{u^2 + v^2}$, signify trade-offs between correlation sign, magnitude, mean alignment, and scale correspondence. For instance, in the intermediate regime where $wu + v > 0$ and $wu < v$, the two candidate slopes $a_+$ and $a_-$ signify positive and negative transformations of the predictor, respectively. The loss favors the positive slope if the correlation is sufficiently high relative to a threshold constructed from the standardized means and standard deviations; otherwise, the negative slope becomes optimal, reversing the relationship sign to exploit the absolute correlation. This threshold explicitly quantifies the point at which the benefit of aligning with the sign of $\rho$ outweighs the bias introduced by the selection of one slope over the other.
\subsubsection{Numerical illustration}
\begin{figure}[htbp]
\centering
\includegraphics{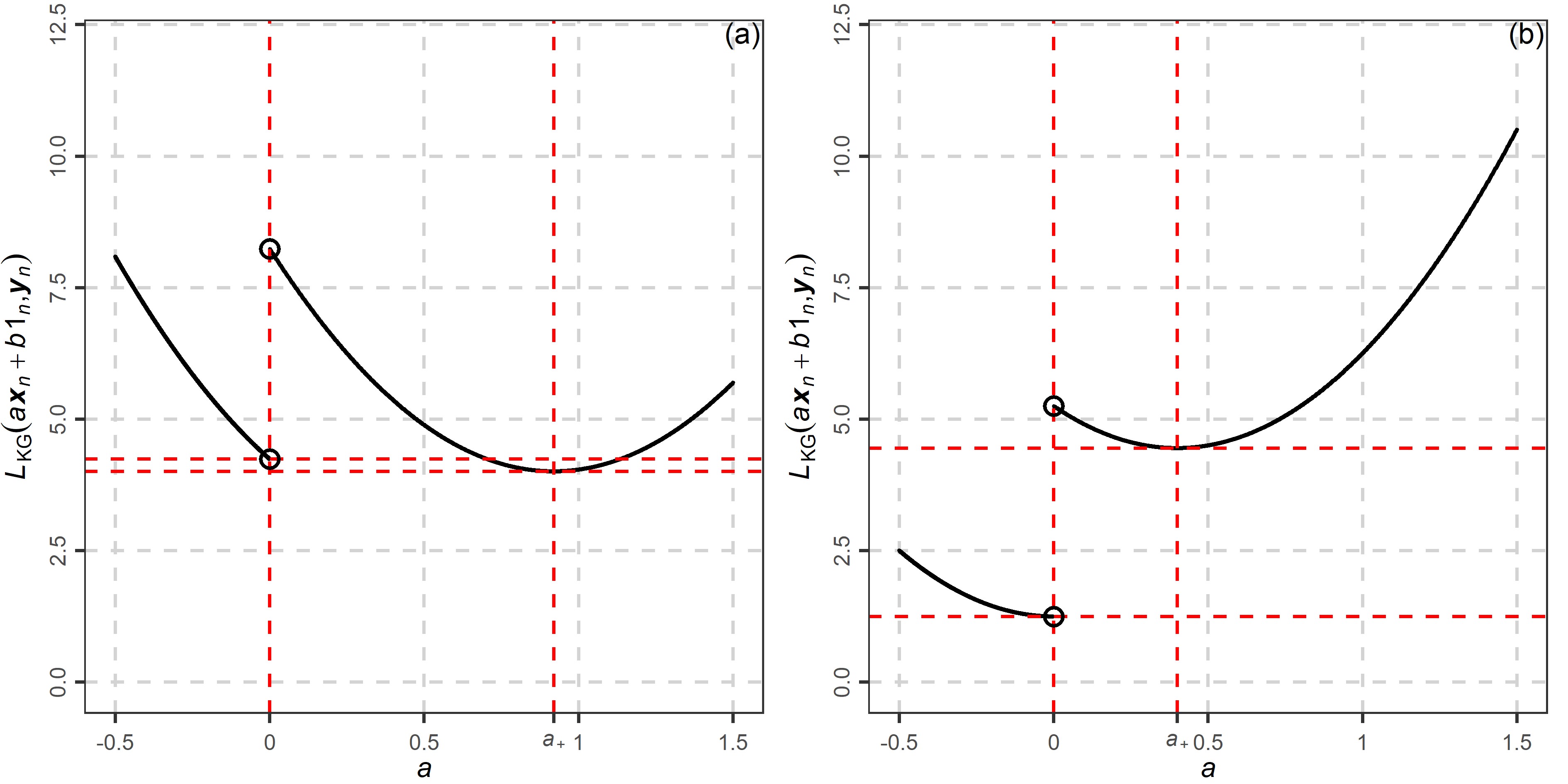}
\caption{Kling-Gupta loss $L_\mathrm{KG}(a \vect{x}_n + b_0 \1_n, \vect{y}_n)$ as a function of $a$ for a single-predictor linear model with sample statistics $\mu(\vect{x}_n) = 2$, $\mu(\vect{y}_n) = 1$, $\sigma(\vect{x}_n) = 1$, $\sigma(\vect{y}_n) = 1$, and $\rho(\vect{x}_n, \vect{y}_n) = -1$: (a) $b_0 = -0.8$, and (b) $b_0 = 0.5$. The positive local optimizer $a_+$ is given by eq.~\eqref{eq:a_positive_candidate_slope}. Open black circles indicate the points where the Kling-Gupta loss is undefined.}
\label{fig:figure4}
\end{figure}

As an example of the fixed-intercept parameter setting, Figure~\ref{fig:figure4} illustrates the Kling-Gupta loss, specifically its behavior as a function of the slope parameter $a$ when the intercept $b_0$ is fixed, using the model $z = a x + b_0$ from eq.~\eqref{eq:single_pred_model}. In this example, we visualize $L_\mathrm{KG}(a \vect{x}_n + b_0 \1_n, \vect{y}_n)$ as a function of $a$. The sample statistics for the predictor random variable $\underline{x}$ and the response random variable $\underline{y}$ are set as follows: $\mu(\vect{x}_n) = 2$, $\mu(\vect{y}_n) = 1$, $\sigma(\vect{x}_n) = 1$, $\sigma(\vect{y}_n) = 1$, and $\rho(\vect{x}_n, \vect{y}_n) = -1$. We plot the loss for two fixed, known values of $b_0$: panel (a) sets $b_0 = -0.8$ and panel (b) sets $b_0 = 0.5$. The loss reduces to eq.~\eqref{eq:KG_with_b_fixed} and is shown in Figure~\ref{fig:figure4}, which demonstrates its non-convex nature and the discontinuity at $a = 0$. The positive local optimizer $a_+$ (from eq.~\eqref{eq:a_positive_candidate_slope}) is also shown, along with the limiting value of the loss as $a \to 0^-$.

In Figure~\ref{fig:figure4}a, where $b_0 = -0.8$, the loss function exhibits a well-defined global minimum at $a_+$. Although a discontinuity exists at $a = 0$, the loss value at $a_+$ lies below the limiting values approached from both the positive and negative sides of $a = 0$.

Conversely, Figure~\ref{fig:figure4}b illustrates the case $b_0 = 0.5$. Here, the loss function decreases monotonically as $a$ approaches zero from the negative side, reaching an infimum lower than any value of the Kling-Gupta loss for $a > 0$. However, the Kling-Gupta metric is undefined at $a = 0$ and discontinuous. Therefore, the loss function cannot attain a global minimum within the permissible parameter space, because the true infimum lies on the boundary of an open set where the estimator breaks down.
\section{Applications}
\label{sec:applications}
This section verifies the theoretical properties established in Section~\ref{sec:kling_gupta_regressions} using observed daily streamflow time series from ten French catchments with complete records. We evaluate the training and test performance of OLS and Kling-Gupta linear regressions using models with either a single lagged predictor or two lagged predictors. Section~\ref{sec:applications_data} details the dataset, Section~\ref{sec:applications_methods} outlines the estimation procedure, and Section~\ref{sec:applications_results} reports the numerical results. Code and data required for replication accompany this manuscript as supplementary material. We performed all computations in the \texttt{R} programming language (version 4.6.0) within the \texttt{RStudio} environment (version 2026.05.0+218); Appendix~\ref{app:software} documents the specific packages used.
\subsection{Data}
\label{sec:applications_data}
We forecast daily mean discharge records ($Q_t$, mm/d) from ten French catchments over the 1999--2018 period, sourced from the \texttt{airGRdatasets} R package \citep{delaigue2025}. For each catchment, we construct two lagged variables, $Q_{t-1}$ and $Q_{t-2}$. To estimate the linear models that forecast $Q_t$ conditional on these lagged variables, we omit the first two observations to remove missing values from the predictor matrix. The remaining series is then partitioned into a training sample spanning 1999--2008 ($n = 3651$ days) and a test sample spanning 2009--2018 ($k = 3652$ days).
\subsection{Application}
\label{sec:applications_methods}
We estimate four linear models for each catchment's training data, taking $Q_t$ as the response variable:

\noindent\makebox[1cm][l]{(i)} OLS linear regression with a single predictor $Q_{t-1}$ (Section~\ref{sec:ols-linear-regression-single-predictor}).

\noindent\makebox[1cm][l]{(ii)} OLS linear regression with two predictors $Q_{t-1}$ and $Q_{t-2}$ (Section~\ref{sec:ols-linear-regression}).

\noindent\makebox[1cm][l]{(iii)} Kling-Gupta linear regression with a single predictor $Q_{t-1}$ (Section~\ref{sec:kling_gupta_linear_regression_single_predictor}).

\noindent\makebox[1cm][l]{(iv)} Kling-Gupta linear regression with two predictors $Q_{t-1}$ and $Q_{t-2}$ (Section~\ref{sec:kling_gupta_linear_regression_multiple_predictors}).

Following parameter estimation on the training sample, we generate both training and test set predictions. For each model and catchment, we evaluate $\mathrm{NSE}$ and $\mathrm{KGE}$ performances. We omit reporting the $\mathrm{MSE}$ and the Kling-Gupta loss ($L_\mathrm{KG}$) due to their monotonic relationships with $\mathrm{NSE}$ and $\mathrm{KGE}$, respectively. Additionally, we compute the bias $1 - \mu(\vect{z}_n)/\mu(\vect{y}_n)$, variability $1 - \sigma(\vect{z}_n)/\sigma(\vect{y}_n)$, and correlation $1 - \rho(\vect{z}_n, \vect{y}_n)$ components of $\mathrm{KGE}$. These are the non-squared counterparts of the terms that appear in $\mathrm{KGE}$.
\subsection{Results}
\label{sec:applications_results}
\begin{figure}[htbp]
\centering
\includegraphics{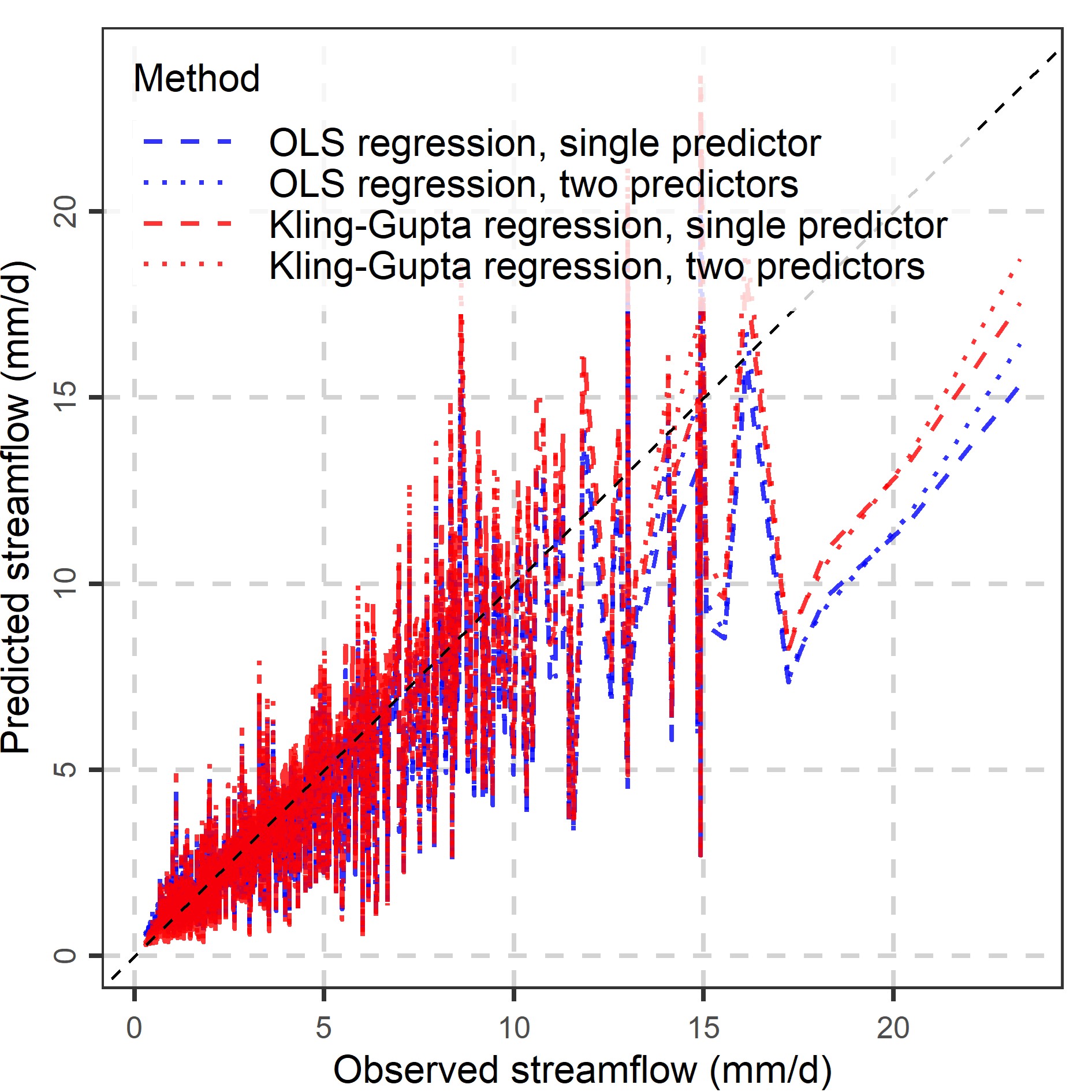}
\caption{Scatterplot comparing observed and predicted streamflow for catchment A273011002 on the test set. Predictions are generated using the four linear regression models detailed in Section~\ref{sec:applications_methods}. The dashed black line represents the 1:1 line of perfect agreement ($z = y$).}
\label{fig:figure5}
\end{figure}

\begin{figure}[htbp]
\centering
\includegraphics{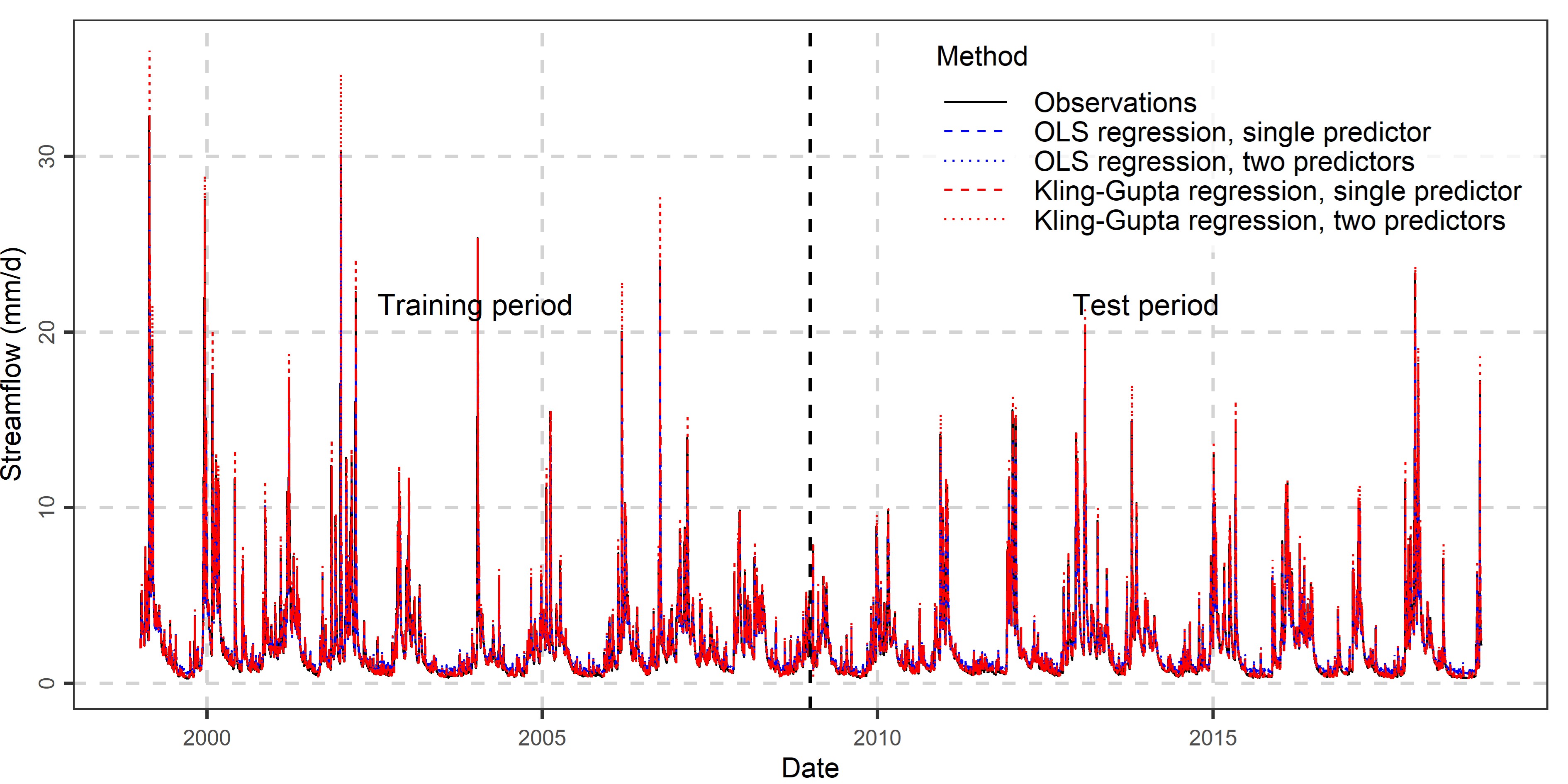}
\caption{Time series plot of observed streamflow and predictions generated by the four linear regression models from Section~\ref{sec:applications_methods} for catchment A273011002. The vertical dashed line marks the transition between the training period (1999--2008) and the test period (2009--2018).}
\label{fig:figure6}
\end{figure}

Figures~\ref{fig:figure5} and~\ref{fig:figure6} illustrate the behavior of the four models for a randomly selected catchment (named as A273011002 in the dataset). In Figure~\ref{fig:figure5}, Kling-Gupta predictions spread more widely compared to OLS predictions. The time series plots in Figure~\ref{fig:figure6} confirms that Kling-Gupta predictions have higher variance than OLS.

\begin{figure}[htbp]
\centering
\includegraphics[scale=1.00]{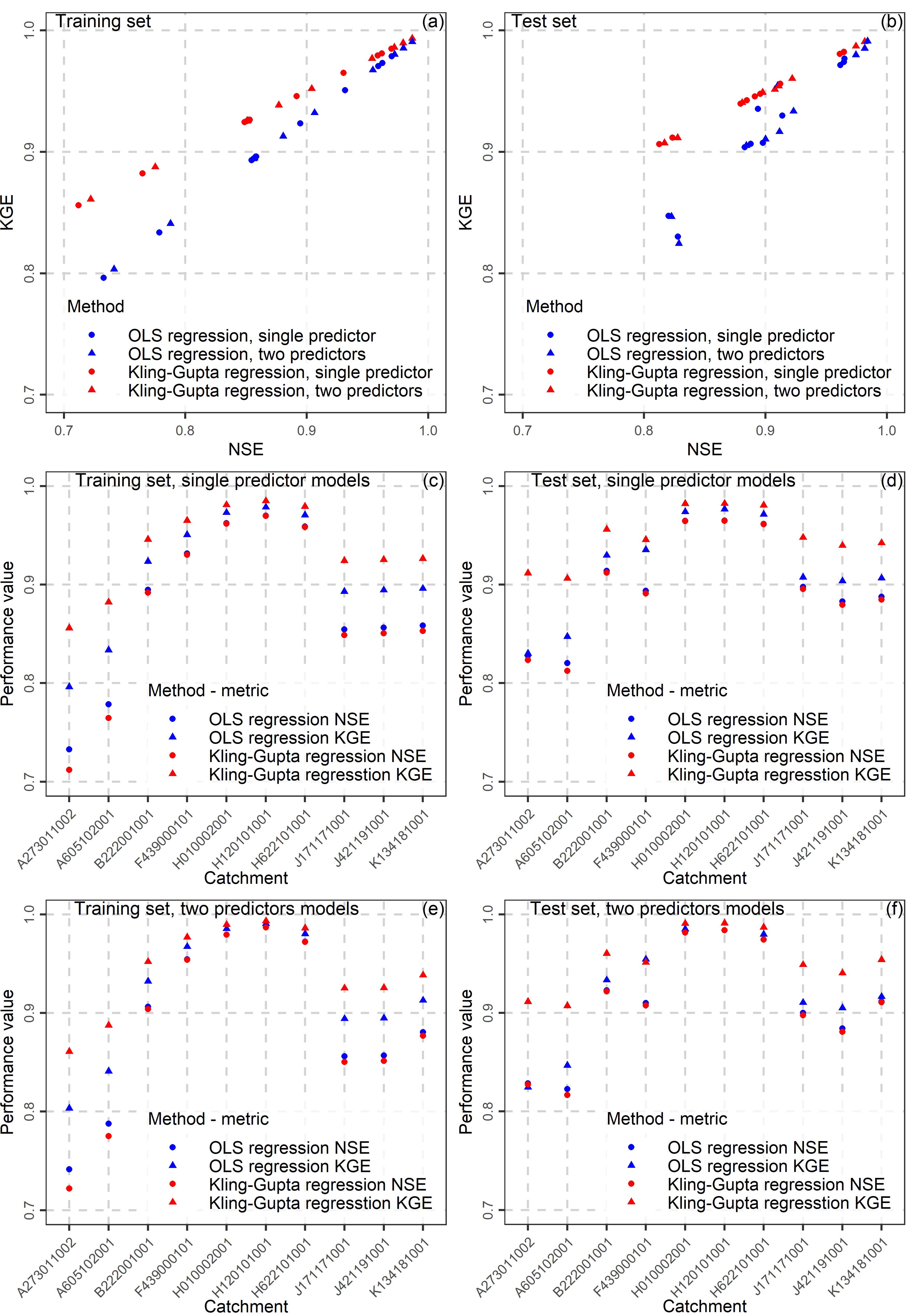}
\caption{Performance comparison of the OLS and Kling-Gupta linear regression models from Section~\ref{sec:applications_methods} across ten French catchments evaluated on both training and test sets: (a, b) $\mathrm{NSE}$ versus $\mathrm{KGE}$, (c, d) $\mathrm{NSE}$ and $\mathrm{KGE}$ for single-predictor models, and (e, f) $\mathrm{NSE}$ and $\mathrm{KGE}$ for two-predictor models by catchment.}
\label{fig:figure7}
\end{figure}

Figure~\ref{fig:figure7} summarizes the $\mathrm{NSE}$ and $\mathrm{KGE}$ values across all ten catchments. Panels (a) and (b) show that, for both training and test sets, the Kling-Gupta models consistently achieve higher $\mathrm{KGE}$ than their OLS counterparts, whereas the OLS models attain higher $\mathrm{NSE}$. The $\mathrm{NSE}$ and $\mathrm{KGE}$ performance on the training set exhibits a linear relationship, as expected from the analytical expressions in Table~\ref{tab:training-performance}. In the test set, the linear relationship is slightly disturbed, especially for OLS regression, due to the inherent error induced when testing trained models on finite samples. Panels (c) and (e) (single-predictor and two-predictor models, respectively) illustrate more clearly that Kling-Gupta regression outperforms OLS regression on the training set in terms of $\mathrm{KGE}$, while OLS regression outperforms Kling-Gupta regression in terms of $\mathrm{NSE}$. The ranking is maintained for every catchment, and the differences become less pronounced as $\mathrm{NSE}$ and $\mathrm{KGE}$ approach unity. The reason for the less pronounced differences is that the predictions begin to approach the observations. The corresponding panels (d) and (f) show a similar pattern for the test set, although the ranking may change slightly, again due to the randomness introduced by testing on a finite test sample.

\begin{figure}[htbp]
\centering
\includegraphics[scale=1.00]{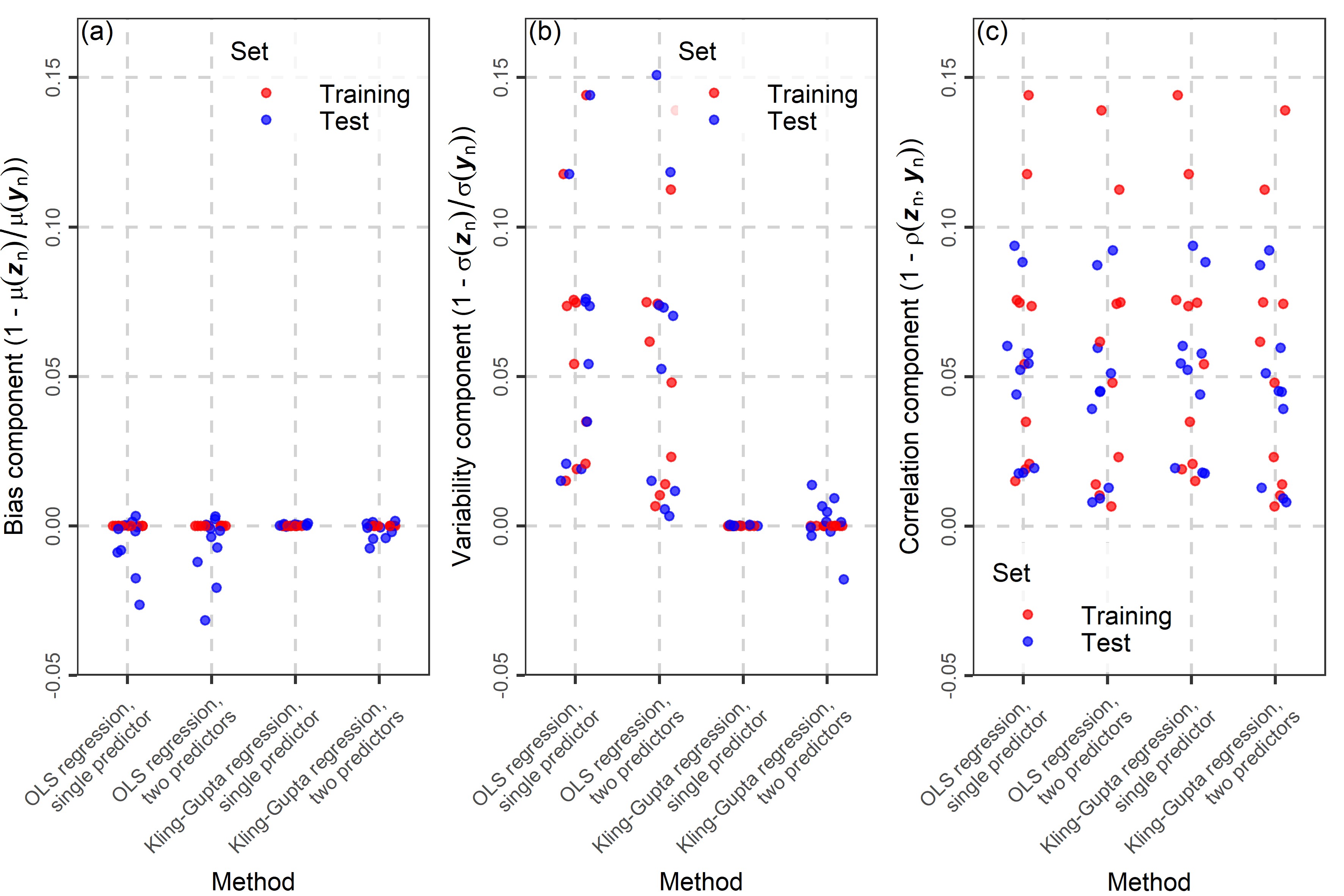}
\caption{Decomposition of the Kling-Gupta loss into its three components for the OLS and Kling-Gupta linear regression models from Section~\ref{sec:applications_methods}: (a) bias component $1 - \frac{\mu(\vect{z}_n)}{\mu(\vect{y}_n)}$, (b) variability component $1 - \frac{\sigma(\vect{z}_n)}{\sigma(\vect{y}_n)}$, and (c) correlation component $1 - \rho(\vect{z}_n, \vect{y}_n)$. Each panel displays both training and test set results across all ten catchments.}
\label{fig:figure8}
\end{figure}

Figure~\ref{fig:figure8} examines the three components of the Kling-Gupta loss: the bias component $1 - \frac{\mu(\vect{z}_n)}{\mu(\vect{y}_n)}$, the variability component $1 - \frac{\sigma(\vect{z}_n)}{\sigma(\vect{y}_n)}$, and the correlation component $1 - \rho(\vect{z}_n, \vect{y}_n)$. On the training set, the bias component (panel a) is near zero for all methods, as expected from theory. On the test set, due to its finite size, the bias component may deviate slightly from zero, although the general pattern remains. The variability component is systematically positive for OLS, indicating variance reduction, and essentially zero for Kling-Gupta (variance preservation). This pattern is maintained on the training set and deviates only slightly on the test set, again because the test set is finite. The correlation component (panel c) is small and similar across methods that use the same information: methods with a single predictor have equal correlation components, methods with two predictors have equal correlation components among themselves, but when compared with methods that use less information, they exhibit smaller correlation components. These empirical findings align with the theoretical properties proved in Sections~\ref{sec:kling_gupta_linear_regression_multiple_predictors} and~\ref{sec:comparative_training_data}.

In most cases, the $\mathrm{NSE}$ and $\mathrm{KGE}$ values on the test set appear higher than those on the training set. This finding is due to the specificities of the paper's dataset and the randomness of results in finite samples; it is not a general rule, as the expected performance on a finite test set should be worse than that on the training set.
\section{Discussion}
\label{sec:discussion}
\subsection{Comparison with extant literature}
Historically, hydrologic literature treated $\mathrm{KGE}$ as an ``informal'' goodness-of-fit metric or empirical diagnostic tool \citep{gupta2009, vrugt2022}. Prior studies analyzed metric sensitivity and parameter uncertainty through synthetic simulations \citep{knoben2019}, resampling architectures like the jackknife or bootstrap \citep{clark2021, vrugt2022}, or numerical algorithms. We departed from these computational approaches after developing a mathematical foundation that framed the negatively oriented Kling-Gupta loss $L_\mathrm{KG}$ as a statistically consistent extremum estimator. Within this parametric setting, we established analytical, closed-form expressions for the regression parameters, eliminating the need for iterative search procedures.

Kling-Gupta and OLS linear regressions represent distinct optimization problems within the extremum estimation framework. Both estimators enforce identical predictive sample means and identical sample correlations with the response variable. However, the Kling-Gupta coefficient vector scales the OLS vector by the variance-inflation factor $\sigma(\vect{y}_n)/\sigma_\mathrm{OLS}$. This scaling counteracts the systematic variance reduction inherent to OLS estimation, ensuring that the model predictions replicate the marginal sample variance of the observations. This structural divergence introduces a trade-off in predictive performance: the Kling-Gupta estimator maximizes the $\mathrm{KGE}$ by eliminating variability errors, but at the cost of decreasing the $\mathrm{NSE}$.

We contextualize the empirical observations reported by \citep{gupta2009}, based on hydrologic model calibrations using both $\mathrm{NSE}$ and $\mathrm{KGE}$, through our linear model framework:

\noindent\makebox[1cm][l]{(i)} \citep{gupta2009} observed empirically that parameter estimates changed minimally when training with either $\mathrm{NSE}$ or $\mathrm{KGE}$. Structurally, for a single-predictor linear model, the OLS parameter estimates converge to a common limit as $\rho(\vect{x}_n, \vect{y}_n) \to \pm 1$ (i.e., when the absolute predictor-response correlation approaches unity). We hypothesize that their findings reflect highly correlated variables in their experimental data.

\noindent\makebox[1cm][l]{(ii)} They reported that both OLS and Kling-Gupta estimators produced similar training-set correlations between predictions and observations, independent of the calibration objective. For linear models, we show that these sample correlations are identical.

\noindent\makebox[1cm][l]{(iii)} They observed that the predictive variance from Kling-Gupta calibration aligns with the observed variance, whereas OLS training reduces it. Within our linear framework, the Kling-Gupta predictive variance equals the sample variance of the observations identically, while the OLS predictive variance remains systematically lower.

\noindent\makebox[1cm][l]{(iv)} Their Figure 8d indicates that the predictive sample means for both estimators align with the observed mean, with closer alignment under $\mathrm{KGE}$. We demonstrate that for any linear model with an unconstrained intercept, both estimation procedures achieve perfect mean alignment.

\noindent\makebox[1cm][l]{(v)} Their Figure 6 illustrates structural links between $\mathrm{NSE}$ and $\mathrm{KGE}$ across different calibration schemes, reporting the training-set relationship $\mathrm{NSE}(\vect{z}_n, \vect{y}_n) = 2 \mathrm{KGE}(\vect{z}_n, \vect{y}_n) - 1$ under $\mathrm{KGE}$ training. For linear regressions, this expression follows immediately from Table~\ref{tab:training-performance} for finite samples and from Table~\ref{tab:infinite-training-performance} asymptotically. Analogous performance identities can be directly formulated from these tables for other metric and estimator configurations.
\subsection{Thinking about the choice between NSE and KGE}
We demonstrated that Kling-Gupta optimization induces a statistical estimator, distinct from OLS regression. Within a linear framework, the Kling-Gupta estimator maintains the direction of the OLS coefficient vector but scales it by a deterministic variance-inflation factor to enforce equality between the empirical variances of the predictions and observations. Therefore, we interpret Kling-Gupta regression as a transformation of OLS that retains both the sample correlation and mean unbiasedness while counteracting the systematic variance reduction inherent to squared-error minimization. This framework clarifies the statistical behavior of $\mathrm{KGE}$-calibrated models and establishes a mathematical foundation explaining why $\mathrm{KGE}$ optimization generates higher variability and augmented peak magnitudes compared to calibration under the $\mathrm{NSE}$.

Although this variance-inflation property appears logical in hindsight, it is counterintuitive a priori. As the $\mathrm{KGE}$ combines bias, variability, and correlation terms, we might expect a complex optimization trade-off that simultaneously alters all three properties. Unexpectedly, applying the $\mathrm{KGE}$ to a linear model leaves the predictive sample mean and correlation unaffected relative to OLS.

We showed that the empirical divergences reported in the hydrologic literature between $\mathrm{NSE}$- and $\mathrm{KGE}$-optimized models are necessary structural consequences of their respective loss functions, rather than data- or model-specific artifacts. Therefore, we view systematic disagreements between the $\mathrm{NSE}$ and $\mathrm{KGE}$ values not as methodological issues, but as the logical outcome of estimators optimized toward distinct targets.

While modelers often ask, ``Which metric is better?'', we frame this question as: ``What property does each metric optimize?'' We do not advocate replacing the $\mathrm{NSE}$ with the $\mathrm{KGE}$, or vice versa; instead, we characterize the statistical consequences of each choice. We demonstrate that selecting between the $\mathrm{NSE}$ and $\mathrm{KGE}$ is a choice between two distinct predictive objectives: variance reduction ($\mathrm{NSE}$) versus variance replication ($\mathrm{KGE}$). It is therefore important that the modeler is given the evaluation metric a priori. Ultimately, the choice of the metric depends on the forecast end-user, based on the specific operational properties required of the predictions.
\subsection{Future work}
We analyzed linear models to establish a baseline where the statistical mechanisms of Kling-Gupta optimization are analytically tractable. Because real-world hydrologic systems are nonlinear, they are not fully represented by these simplified assumptions. We present these results as a focused exposition of how $\mathrm{KGE}$ objectives alter predictive behavior in a stylized environment. Nonetheless, the empirical behavior of nonlinear hydrologic models estimated with $\mathrm{KGE}$ frequently mirrors the variance-inflation and hydrograph expansion of our linear framework. This consistency suggests that our relationships reflect general structural properties of variance-constrained loss functions. Formal extension of these results to nonlinear estimators remains an open problem and a natural direction for future research.

Another key open question concerns the conditional statistical functional estimated by the Kling-Gupta loss. While OLS regression is known to target the conditional mean, the functional optimized by $\mathrm{KGE}$ is less transparent, although we have shown that Kling-Gupta linear regression identifies the marginal mean. Future work should investigate whether $\mathrm{KGE}$ optimization can be formally interpreted as targeting a specific functional of the conditional distribution. Answering this question would place $\mathrm{KGE}$ on firmer theoretical ground and better clarify its use in semiparametric prediction problems.
\section{Conclusions}
\label{sec:conclusions}
We formulated a theoretical framework for the Kling-Gupta efficiency ($\mathrm{KGE}$), characterizing it as an extremum estimator for linear models. By defining the negatively oriented loss $L_\mathrm{KG} = (1 - \mathrm{KGE})^2$, we established closed-form parameter estimates showing that Kling-Gupta linear regression scales the OLS coefficient vector by the variance-inflation factor $\sigma(\vect{y}_n)/\sigma_\mathrm{OLS}$. This variance-inflation mechanism ensures that predictions on the training set reproduce the marginal sample variance of the observations while maintaining the predictive mean and correlation properties of OLS. Therefore, the estimator introduces a trade-off, maximizing the $\mathrm{KGE}$ while decreasing the $\mathrm{NSE}$. We proved that $L_\mathrm{KG}$ attains a unique global minimum of zero if and only if $\vect{z}_n = \vect{y}_n$. We analytically computed the finite training set performance and showed that the asymptotic performances on the training and test sets are equal.

These findings have practical value for hydrologic science. The closed-form parameter estimates eliminate the need for iterative numerical calibration searches for linear models, offering computational efficiency and interpretability. The trade-offs between the $\mathrm{NSE}$ and $\mathrm{KGE}$ explain empirical literature reports of metric disagreement, offering modelers a principled foundation for selecting loss functions tailored to their needs. Furthermore, our analysis of constrained estimation, specifically the fixed-intercept setting, indicates that a global minimum for the slope parameter does not always exist, highlighting a previously unaddressed limitation of $\mathrm{KGE}$-based optimization when mean alignment is constrained.
\appendix
\titleformat{\section}
  {\bfseries\fontsize{12}{18}\selectfont}
  {Appendix \thesection}
  {0.5em}
  {}
\titlespacing*{\section}{0pt}{6pt}{6pt}

\section{Vector, matrix and statistical notation}
\label{app:notation}
This appendix establishes the notation for vector and matrix operations, norms, and statistical functions employed throughout the manuscript.

\textbf{\underline{Basic vector definitions}}

Vector notation:
\begin{equation}
\label{eq:vecdef}
\vect{x}_n := (x_1, \ldots, x_n)^\mathsf{T}
\end{equation}
where the superscript $\mathsf{T}$ indicates the transpose.

Zero vector (all elements zero):
\begin{equation}
\label{eq:zerovec}
\0_n := (0, \ldots, 0)^\mathsf{T}
\end{equation}
All-ones vector (all elements unity):
\begin{equation}
\label{eq:onesvec}
\1_n := (1, \ldots, 1)^\mathsf{T}
\end{equation}

\textbf{\underline{Basic matrix definitions}}

Matrix notation:
\begin{equation}
\label{eq:matrix}
\mat{X}_{n \times m} := [x_{ij}], i = 1, \ldots, n, j = 1, \ldots, m
\end{equation}
Matrix row:
\begin{equation}
\label{eq:matrixrow}
\mat{X}_{i, \bullet} := (x_{i1}, \ldots, x_{im}), i = 1, \ldots, n
\end{equation}
Matrix column:
\begin{equation}
\label{eq:matrixcol}
\mat{X}_{\bullet ,j} := (x_{1j}, \ldots, x_{nj})^\mathsf{T}, j = 1, \ldots, m
\end{equation}
Identity matrix ($\delta_{ij}$ is the Kronecker delta):
\begin{equation}
\mat{I}_n := [\delta_{ij}], \delta_{ij} := \begin{cases} 1, & \text{if } 1 \leq i = j \leq n \\ 0 & \text{otherwise} \end{cases}
\end{equation}

\textbf{\underline{Element-wise operations on vectors}}

Element-wise inequality:
\begin{equation}
\label{eq:vecne}
\vect{x}_n \neq \vect{y}_n \Leftrightarrow \exists i \in \{1, \ldots, n \}: x_i \neq y_i
\end{equation}

\textbf{\underline{Inner products and norms}}

Euclidean inner product:
\begin{equation}
\left\langle \vect{x}_n, \vect{y}_n \right\rangle := \vect{x}_n^\mathsf{T} \vect{y}_n = \sum_{i = 1}^n x_i y_i = \left\langle \vect{y}_n, \vect{x}_n \right\rangle
\end{equation}
\begin{equation}
\left\langle a \vect{x}_n, \vect{y}_n \right\rangle = a \left\langle \vect{x}_n, \vect{y}_n \right\rangle
\end{equation}
Euclidean norm:
\begin{equation}
\left\|\vect{x}_n \right\|_2 := \sqrt{(\sum_{i = 1}^n x_i^2)}
\end{equation}
Squared Euclidean norm:
\begin{equation}
\left\|\vect{x}_n \right\|_2^2 = \left\langle \vect{x}_n, \vect{x}_n \right\rangle = \vect{x}_n^\mathsf{T} \vect{x}_n = \sum_{i = 1}^n x_i^2
\end{equation}
Homogeneity:
\begin{equation}
\left\|a \vect{x}_n \right\|_2^2 = a^2 \left\|\vect{x}_n \right\|_2^2
\end{equation}
Norm equivalences:
\begin{equation}
\left\|\left|\vect{x}_n \right| \right\|_2 = \left\|\vect{x}_n \right\|_2
\end{equation}
\begin{equation}
\left\langle|\vect{x}_n|, |\vect{x}_n| \right\rangle = \left\langle \vect{x}_n,\vect{x}_n \right\rangle = \left\|\vect{x}_n \right\|_2^2
\end{equation}

\textbf{\underline{Sign function}}

Scalar $\operatorname{sign}$ function:
\begin{equation}
\label{eq:sign}
\operatorname{sign}(x) := \begin{cases}
1, & \text{if } x > 0 \\
-1, & \text{if } x < 0
\end{cases}
\end{equation}

\textbf{\underline{Statistical functions for vectors}}

Sample mean:
\begin{equation}
\label{eq:samplemean}
\mu(\vect{x}_n) := (1/n) \sum_{i = 1}^n x_i = (1/n) \left\langle \vect{x}_n, \1_n \right\rangle = (1/n) \1_n^\mathsf{T} \vect{x}_n
\end{equation}
Matrix sample mean:
\begin{equation}
\vect{\mu}(\mat{X}_{n \times m}) := (\mu(\mat{X}_{\bullet ,1}), \ldots, \mu(\mat{X}_{\bullet, m}))^\mathsf{T}
\label{eq:matrixmean}
\end{equation}
Sample standard deviation:
\begin{equation}
\label{eq:samplestd}
\begin{aligned}
\sigma(\vect{x}_n) &:= \sqrt{(\left\|\vect{x}_n \right\|_2^2 - n \mu^2(\vect{x}_n))/n} = (\|\1_n \mu(\vect{x}_n) - \vect{x}_n\|_2)/\sqrt{n} = \sqrt{(\vect{x}_{n \mathrm{c}}^\mathsf{T} \vect{x}_{n \mathrm{c}})/n} \\
&= \sqrt{n \vect{x}_n^\mathsf{T} \vect{x}_n - (\1_n^\mathsf{T} \vect{x}_n)^2}/n = \sqrt{(1/n) \vect{x}_n^\mathsf{T} \vect{x}_n - \mu^2(\vect{x}_n)}
\end{aligned}
\end{equation}
Sample variance:
\begin{equation}
\begin{aligned}
\sigma^2(\vect{x}_n) &:= (\left\|\vect{x}_n \right\|_2^2 - n \mu^2(\vect{x}_n))/n = \left\|\1_n \mu(\vect{x}_n) - \vect{x}_n \right\|_2^2/n = \vect{x}_{n \mathrm{c}}^\mathsf{T} \vect{x}_{n \mathrm{c}}/n \\
&= (n \vect{x}_n^\mathsf{T} \vect{x}_n - (\1_n^\mathsf{T} \vect{x}_n)^2)/n^2 = (1/n) \vect{x}_n^\mathsf{T} (\mat{I}_n - (1/n) \1_n \1_n^\mathsf{T}) \vect{x}_n
\end{aligned}
\end{equation}

\textbf{\underline{Centered vectors, matrices and key identities}}

Centered vector (mean-zero):
\begin{equation}
\label{eq:centvec}
\vect{x}_{n \mathrm{c}} := \vect{x}_n - \1_n \mu(\vect{x}_n) = (\mat{I}_n - (1/n) \1_n \1_n^\mathsf{T}) \vect{x}_n = \mat{P}_{n \times n}\vect{x}_n
\end{equation}
where $\mat{P}_{n \times n}$ is symmetric and idempotent:
\begin{equation}
\label{eq:centeringmat}
\mat{P}_{n \times n} := \mat{I}_n - (1/n) \1_n \1_n^\mathsf{T}
\end{equation}
\begin{equation}
\mat{P}_{n \times n}^\mathsf{T}\mat{P}_{n \times n} = \mat{P}_{n \times n}
\end{equation}
\begin{equation}
\mat{P}_{n \times n} \1_n = \mat{I}_n \1_n - (1/n) \1_n \1_n^\mathsf{T} \1_n = \0_n
\end{equation}
Centered matrix:
\begin{equation}
\label{eq:centmat}
\mat{X}_{n \times m, \mathrm{c}} := \mat{P}_{n \times n}\mat{X}_{n \times m} = (\mat{I}_n - (1/n) \1_n \1_n^\mathsf{T}) \mat{X}_{n \times m}
\end{equation}
Mean of centered vector:
\begin{equation}
\label{eq:meancentvec}
\mu(\vect{x}_{n \mathrm{c}}) = 0
\end{equation}
Variance of centered vector:
\begin{equation}
\label{eq:variance_centered_vector}
\sigma^2(\vect{x}_{n \mathrm{c}}) = \sigma^2(\vect{x}_n) = \left\|\vect{x}_{n \mathrm{c}} \right\|_2^2/n
\end{equation}
Euclidean product of centered vectors:
\begin{equation}
\left\langle \vect{x}_{n \mathrm{c}}, \vect{y}_{n \mathrm{c}} \right\rangle = (n \vect{x}_n^\mathsf{T} \vect{y}_n - (\1_n^\mathsf{T} \vect{x}_n) (\1_n^\mathsf{T} \vect{y}_n))/n = \vect{x}_n^\mathsf{T} \vect{y}_n - n \mu(\vect{x}_n) \mu(\vect{y}_n)
\end{equation}
Pythagorean decomposition:
\begin{equation}
\left\|\vect{x}_n \right\|_2^2 = \left\|\1_n \mu(\vect{x}_n) \right\|_2^2 + \left\|\vect{x}_{n \mathrm{c}} \right\|_2^2 = n (\mu^2(\vect{x}_n) + \sigma^2(\vect{x}_n))
\end{equation}
Sum of squared deviations:
\begin{equation}
\label{eq:sum_squared_deviations}
\left\|\vect{x}_{n \mathrm{c}} \right\|_2^2 = \sum_{i = 1}^n (x_i - \mu(\vect{x}_n))^2 = (n \vect{x}_n^\mathsf{T} \vect{x}_n - (\1_n^\mathsf{T} \vect{x}_n)^2)/n
\end{equation}
Norm expansion:
\begin{equation}
\left\|\vect{x}_n + \vect{y}_n \right\|_2^2 = \left\|\vect{x}_n \right\|_2^2 + \left\|\vect{y}_n \right\|_2^2 + 2 \left\langle \vect{x}_n, \vect{y}_n \right\rangle
\end{equation}
Orthogonality to the all-ones vector:
\begin{equation}
\left\langle \vect{x}_{n \mathrm{c}}, \1_n \right\rangle = 0
\end{equation}
Sample Pearson correlation $\rho(\vect{x}_n, \vect{y}_n) \in [-1, 1]$:
\begin{equation}
\label{eq:pearson}
\begin{aligned}
\rho(\vect{x}_n, \vect{y}_n) &:= \frac{\left\langle \vect{x}_{n \mathrm{c}}, \vect{y}_{n \mathrm{c}} \right\rangle}{\left\|\vect{x}_{n \mathrm{c}} \right\|_2 \left\|\vect{y}_{n \mathrm{c}} \right\|_2} = \frac{\vect{x}_{n \mathrm{c}}^\mathsf{T} \vect{y}_{n \mathrm{c}}}{\sqrt{(\vect{x}_{n \mathrm{c}}^\mathsf{T} \vect{x}_{n \mathrm{c}}) (\vect{y}_{n \mathrm{c}}^\mathsf{T} \vect{y}_{n \mathrm{c}})}} \\
&= \frac{\vect{x}_{n \mathrm{c}}^\mathsf{T} \vect{y}_{n \mathrm{c}}}{n \sigma(\vect{x}_n) \sigma(\vect{y}_n)} = \frac{n \vect{x}_n^\mathsf{T} \vect{y}_n - (\1_n^\mathsf{T} \vect{x}_n) (\1_n^\mathsf{T} \vect{y}_n)}{\sqrt{(n \vect{x}_n^\mathsf{T} \vect{x}_n - (\1_n^\mathsf{T} \vect{x}_n)^2) (n \vect{y}_n^\mathsf{T} \vect{y}_n - (\1_n^\mathsf{T} \vect{y}_n)^2)}}
\end{aligned}
\end{equation}

\textbf{\underline{Mean Squared Error ($\mathrm{MSE}$)}}

Using eq.~\eqref{eq:mean_squared_error}, the $\mathrm{MSE}$ can be expressed in terms of Euclidean norms and sample statistics as:
\begin{equation}
\label{eq:MSE_general_1}
\begin{aligned}
\mathrm{MSE}(\vect{x}_n, \vect{y}_n) &= \frac{1}{n} \left\|\vect{x}_n - \vect{y}_n \right\|_2^2 = \frac{1}{n} (\left\|\vect{x}_n \right\|_2^2 + \left\|\vect{y}_n \right\|_2^2 - 2 \left\langle \vect{x}_n, \vect{y}_n \right\rangle) \\
&= \mu^2(\vect{x}_n) + \mu^2(\vect{y}_n) + \sigma^2(\vect{x}_n) + \sigma^2(\vect{y}_n) - \frac{2}{n} \left\langle \vect{x}_n, \vect{y}_n \right\rangle
\end{aligned}
\end{equation}
For $\sigma(\vect{x}_n), \sigma(\vect{y}_n) \in \R \backslash \{0\}$, the $\mathrm{MSE}$ can also be expressed as:
\begin{equation}
\label{eq:MSE_general_2}
\mathrm{MSE}(\vect{x}_n, \vect{y}_n) = (\mu(\vect{x}_n) - \mu(\vect{y}_n))^2 + \sigma^2(\vect{x}_n) + \sigma^2(\vect{y}_n) - 2 \rho(\vect{x}_n, \vect{y}_n) \sigma(\vect{x}_n) \sigma(\vect{y}_n)
\end{equation}
For the special case where $\vect{x}_n = \mu(\vect{y}_n) \1_n$, combining eqs.~\eqref{eq:variance_centered_vector} and \eqref{eq:sum_squared_deviations} gives:
\begin{equation}
\label{eq:MSE_constant_x_n}
\mathrm{MSE}(\mu(\vect{y}_n) \1_n, \vect{y}_n) = \sigma^2(\vect{y}_n)
\end{equation}

\textbf{\underline{Functionals and sample statistics convergence}}

The expectation of a random variable $\underline{x}$ with distribution $F_{\underline{x}}$ is defined as
\begin{equation}
\mathbb{E}_{F_{\underline{x}}}[\underline{x}] := \int x\, \mathrm{d}F_{\underline{x}}(x)
\end{equation}
For a random vector $\underline{\vect{x}}_n$, the component-wise expectation is:
\begin{equation}
\label{eq:component_wise_mean}
\mathbb{E}_{F_{\underline{\vect{x}}_n}}[\underline{\vect{x}}_n] := (\mathbb{E}_{F_{\underline{\vect{x}}_n}}[\underline{x}_1], \ldots, \mathbb{E}_{F_{\underline{\vect{x}}_n}}[\underline{x}_n])^\mathsf{T}
\end{equation}
For a random matrix, the expectation is defined component-wise and is itself a matrix of the same dimensions.

The variance of a random variable $\underline{x}$ with distribution $F_{\underline{x}}$ is:
\begin{equation}
\label{eq:variance}
\mathrm{Var}_{F_{\underline{x}}}[\underline{x}] := \mathbb{E}_{F_{\underline{x}}}[(\underline{x} - \mathbb{E}_{F_{\underline{x}}}[\underline{x}])^2]
\end{equation}
The variance-covariance matrix of a random vector $\underline{\vect{x}}_n$ is then defined as:
\begin{equation}
\label{eq:variance_covariance_matrix}
\mathrm{Var}_{F_{\underline{\vect{x}}_n}}[\underline{\vect{x}}_n] := \mathbb{E}_{F_{\underline{\vect{x}}_n}}[(\underline{\vect{x}}_n - \mathbb{E}_{F_{\underline{\vect{x}}_n}}[\underline{\vect{x}}_n]) (\underline{\vect{x}}_n - \mathbb{E}_{F_{\underline{\vect{x}}_n}}[\underline{\vect{x}}_n])^\mathsf{T}]
\end{equation}
The covariance between random variables $\underline{x}$ and $\underline{y}$ is:
\begin{equation}
\label{eq:covariance}
\mathrm{Cov}_{F_{\underline{x}, \underline{y}}}(\underline{x}, \underline{y}) = \mathbb{E}_{F_{\underline{x}, \underline{y}}}[(\underline{x} - \mathbb{E}_{F_{\underline{x}}}[\underline{x}]) (\underline{y} - \mathbb{E}_{F_{\underline{y}}}[\underline{y}])]
\end{equation}
The correlation between random variables $\underline{x}$ and $\underline{y}$ is:
\begin{equation}
\mathrm{Corr}_{F_{\underline{x}, \underline{y}}}(\underline{x}, \underline{y}) = \frac{\mathrm{Cov}_{F_{\underline{x}, \underline{y}}}(\underline{x}, \underline{y})}{\sqrt{\mathrm{Var}_{F_{\underline{x}}}[\underline{x}] \mathrm{Var}_{F_{\underline{y}}}[\underline{y}]}}
\end{equation}
The vector of covariances between each component of $\underline{\vect{x}}_n$ and a scalar random variable $\underline{y}$ is
\begin{equation}
\label{eq:covariances_vector}
\mathrm{Cov}_{F_{\underline{\vect{x}}_n, \underline{y}}}(\underline{\vect{x}}_n, \underline{y}) := (\mathrm{Cov}_{F_{\underline{\vect{x}}_n, \underline{y}}}(\underline{x}_1, \underline{y}), \ldots, \mathrm{Cov}_{F_{\underline{\vect{x}}_n, \underline{y}}}(\underline{x}_n, \underline{y}))^\mathsf{T}
\end{equation}
Analogously, the vector of Pearson correlations is:
\begin{equation}
\mathrm{Corr}_{F_{\underline{\vect{x}}_n, \underline{y}}}(\underline{\vect{x}}_n, \underline{y}) := (\mathrm{Corr}_{F_{\underline{\vect{x}}_n, \underline{y}}}(\underline{x}_1, \underline{y}), \ldots, \mathrm{Corr}_{F_{\underline{\vect{x}}_n, \underline{y}}}(\underline{x}_n, \underline{y}))^\mathsf{T}
\end{equation}
The sample statistics that appear in the estimators; specifically the sample means, variances and covariances, converge almost surely to their population counterparts by the strong law of large numbers. Under the assumption of i.i.d. variables, the relevant sample statistics converge as $n \to \infty$ as follows:

The sample mean $\mu(\underline{\vect{x}}_n)$, defined in eq.~\eqref{eq:samplemean}, converges almost surely to the expectation (by the strong law of large numbers):
\begin{equation}
\label{eq:sample_mean_convergence}
\mu(\underline{\vect{x}}_n) \overset{\mathrm{a.s.}}{\longrightarrow} \mathbb{E}_{F_{\underline{x}}}[\underline{x}]
\end{equation}
The following convergence results follow immediately from the strong law of large numbers and the continuous mapping theorem.

The sample standard deviation $\sigma(\underline{\vect{x}}_n)$ (eq.~\eqref{eq:samplestd}) converges almost surely to the square root of the variance:
\begin{equation}
\label{eq:sample_standard_deviation_convergence}
\sigma(\underline{\vect{x}}_n) \overset{\mathrm{a.s.}}{\longrightarrow} \sqrt{\mathrm{Var}_{F_{\underline{x}}}[\underline{x}]}
\end{equation}
The sample Pearson correlation $\rho(\underline{\vect{x}}_n, \underline{\vect{y}}_n)$ (eq.~\eqref{eq:pearson}) converges almost surely to the correlation:
\begin{equation}
\label{eq:sample_Pearson_correlation_convergence}
\rho(\underline{\vect{x}}_n, \underline{\vect{y}}_n) \overset{\mathrm{a.s.}}{\longrightarrow} \mathrm{Corr}_{F_{\underline{x}, \underline{y}}}(\underline{x}, \underline{y})
\end{equation}
The vector of sample means of the predictor matrix, $\vect{\mu}(\underline{\mat{X}}_{n \times p})$ (eq.~\eqref{eq:matrixmean}), converges almost surely to the element-wise expectations:
\begin{equation}
\label{eq:matrix_sample_mean_convergence}
\vect{\mu}(\underline{\mat{X}}_{n \times p}) \overset{\mathrm{a.s.}}{\longrightarrow} \mathbb{E}_{F_{\underline{\vect{x}}_p}}[\underline{\vect{x}}_p]
\end{equation}
The sample covariance matrix $\underline{\mat{S}}_{p \times p}$ (eq.~\eqref{eq:covariance_matrix}) converges almost surely to the variance-covariance matrix:
\begin{equation}
\label{eq:sample_covariance_matrix_convergence}
\underline{\mat{S}}_{p \times p} \overset{\mathrm{a.s.}}{\longrightarrow} \mathrm{Var}_{F_{\underline{\vect{x}}_p}}[\underline{\vect{x}}_p]
\end{equation}
The sample cross-covariance vector $\underline{\vect{s}}_p$ (eq.~\eqref{eq:cross_covariance_vector}) converges almost surely to the cross-covariance vector:
\begin{equation}
\label{eq:sample_cross_covariance_vector_convergence}
\underline{\vect{s}}_p \overset{\mathrm{a.s.}}{\longrightarrow} \mathrm{Cov}_{F_{\underline{\vect{x}}_p, \underline{y}}}(\underline{\vect{x}}_p, \underline{y})
\end{equation}

\textbf{\underline{Inequalities}}

Cauchy-Schwarz inequality \cite[p.~36]{gentle2024}:
\begin{equation}
\label{eq:Cauchy_Schwarz_inequality}
\left|\left\langle \vect{x}_n, \vect{y}_n \right\rangle \right| \leq \left\|\vect{x}_n \right\|_2 \left\|\vect{y}_n \right\|_2
\end{equation}
with equality \cite[p.~34]{gentle2024}:
\begin{equation}
\left|\left\langle \vect{x}_n, \vect{y}_n \right\rangle \right| = \left\|\vect{x}_n \right\|_2 \left\|\vect{y}_n \right\|_2, \ \text{iff}\ \vect{x}_n = \0_n, \ \text{or}\ \vect{y}_n = \0_n,\ \text{or}\ \vect{x}_n = a \vect{y}_n
\end{equation}
Because $\left\langle a \vect{y}_n, \vect{y}_n \right\rangle = a \left\langle \vect{y}_n,\vect{y}_n \right\rangle$ and $\left\langle a \vect{y}_n,\vect{y}_n \right\rangle \leq \left|\left\langle a \vect{y}_n,\vect{y}_n \right\rangle \right|$:
\begin{equation}
\label{eq:Cauchy_Schwarz_equality_condition}
\left\langle \vect{x}_n, \vect{y}_n \right\rangle = \left\|\vect{x}_n \right\|_2 \left\| \vect{y}_n \right\|_2, \ \text{iff}\ \vect{x}_n = \0_n, \ \text{or}\ \vect{y}_n = \0_n, \ \text{or}\ \vect{x}_n = a \vect{y}_n, a > 0
\end{equation}

\textbf{\underline{The normal (Gaussian) distribution}}

Let $\underline{y}$ be a Gaussian random variable $\underline{y} \sim N(\mu, \sigma)$. Its probability density function (PDF) is:
\begin{equation}
\label{eq:normal_pdf}
f_N(y; \mu, \sigma) = \frac{1}{\sigma \sqrt{2 \pi}} \exp(-\frac{(y - \mu)^2}{2 \sigma^2}), \mu \in \R, \sigma > 0
\end{equation}
\section{Proofs}
\label{app:proofs}
This appendix supplies the proofs for all theoretical results presented in the manuscript. \hyperref[proof:b1]{Proof~B.1} establishes that the Kling-Gupta loss attains its global minimum uniquely at $\vect{z}_n = \vect{y}_n$ (see Section~\ref{sec:kling_gupta_loss_sec}). \hyperref[proof:b2]{Proof~B.2} gives the parameter estimates for the linear model with multiple predictors estimated by minimizing the Kling-Gupta loss (Section~\ref{sec:kling_gupta_linear_regression_multiple_predictors}). \hyperref[proof:b3]{Proof~B.3} computes the training-set performance metrics appearing in Table~\ref{tab:training-performance} (Section~\ref{sec:comparative_training_data}). \hyperref[proof:b4]{Proof~B.4} gives the asymptotic limits of these metrics on an independent test set (Section~\ref{sec:comparative_infinite_training_test_data}). \hyperref[proof:b5]{Proof~B.5} treats the intercept estimate when the regression coefficients are fixed (Section~\ref{sec:fixed_regression_coefficient}), while \hyperref[proof:b6]{Proof~B.6} treats the complementary case where the intercept is fixed and the slope is estimated, restricted to the single-predictor model (Section~\ref{sec:fixed_intercept_kg_linear_regression}). All proofs build on the sample statistics, vector and matrix relationships summarized in Appendix~\ref{app:notation}.
\prooftitle{Proof B.1}{(Uniqueness of the global minimum of the Kling-Gupta loss, Section~\ref{sec:kling_gupta_loss_sec})}{proof:b1}
We prove that the Kling-Gupta loss $L_\mathrm{KG}(\vect{z}_n, \vect{y}_n)$ attains a unique global minimum at $\vect{z}_n = \vect{y}_n$. Because $L_\mathrm{KG}$ is a sum of three non-negative terms, it follows immediately that $L_\mathrm{KG}(\vect{z}_n, \vect{y}_n) \geq 0$ and that the lower bound is attained when $\vect{z}_n = \vect{y}_n$. To establish uniqueness, we show that $L_\mathrm{KG}(\vect{z}_n, \vect{y}_n) = 0$ implies $\vect{z}_n = \vect{y}_n$.

The loss equals zero if and only if each of its three component terms equals zero. We assume $\mu(\vect{y}_n), \sigma(\vect{y}_n), \sigma(\vect{z}_n) \in \R \backslash \{0\}$, which in turn (from eq.~\eqref{eq:variance_centered_vector}) ensures that the centered vectors satisfy $\vect{y}_{n \mathrm{c}} \neq \0_n$ and $\vect{z}_{n \mathrm{c}} \neq \0_n$.

First, setting the correlation term to zero and applying the Cauchy-Schwarz equality condition (eq.~\eqref{eq:Cauchy_Schwarz_equality_condition}) gives:
\begin{equation}
\label{eq:cor_equality_KG_minima_unique}
\rho(\vect{z}_n, \vect{y}_n) = 1 \Leftrightarrow \left\langle \vect{z}_{n \mathrm{c}}, \vect{y}_{n \mathrm{c}} \right\rangle = \left\|\vect{z}_{n \mathrm{c}} \right\|_2 \left\|\vect{y}_{n \mathrm{c}} \right\|_2 \Leftrightarrow \vect{z}_{n \mathrm{c}} = a \vect{y}_{n \mathrm{c}},a > 0
\end{equation}
Second, setting the variability term to zero and using eq.~\eqref{eq:variance_centered_vector} implies:
\begin{equation}
\label{eq:sigma_equality_KG_minima_unique}
\sigma(\vect{z}_n) = \sigma(\vect{y}_n) \Leftrightarrow \left\|\vect{z}_{n \mathrm{c}} \right\|_2 = \left\|\vect{y}_{n \mathrm{c}} \right\|_2 \Leftrightarrow |a| \left\|\vect{y}_{n \mathrm{c}} \right\|_2 = \left\|\vect{y}_{n \mathrm{c}} \right\|_2 \Leftrightarrow a = 1
\end{equation}
The final step in eq.~\eqref{eq:sigma_equality_KG_minima_unique} follows from $a > 0$, which was established in eq.~\eqref{eq:cor_equality_KG_minima_unique}; together with $\vect{z}_{n \mathrm{c}} = a \vect{y}_{n \mathrm{c}}$ from eq.~\eqref{eq:cor_equality_KG_minima_unique}, this gives $\vect{z}_{n \mathrm{c}} = \vect{y}_{n \mathrm{c}}$. Third, setting the bias term to zero together with eq.~\eqref{eq:centvec} leads to:
\begin{equation}
\mu(\vect{z}_n) = \mu(\vect{y}_n) \Leftrightarrow \1_n \mu(\vect{z}_n) = \1_n \mu(\vect{y}_n) \Leftrightarrow \vect{z}_n - \vect{z}_{n \mathrm{c}} = \vect{y}_n - \vect{y}_{n \mathrm{c}} \Leftrightarrow \vect{z}_n = \vect{y}_n
\end{equation}
Thus, $L_\mathrm{KG}(\vect{z}_n, \vect{y}_n) = 0$ if and only if $\vect{z}_n = \vect{y}_n$, which completes the proof. \qed

\prooftitle{Proof B.2}{(Parameter estimates for the multiple-predictor linear model estimated with the Kling-Gupta loss, Section~\ref{sec:kling_gupta_linear_regression_multiple_predictors})}{proof:b2}

We determine the parameter estimates for the linear model $z = \vect{a}_p^\mathsf{T} \vect{x}_p + b_p$ specified in eq.~\eqref{eq:linear_model} when it is trained by minimizing the Kling-Gupta loss function $L_\mathrm{KG}$ defined in eq.~\eqref{eq:Kling_Gupta_loss}, given observations $\vect{y}_n$ of the response variable and $\mat{X}_{n \times p}$ of the predictor variables, and following the notational conventions of Section~\ref{sec:linear_model_multiple_predictors}. We prove that the model parameter estimates are given by eqs.~\eqref{eq:kg_slope_multi}, \eqref{eq:kg_ellipsoid_condition} and \eqref{eq:kg_intercept_multi}.

The standard conditions required for the definition of $L_\mathrm{KG}$ are assumed: $\mu(\vect{y}_n) \neq 0$, $\sigma(\vect{y}_n) \neq 0$ and $\sigma(\vect{z}_n) \neq 0$ (eq.~\eqref{eq:Kling_Gupta_loss}). Furthermore, as indicated in Section~\ref{sec:linear_model_multiple_predictors}, the sample size satisfies $n \geq p$, the predictor matrix $\mat{X}_{n \times p}$ has full column rank, and the all-ones vector $\1_n$ does not lie in its column space. These conditions ensure that the sample covariance matrix $\mat{S}_{p \times p}$ (defined by eq.~\eqref{eq:covariance_matrix}) is positive definite (Section~\ref{sec:linear_model_multiple_predictors}) and that the OLS regression parameter estimate $\widehat{\vect{a}}_{p, \mathrm{OLS}}$ (given by eq.~\eqref{eq:ols_slope_multi}) is uniquely defined (Section~\ref{sec:ols-linear-regression}). The restriction $\sigma(\vect{z}_n) \neq 0$, together with eq.~\eqref{eq:sigma_predictions} and the positive definiteness of $\mat{S}_{p \times p}$, implies that any admissible estimate must satisfy $\widehat{\vect{a}}_{p, \mathrm{KG}} \neq \0_p$.

The objective is to minimize $L_\mathrm{KG}$ with respect to the parameter vector $\vect{\theta}_{1 \times (p + 1)} = (a_1, \ldots, a_p, b_p)$, as formulated in eq.~\eqref{eq:kg_estimator_multi_def}. Substituting the linear model predictions $\vect{z}_n$ from eq.~\eqref{eq:lm_training_predictions} into $L_\mathrm{KG}$, we observe that the intercept $b_p$ appears exclusively in the bias term, whereas the variability and correlation terms depend only on the slope vector $\vect{a}_p$ (see eqs.~\eqref{eq:sigma_predictions} and \eqref{eq:correlation_predictions_multi}). Therefore, the bias term can be minimized independently of the other components. It attains its unique minimum of zero when the sample means of the predictions and the observations are equal, i.e., $\mu(\vect{z}_n) = \mu(\vect{y}_n)$. Using the expression for the predictive mean of the linear model from eq.~\eqref{eq:mean_predictions_multi}, we solve for the intercept, which gives:
\begin{equation}
\label{eq:KG_initial_intercept_estimate}
\widehat{b}_{p, \mathrm{KG}} = \mu(\vect{y}_n) - (\vect{\mu}(\mat{X}_{n \times p}))^\mathsf{T} \vect{a}_p
\end{equation}
With this choice of intercept, the bias term vanishes. For any $\vect{a}_p \neq \0_p$, $L_\mathrm{KG}$ reduces to:
\begin{equation}
\label{eq:L_ap_function}
L(\vect{a}_p) = (1 - \frac{\sigma(\vect{z}_n)}{\sigma(\vect{y}_n)})^2 + (1 - \rho(\vect{z}_n, \vect{y}_n))^2, \vect{a}_p \neq \0_p
\end{equation}
The function $L(\vect{a}_p)$ depends on $\vect{a}_p$ through $\sigma(\vect{z}_n)$ (eq.~\eqref{eq:sigma_predictions}) and $\rho(\vect{z}_n, \vect{y}_n)$ (eq.~\eqref{eq:correlation_predictions_multi}); substituting these expressions gives:
\begin{equation}
L(\vect{a}_p) = (1 - \frac{\sqrt{\vect{a}_p^\mathsf{T} \mat{S}_{p \times p} \vect{a}_p}}{\sigma(\vect{y}_n)})^2 + (1 - \frac{\vect{a}_p^\mathsf{T} \vect{s}_p}{\sqrt{\vect{a}_p^\mathsf{T} \mat{S}_{p \times p}\vect{a}_p} \sigma(\vect{y}_n)})^2, \vect{a}_p \neq \0_p
\end{equation}
First, we demonstrate that:
\begin{equation}
\label{eq:limit_infinity_L_ap}
\lim_{\|\vect{a}_p\|_2 \to \infty} L(\vect{a}_p) = \infty
\end{equation}
This follows because (i) $-1 \leq \rho(\vect{z}_n, \vect{y}_n) \leq 1$ and (ii) as $\|\vect{a}_p\|_2 \to \infty$, the positive definiteness of $\mat{S}_{p \times p}$ together with eq.~\eqref{eq:sigma_predictions} implies $\sigma(\vect{z}_n) = (\vect{a}_p^\mathsf{T} \mat{S}_{p \times p} \vect{a}_p)^{1/2} \to \infty$.

We proceed by distinguishing two mutually exclusive cases based on the sample cross-covariance vector $\vect{s}_p$ defined in eq.~\eqref{eq:cross_covariance_vector}:

\noindent\makebox[1cm][l]{(i)} \textbf{Case \#1} ($\vect{s}_p \neq \0_p$): We first determine a theoretical lower bound for $L(\vect{a}_p)$. Eq.~\eqref{eq:limit_infinity_L_ap} shows that this lower bound is not attained as $\|\vect{a}_p\|_2 \to \infty$. We examine the correlation $\rho(\vect{z}_n, \vect{y}_n)$, which, from eqs.~\eqref{eq:correlation_predictions_multi} and \eqref{eq:ols_cross_covariance}, can be expressed as:
\begin{equation}
\label{eq:correlation_for_proof_B2}
\rho(\vect{z}_n, \vect{y}_n) = \frac{\vect{a}_p^\mathsf{T} \mat{S}_{p \times p} \widehat{\vect{a}}_{p, \mathrm{OLS}}}{\sqrt{\vect{a}_p^\mathsf{T} \mat{S}_{p \times p} \vect{a}_p} \sigma(\vect{y}_n)}, \vect{a}_p \neq \0_p
\end{equation}
From $\vect{s}_p \neq \0_p$ and eq.~\eqref{eq:ols_cross_covariance} it follows that $\widehat{\vect{a}}_{p, \mathrm{OLS}} \neq \0_p$; eq.~\eqref{eq:ols_sigma_expansion} then implies $\sigma_\mathrm{OLS} = (\widehat{\vect{a}}_{p, \mathrm{OLS}}^\mathsf{T} \mat{S}_{p \times p} \widehat{\vect{a}}_{p, \mathrm{OLS}})^{1/2} > 0$, and together with eq.~\eqref{eq:ols_correlation_deriv} this leads to $\rho_\mathrm{OLS} = \frac{\sigma_\mathrm{OLS}}{\sigma(\vect{y}_n)} > 0$. The matrix $\mat{S}_{p \times p}$ is symmetric positive definite and therefore invertible, so it possesses a unique symmetric positive definite square root $\mat{S}_{p \times p}^{1/2}$ \cite[p.~187, 244]{gentle2024}. Applying the Cauchy-Schwarz inequality (eq.~\eqref{eq:Cauchy_Schwarz_inequality}) to the inner product of $\mat{S}_{p \times p}^{1/2} \vect{a}_p$ and $\mat{S}_{p \times p}^{1/2} \widehat{\vect{a}}_{p, \mathrm{OLS}}$ gives:
\begin{equation}
\label{eq:CS_inequality_proof_KG_estimates}
\begin{aligned}
\vect{a}_p^\mathsf{T} \mat{S}_{p \times p} \widehat{\vect{a}}_{p, \mathrm{OLS}}
&= \vect{a}_p^\mathsf{T} \mat{S}_{p \times p}^{1/2} \mat{S}_{p \times p}^{1/2} \widehat{\vect{a}}_{p, \mathrm{OLS}} \\
&= \left\langle \mat{S}_{p \times p}^{1/2} \vect{a}_p, \mat{S}_{p \times p}^{1/2} \widehat{\vect{a}}_{p, \mathrm{OLS}} \right\rangle \\
&\leq \left|\left\langle \mat{S}_{p \times p}^{1/2} \vect{a}_p, \mat{S}_{p \times p}^{1/2} \widehat{\vect{a}}_{p, \mathrm{OLS}} \right\rangle \right| \\
&\leq \left\|\mat{S}_{p \times p}^{1/2} \vect{a}_p \right\|_{2} \left\|\mat{S}_{p \times p}^{1/2} \widehat{\vect{a}}_{p, \mathrm{OLS}} \right\|_2 \\
&= (\vect{a}_p^\mathsf{T} \mat{S}_{p \times p} \vect{a}_p)^{1/2} (\widehat{\vect{a}}_{p, \mathrm{OLS}}^\mathsf{T} \mat{S}_{p \times p} \widehat{\vect{a}}_{p, \mathrm{OLS}})^{1/2}, \vect{a}_p \neq \0_p
\end{aligned}
\end{equation}
Using this bound in eq.~\eqref{eq:correlation_for_proof_B2} together with the definitions of $\sigma_\mathrm{OLS}$ from eq.~\eqref{eq:ols_sigma_expansion} and $\rho_\mathrm{OLS}$ from eq.~\eqref{eq:ols_correlation_deriv} gives:
\begin{equation}
\rho(\vect{z}_n, \vect{y}_n) \leq \frac{\sqrt{\widehat{\vect{a}}_{p, \mathrm{OLS}}^\mathsf{T} \mat{S}_{p \times p} \widehat{\vect{a}}_{p, \mathrm{OLS}}}}{\sigma(\vect{y}_n)} = \frac{\sigma_\mathrm{OLS}}{\sigma(\vect{y}_n)} = \rho_\mathrm{OLS} \leq 1
\end{equation}
Thus:
\begin{equation}
1 - \rho(\vect{z}_n, \vect{y}_n) \geq 1 - \rho_\mathrm{OLS} \geq 0
\end{equation}
It follows that $L(\vect{a}_p)$ satisfies:
\begin{equation}
L(\vect{a}_p) = (1 - \frac{\sigma(\vect{z}_n)}{\sigma(\vect{y}_n)})^2 + (1 - \rho(\vect{z}_n, \vect{y}_n))^2 \geq (1 - \rho(\vect{z}_n, \vect{y}_n))^2 \geq (1 - \rho_\mathrm{OLS})^2 \geq 0
\end{equation}
Equality $L(\vect{a}_p) = (1 - \rho_\mathrm{OLS})^2$ is attained if and only if (i) the Cauchy-Schwarz inequality in eq.~\eqref{eq:CS_inequality_proof_KG_estimates} is an equality, and (ii) the variability term vanishes, i.e., $\sigma(\vect{z}_n) = \sigma(\vect{y}_n)$.

The Cauchy-Schwarz inequality in eq.~\eqref{eq:CS_inequality_proof_KG_estimates} becomes an equality if and only if (see eq.~\eqref{eq:Cauchy_Schwarz_equality_condition} and \cite[p.~258]{seber2007}, together with the condition $\vect{a}_p \neq \0_p$):
\begin{equation}
\mat{S}_{p \times p}^{1/2} \vect{a}_p = c \mat{S}_{p \times p}^{1/2} \widehat{\vect{a}}_{p, \mathrm{OLS}}, c > 0
\end{equation}
The invertibility of $\mat{S}_{p \times p}^{1/2}$ implies:
\begin{equation}
\vect{a}_p = c \widehat{\vect{a}}_{p, \mathrm{OLS}}, c > 0
\end{equation}
From eqs.~\eqref{eq:sigma_predictions} and \eqref{eq:ols_sigma_expansion} we have:
\begin{equation}
\sigma(\vect{z}_n) = \sqrt{c^2 \widehat{\vect{a}}_{p, \mathrm{OLS}}^\mathsf{T} \mat{S}_{p \times p} \widehat{\vect{a}}_{p, \mathrm{OLS}}} = c \sigma_\mathrm{OLS}
\end{equation}
The condition $\sigma(\vect{z}_n) = \sigma(\vect{y}_n)$ then requires:
\begin{equation}
c = \frac{\sigma(\vect{y}_n)}{\sigma_\mathrm{OLS}}
\end{equation}
Therefore, the unique global minimizer is:
\begin{equation}
\widehat{\vect{a}}_{p, \mathrm{KG}} = \frac{\sigma(\vect{y}_n)}{\sigma_\mathrm{OLS}} \widehat{\vect{a}}_{p, \mathrm{OLS}}
\end{equation}
The corresponding intercept estimate, from eq.~\eqref{eq:KG_initial_intercept_estimate}, is:
\begin{equation}
\widehat{b}_{p, \mathrm{KG}} = \mu(\vect{y}_n) - (\vect{\mu}(\mat{X}_{n \times p}))^\mathsf{T} \widehat{\vect{a}}_{p, \mathrm{KG}}
\end{equation}
At this minimum the loss equals:
\begin{equation}
L_\mathrm{KG}(\mat{X}_{n \times p} \widehat{\vect{a}}_{p, \mathrm{KG}} + \widehat{b}_{p, \mathrm{KG}} \1_n, \vect{y}_n) = (1 - \rho_\mathrm{OLS})^2
\end{equation}
To confirm that this point is a global minimum and is not excluded by the discontinuity at $\vect{a}_p = \0_p$, we examine the behavior along sequences that converge to the origin. For any sequence $\vect{a}_p^{(k)} \to \0_p$, the prediction variance satisfies $\sigma(\vect{z}_n^{(k)}) = \sqrt{(\vect{a}_p^{(k)})^\mathsf{T} \mat{S}_{p \times p} \vect{a}_p^{(k)}} \to 0$, so the variability term $(1 - \frac{\sigma(\vect{z}_n^{(k)})}{\sigma(\vect{y}_n)})^2 \to 1$. Because the correlation term is non-negative, the loss satisfies $L(\vect{a}_p) \geq (1 - \frac{\sigma(\vect{z}_n)}{\sigma(\vect{y}_n)})^2$. Therefore, the limit inferior is at least unity, as the variability term converges to unity while the correlation term remains non-negative:
\begin{equation}
\liminf_{\vect{a}_p \to \0_p} L(\vect{a}_p) \geq 1
\end{equation}
From $\rho_\mathrm{OLS} > 0$, we have $L_\mathrm{KG}(\mat{X}_{n \times p} \widehat{\vect{a}}_{p, \mathrm{KG}} + \widehat{b}_{p, \mathrm{KG}} \1_n, \vect{y}_n) = (1 - \rho_\mathrm{OLS})^2 < 1$. Thus the value at the global minimum is strictly less than the infimum near the discontinuity, confirming that the minimum is indeed globally optimal.

\noindent\makebox[1cm][l]{(ii)} \textbf{Case \#2} ($\vect{s}_p = \0_p$): From eqs.~\eqref{eq:ols_sigma_expansion} and \eqref{eq:ols_cross_covariance} we have $\sigma_\mathrm{OLS} = 0$ and $\widehat{\vect{a}}_{p, \mathrm{OLS}} = \0_p$. Eq.~\eqref{eq:correlation_predictions_multi} shows that for every non-zero slope the correlation vanishes identically, $\rho(\vect{z}_n, \vect{y}_n) = 0$ for all $\vect{a}_p \neq \0_p$. The loss function $L(\vect{a}_p)$ in eq.~\eqref{eq:L_ap_function} therefore simplifies to:
\begin{equation}
L(\vect{a}_p) = (1 - \frac{\sigma(\vect{z}_n)}{\sigma(\vect{y}_n)})^2 + 1
\end{equation}
Because eq.~\eqref{eq:limit_infinity_L_ap} gives $L(\vect{a}_p)\to\infty$ as $\|\vect{a}_p\|_2\to\infty$, a finite global minimum exists. The minimum value $1$ is attained when the variability term is zero, i.e., when $\sigma(\vect{z}_n) = \sigma(\vect{y}_n)$. Using eq.~\eqref{eq:sigma_predictions} this condition becomes:
\begin{equation}
\label{eq:KG_ellipsoid}
\widehat{\vect{a}}_{p, \mathrm{KG}}^\mathsf{T} \mat{S}_{p \times p} \widehat{\vect{a}}_{p, \mathrm{KG}} = \sigma^2(\vect{y}_n)
\end{equation}
Any vector $\widehat{\vect{a}}_{p, \mathrm{KG}}$ satisfying this ellipsoid constraint gives $L_\mathrm{KG} = 1$, which is the global minimum over the domain $\vect{a}_p \neq \0_p$.

The intercept estimate continues to be given by eq.~\eqref{eq:KG_initial_intercept_estimate}:
\begin{equation}
\widehat{b}_{p, \mathrm{KG}} = \mu(\vect{y}_n) - (\vect{\mu}(\mat{X}_{n \times p}))^\mathsf{T} \widehat{\vect{a}}_{p, \mathrm{KG}}
\end{equation}
Because $\widehat{\vect{a}}_{p, \mathrm{KG}} \neq \0_p$ (a consequence of eq.~\eqref{eq:KG_ellipsoid} and the positive definiteness of $\mat{S}_{p \times p}$), the value of the loss at any such global minimizer is:
\begin{equation}
L_\mathrm{KG}(\mat{X}_{n \times p} \widehat{\vect{a}}_{p, \mathrm{KG}} + \widehat{b}_{p, \mathrm{KG}} \1_n, \vect{y}_n) = 1
\end{equation}
Thus the lower bound $1$ is a global minimum: for every $\vect{a}_p \neq \0_p$ we have $L(\vect{a}_p) > 1$ unless $\vect{a}_p$ satisfies the ellipsoid condition~\eqref{eq:KG_ellipsoid}.

Finally, the limit inferior of $L(\vect{a}_p)$ as $\vect{a}_p \to \0_p$ equals $2$, for the following reasons:

\noindent\makebox[1cm][l]{(ii.1)} The variability term converges to unity.

\noindent\makebox[1cm][l]{(ii.2)} The correlation term also converges to unity, because $\rho(\vect{z}_n, \vect{y}_n)$ is discontinuous at $\vect{a}_p = \0_p$, and, by eq.~\eqref{eq:correlation_predictions_multi} its limit as $\vect{a}_p \to \0_p$ is zero. This is established because $\vect{s}_p = \0_p$ and the Euclidean norm of $\frac{\vect{a}_p^\mathsf{T}}{\sqrt{\vect{a}_p^\mathsf{T} \mat{S}_{p \times p} \vect{a}_p}}$ (which appears in eq.~\eqref{eq:correlation_predictions_multi}) is bounded. From the properties of the Rayleigh quotient $\frac{\vect{a}_p^\mathsf{T} \mat{S}_{p \times p} \vect{a}_p}{\|\vect{a}_p\|_2^2}$, there exists a constant $\lambda_{\min} > 0$ such that $\frac{\vect{a}_p^\mathsf{T} \mat{S}_{p \times p} \vect{a}_p}{\|\vect{a}_p\|_2^2} \geq \lambda_{\min}$. Therefore, $\|\frac{\vect{a}_p^\mathsf{T}}{\sqrt{\vect{a}_p^\mathsf{T} \mat{S}_{p \times p} \vect{a}_p}}\|_2 \leq \frac{1}{\sqrt{\lambda_{\min}}}$. \qed
\prooftitle{Proof B.3}{(Training set performance, Section~\ref{sec:comparative_training_data})}{proof:b3}
We evaluate the training set performance metrics reported in Table~\ref{tab:training-performance}.

\noindent\makebox[1cm][l]{(i)} \textbf{OLS predictions}: For the multiple-predictor OLS linear regression model, the training set predictions have the following properties: $\mu(\vect{z}_{n, \mathrm{OLS}}) = \mu(\vect{y}_n)$ (eq.~\eqref{eq:ols_mean_property}), $\sigma(\vect{z}_{n, \mathrm{OLS}}) = \sigma_\mathrm{OLS}$ (eq.~\eqref{eq:ols_sigma_def}), and, for $\vect{s}_p \neq \0_p$, $\rho(\vect{z}_{n, \mathrm{OLS}}, \vect{y}_n) = \rho_\mathrm{OLS}$ (eq.~\eqref{eq:ols_rho_def}). To treat the zero cross-covariance case $\vect{s}_p = \0_p$ uniformly, we use the extended correlation $\rho_{\mathrm{OLS}_{(0)}}$ defined in eq.~\eqref{eq:rho_ols_star_def}, which equals $\rho_\mathrm{OLS}$ when $\vect{s}_p \neq \0_p$ and is $0$ when $\vect{s}_p = \0_p$.

For $\vect{s}_p \neq \0_p$, the $\mathrm{MSE}$ follows from eq.~\eqref{eq:MSE_general_2} and the properties above:
\begin{equation}
\mathrm{MSE}(\vect{z}_{n, \mathrm{OLS}}, \vect{y}_n) = \sigma_\mathrm{OLS}^2 + \sigma^2(\vect{y}_n) - 2 \sigma_\mathrm{OLS} \sigma(\vect{y}_n) \rho_\mathrm{OLS}
\end{equation}
Substituting the identity $\sigma_\mathrm{OLS} = \rho_\mathrm{OLS} \sigma(\vect{y}_n)$ from eq.~\eqref{eq:ols_correlation_deriv} gives:
\begin{equation}
\label{eq:OSL_training_MSE_preliminary}
\mathrm{MSE}(\vect{z}_{n, \mathrm{OLS}}, \vect{y}_n) = \sigma^2(\vect{y}_n) (1 - \rho_\mathrm{OLS}^2)
\end{equation}
For $\vect{s}_p = \0_p$, Section~\ref{sec:ols-linear-regression} states $\vect{z}_{n, \mathrm{OLS}} = \mu(\vect{y}_n) \1_n$, and eq.~\eqref{eq:MSE_constant_x_n} then gives $\mathrm{MSE}(\mu(\vect{y}_n) \1_n, \vect{y}_n) = \sigma^2(\vect{y}_n)$. Therefore, the result of eq.~\eqref{eq:OSL_training_MSE_preliminary} can be unified for both cases of $\vect{s}_p$ as:
\begin{equation}
\label{eq:OSL_training_MSE}
\mathrm{MSE}(\vect{z}_{n, \mathrm{OLS}}, \vect{y}_n) = \sigma^2(\vect{y}_n) (1 - \rho_{\mathrm{OLS}_{(0)}}^2)
\end{equation}
From the definition of $\mathrm{NSE}$ (eq.~\eqref{eq:Nash_Sutcliffe_Efficiency}) together with eqs.~\eqref{eq:MSE_constant_x_n} and \eqref{eq:OSL_training_MSE}, it follows that:
\begin{equation}
\mathrm{NSE}(\vect{z}_{n, \mathrm{OLS}}, \vect{y}_n) = 1 - \frac{\sigma^2(\vect{y}_n) (1 - \rho_{\mathrm{OLS}_{(0)}}^2)}{\sigma^2(\vect{y}_n)} = \rho_{\mathrm{OLS}_{(0)}}^2
\end{equation}
For the Kling-Gupta loss of the OLS model we restrict attention to the case $\vect{s}_p \neq \0_p$ (the loss is undefined when $\vect{s}_p = \0_p$ because eqs.~\eqref{eq:ols_sigma_def}, \eqref{eq:ols_sigma_expansion} and \eqref{eq:ols_cross_covariance} together imply $\sigma(\vect{z}_{n, \mathrm{OLS}}) = 0$). Substituting the properties of the OLS predictions into eq.~\eqref{eq:Kling_Gupta_loss} gives
\begin{equation}
L_\mathrm{KG}(\vect{z}_{n, \mathrm{OLS}}, \vect{y}_n) = (1 - \frac{\mu(\vect{z}_{n, \mathrm{OLS}})}{\mu(\vect{y}_n)})^2 + (1 - \frac{\sigma_{\mathrm{OLS}}}{\sigma(\vect{y}_n)})^2 + (1 - \rho_{\mathrm{OLS}})^2
\end{equation}
Because $\mu(\vect{z}_{n, \mathrm{OLS}}) = \mu(\vect{y}_n)$, the bias term vanishes. The relation $\sigma_\mathrm{OLS}/\sigma(\vect{y}_n) = \rho_\mathrm{OLS}$ taken from eq.~\eqref{eq:ols_correlation_deriv} turns the variability term into the same form as the correlation term, giving:
\begin{equation}
\label{eq:OSL_training_LKG_preliminary}
L_\mathrm{KG}(\vect{z}_{n, \mathrm{OLS}}, \vect{y}_n) = 0 + (1 - \rho_\mathrm{OLS})^2 + (1 - \rho_\mathrm{OLS})^2 = 2 (1 - \rho_\mathrm{OLS})^2
\end{equation}
To write the result with $\rho_{\mathrm{OLS}_{(0)}}$ while retaining the restriction $\vect{s}_p \neq \0_p$, we write eq.~\eqref{eq:OSL_training_LKG_preliminary} as:
\begin{equation}
\label{eq:OSL_training_LKG}
L_\mathrm{KG}(\vect{z}_{n, \mathrm{OLS}}, \vect{y}_n) = 2 (1 - \rho_{\mathrm{OLS}_{(0)}})^2, \vect{s}_p \neq \0_p
\end{equation}
The $\mathrm{KGE}$ follows directly from the relation between $L_\mathrm{KG}$ and $\mathrm{KGE}$ (eq.~\eqref{eq:KGE_LKG_relation}) together with eq.~\eqref{eq:OSL_training_LKG}:
\begin{equation}
\mathrm{KGE}(\vect{z}_{n, \mathrm{OLS}}, \vect{y}_n) = 1 - \sqrt{2 (1 - \rho_{\mathrm{OLS}_{(0)}})^2} = \sqrt{2} \rho_{\mathrm{OLS}_{(0)}} + 1 - \sqrt{2}, \vect{s}_p \neq \0_p
\end{equation}
\noindent\makebox[1cm][l]{(ii)} \textbf{Kling-Gupta predictions}: For the Kling-Gupta linear regression model, the predictions on the training set satisfy $\mu(\vect{z}_{n, \mathrm{KG}}) = \mu(\vect{y}_n)$ (eq.~\eqref{eq:KG_mu_predictions}), $\sigma(\vect{z}_{n, \mathrm{KG}}) = \sigma(\vect{y}_n)$ (eq.~\eqref{eq:KG_sigma_predictions}), and $\rho(\vect{z}_{n, \mathrm{KG}}, \vect{y}_n) = \rho_{\mathrm{OLS}_{(0)}}$ (eq.~\eqref{eq:KG_rho_predictions}). The $\mathrm{MSE}$ then follows from eq.~\eqref{eq:MSE_general_2} together with these properties:
\begin{equation}
\label{eq:KG_training_MSE}
\mathrm{MSE}(\vect{z}_{n, \mathrm{KG}}, \vect{y}_n) = \sigma^2(\vect{y}_n) + \sigma^2(\vect{y}_n) - 2 \rho_{\mathrm{OLS}_{(0)}} \sigma^2(\vect{y}_n) = 2 \sigma^2(\vect{y}_n) (1 - \rho_{\mathrm{OLS}_{(0)}})
\end{equation}
Therefore, from eqs.~\eqref{eq:Nash_Sutcliffe_Efficiency}, \eqref{eq:MSE_constant_x_n} and \eqref{eq:KG_training_MSE}, the $\mathrm{NSE}$ is
\begin{equation}
\mathrm{NSE}(\vect{z}_{n, \mathrm{KG}}, \vect{y}_n) = 1 - \frac{2 \sigma^2(\vect{y}_n) (1 - \rho_{\mathrm{OLS}_{(0)}})}{\sigma^2(\vect{y}_n)} = 2 \rho_{\mathrm{OLS}_{(0)}} - 1
\end{equation}
Substituting the statistical properties of the predictions, eq.~\eqref{eq:Kling_Gupta_loss} simplifies to:
\begin{equation}
\label{eq:KG_training_LKG}
L_\mathrm{KG}(\vect{z}_{n, \mathrm{KG}}, \vect{y}_n) = (1 - 1)^2 + (1 - 1)^2 + (1 - \rho_{\mathrm{OLS}_{(0)}})^2 = (1 - \rho_{\mathrm{OLS}_{(0)}})^2
\end{equation}
From eqs.~\eqref{eq:KGE_LKG_relation} and \eqref{eq:KG_training_LKG}, the $\mathrm{KGE}$ follows immediately:
\begin{equation}
\mathrm{KGE}(\vect{z}_{n, \mathrm{KG}}, \vect{y}_n) = 1 - \sqrt{(1 - \rho_{\mathrm{OLS}_{(0)}})^2} = \rho_{\mathrm{OLS}_{(0)}}
\end{equation}
\noindent\makebox[1cm][l]{(iii)} \textbf{Single-predictor case}: For the single-predictor special case ($p = 1$), the same steps apply. It suffices to substitute $\rho(\vect{z}_{n, \mathrm{KG}}, \vect{y}_n) = |\rho(\vect{x}_n, \vect{y}_n)|$ (see eq.~\eqref{eq:rho_single_ols_predictions}) and to replace $\vect{s}_p = \0_p$ by $\rho(\vect{x}_n, \vect{y}_n) = 0$.

Thus, all entries in Table~\ref{tab:training-performance} are now verified. \qed
\prooftitle{Proof B.4}{(Asymptotic performance metrics on the test set, Section~\ref{sec:comparative_infinite_training_test_data})}{proof:b4}
Consider the asymptotic setting where the training sample size $n \to \infty$. In this limit, the OLS and Kling-Gupta parameter estimates for the linear model in eq.~\eqref{eq:linear_model} converge almost surely to the population limits shown in Section~\ref{sec:asymptotic_properties}. We now evaluate the performance of these limiting models on an independent test set with $k$ i.i.d.\ observations $(\mat{X}_{k \times p, \mathrm{test}}, \vect{y}_{k, \mathrm{test}})$ generated from the same joint distribution $F_{\underline{\vect{x}}_p, \underline{y}}$. The predictions for the two estimated models are (eqs.~\eqref{eq:k_OLS_predictions} and \eqref{eq:k_KG_predictions}):
\begin{equation}
\label{eq:k_OLS_predictions}
\vect{z}_{k, \mathrm{OLS}, \mathrm{test}} = \mat{X}_{k \times p, \mathrm{test}} \vect{a}_{p, \mathrm{OLS}^*} + b_{p, \mathrm{OLS}^*} \1_k
\end{equation}
\begin{equation}
\label{eq:k_KG_predictions}
\vect{z}_{k, \mathrm{KG}, \mathrm{test}} = \mat{X}_{k \times p, \mathrm{test}} \vect{a}_{p, \mathrm{KG}^*} + b_{p, \mathrm{KG}^*} \1_k
\end{equation}
In the analysis that follows, limits are taken as the test set size $k \to \infty$. Using the strong law of large numbers and the continuous mapping theorem, we determine the limiting values of the performance metrics. Specifically, we prove that as $k \to \infty$, each metric converges to the asymptotic training set values reported in Table~\ref{tab:infinite-training-performance}.

\noindent\makebox[1cm][l]{(i)} \textbf{Mean Squared Error}

\noindent\makebox[1cm][l]{(i.1)} \textit{OLS predictions}: The test set $\mathrm{MSE}$ of the OLS predictions, from eq.~\eqref{eq:MSE_general_1}, is:
\begin{equation}
\label{eq:MSE_OLS_finite}
\underline{\mathrm{MSE}}_\mathrm{OLS}^{(k)} = \mathrm{MSE}(\underline{\vect{z}}_{k, \mathrm{OLS}, \mathrm{test}}, \underline{\vect{y}}_{k, \mathrm{test}}) = \frac{1}{k} \left\|\underline{\vect{y}}_{k, \mathrm{test}} - \underline{\mat{X}}_{k \times p, \mathrm{test}} \vect{a}_{p, \mathrm{OLS}^*} - b_{p, \mathrm{OLS}^*} \1_k \right\|_2^2
\end{equation}
Expanding the squared Euclidean norm in eq.~\eqref{eq:MSE_OLS_finite} gives:
\begin{equation}
\label{eq:MSE_OLS_expand}
\begin{aligned}
\underline{\mathrm{MSE}}_{\mathrm{OLS}}^{(k)} &= \frac{1}{k} (\underline{\vect{y}}_{k, \mathrm{test}}^\mathsf{T} \underline{\vect{y}}_{k, \mathrm{test}} - 2 \underline{\vect{y}}_{k, \mathrm{test}}^\mathsf{T} (\underline{\mat{X}}_{k \times p, \mathrm{test}} \vect{a}_{p, \mathrm{OLS}^*} + b_{p, \mathrm{OLS}^*} \1_k) \\
&\qquad\quad + (\underline{\mat{X}}_{k \times p, \mathrm{test}} \vect{a}_{p, \mathrm{OLS}^*} + b_{p, \mathrm{OLS}^*} \1_k)^\mathsf{T} (\underline{\mat{X}}_{k \times p, \mathrm{test}} \vect{a}_{p, \mathrm{OLS}^*} + b_{p, \mathrm{OLS}^*} \1_k)) \\[2pt]
&= \frac{\underline{\vect{y}}_{k, \mathrm{test}}^\mathsf{T} \underline{\vect{y}}_{k, \mathrm{test}}}{k} - 2 \frac{\underline{\vect{y}}_{k, \mathrm{test}}^\mathsf{T} \underline{\mat{X}}_{k \times p, \mathrm{test}}}{k} \vect{a}_{p, \mathrm{OLS}^*} - 2 b_{p, \mathrm{OLS}^*}\frac{\underline{\vect{y}}_{k, \mathrm{test}}^\mathsf{T} \1_k}{k} \\
&\quad + \vect{a}_{p, \mathrm{OLS}^*}^\mathsf{T} \frac{\underline{\mat{X}}_{k \times p, \mathrm{test}}^\mathsf{T} \underline{\mat{X}}_{k \times p, \mathrm{test}}}{k} \vect{a}_{p, \mathrm{OLS}^*} + 2 b_{p, \mathrm{OLS}^*} \frac{\1_k^\mathsf{T} \underline{\mat{X}}_{k \times p, \mathrm{test}}}{k} \vect{a}_{p, \mathrm{OLS}^*} + b_{p, \mathrm{OLS}^*}^2
\end{aligned}
\end{equation}
By the continuous mapping theorem, each term converges almost surely as $k \to \infty$ to the population counterpart in Appendix~\ref{app:notation}. Thus:
\begin{equation}
\label{eq:lim_y2}
\frac{\underline{\vect{y}}_{k, \mathrm{test}}^\mathsf{T} \underline{\vect{y}}_{k, \mathrm{test}}}{k} \overset{\mathrm{a.s.}}{\longrightarrow} \mathbb{E}_{F_{\underline{y}}}[\underline{y}^2] = \mathrm{Var}_{F_{\underline{y}}}[\underline{y}] + (\mathbb{E}_{F_{\underline{y}}}[\underline{y}])^2
\end{equation}
\begin{equation}
\label{eq:lim_yX}
\frac{\underline{\vect{y}}_{k, \mathrm{test}}^\mathsf{T} \underline{\mat{X}}_{k \times p, \mathrm{test}}}{k} \overset{\mathrm{a.s.}}{\longrightarrow} \mathbb{E}_{F_{\underline{\vect{x}}_p, \underline{y}}}[\underline{y}\, \underline{\vect{x}}_p^\mathsf{T}] = (\mathrm{Cov}_{F_{\underline{\vect{x}}_p, \underline{y}}}(\underline{\vect{x}}_p, \underline{y}))^\mathsf{T} + \mathbb{E}_{F_{\underline{y}}}[\underline{y}] (\mathbb{E}_{F_{\underline{\vect{x}}_p}}[\underline{\vect{x}}_p])^\mathsf{T}
\end{equation}
\begin{equation}
\label{eq:lim_XTX}
\frac{\underline{\mat{X}}_{k \times p, \mathrm{test}}^\mathsf{T} \underline{\mat{X}}_{k \times p, \mathrm{test}}}{k} \overset{\mathrm{a.s.}}{\longrightarrow} \mathbb{E}_{F_{\underline{\vect{x}}_p}}[\underline{\vect{x}}_p \underline{\vect{x}}_p^\mathsf{T}] = \mathrm{Var}_{F_{\underline{\vect{x}}_p}}[\underline{\vect{x}}_p] + \mathbb{E}_{F_{\underline{\vect{x}}_p}}[\underline{\vect{x}}_p] (\mathbb{E}_{F_{\underline{\vect{x}}_p}}[\underline{\vect{x}}_p])^\mathsf{T}
\end{equation}
\begin{equation}
\label{eq:lim_y1}
\frac{\underline{\vect{y}}_{k, \mathrm{test}}^\mathsf{T} \1_k}{k} \overset{\mathrm{a.s.}}{\longrightarrow} \mathbb{E}_{F_{\underline{y}}}[\underline{y}]
\end{equation}
\begin{equation}
\label{eq:lim_means}
\frac{\1_k^\mathsf{T} \underline{\mat{X}}_{k \times p, \mathrm{test}}}{k} \overset{\mathrm{a.s.}}{\longrightarrow} (\mathbb{E}_{F_{\underline{\vect{x}}_p}}[\underline{\vect{x}}_p])^\mathsf{T}
\end{equation}
With eqs.~\eqref{eq:lim_y2}, \eqref{eq:lim_yX}, \eqref{eq:lim_XTX}, \eqref{eq:lim_y1}, and \eqref{eq:lim_means}, eq.~\eqref{eq:MSE_OLS_expand} becomes:
\begin{equation}
\label{eq:MSE_OLS_lim_with_means}
\begin{aligned}
\underline{\mathrm{MSE}}_{\mathrm{OLS}}^{(k)} \overset{\mathrm{a.s.}}{\longrightarrow} & (\mathrm{Var}_{F_{\underline{y}}}[\underline{y}] + (\mathbb{E}_{F_{\underline{y}}}[\underline{y}])^2) - 2 ((\mathrm{Cov}_{F_{\underline{\vect{x}}_p, \underline{y}}}(\underline{\vect{x}}_p, \underline{y}))^\mathsf{T} + \mathbb{E}_{F_{\underline{y}}}[\underline{y}] (\mathbb{E}_{F_{\underline{\vect{x}}_p}}[\underline{\vect{x}}_p])^\mathsf{T}) \vect{a}_{p, \mathrm{OLS}^*} \\
& - 2 b_{p, \mathrm{OLS}^*} \mathbb{E}_{F_{\underline{y}}}[\underline{y}] + \vect{a}_{p, \mathrm{OLS}^*}^\mathsf{T} (\mathrm{Var}_{F_{\underline{\vect{x}}_p}}[\underline{\vect{x}}_p] + \mathbb{E}_{F_{\underline{\vect{x}}_p}}[\underline{\vect{x}}_p] (\mathbb{E}_{F_{\underline{\vect{x}}_p}}[\underline{\vect{x}}_p])^\mathsf{T}) \vect{a}_{p, \mathrm{OLS}^*}\\
& + 2 b_{p, \mathrm{OLS}^*} (\mathbb{E}_{F_{\underline{\vect{x}}_p}}[\underline{\vect{x}}_p])^\mathsf{T} \vect{a}_{p, \mathrm{OLS}^*} + b_{p, \mathrm{OLS}^*}^2
\end{aligned}
\end{equation}
The asymptotic OLS intercept estimate is $b_{p, \mathrm{OLS}^*} = \mathbb{E}_{F_{\underline{y}}}[\underline{y}] - (\mathbb{E}_{F_{\underline{\vect{x}}_p}}[\underline{\vect{x}}_p])^\mathsf{T} \vect{a}_{p, \mathrm{OLS}^*}$ (eq.~\eqref{eq:OLS_asymptotic_intercept_estimate}). With this expression, the terms with $b_{p, \mathrm{OLS}^*}$ expand to:
\begin{equation}
\label{eq:b2_expand}
b_{p, \mathrm{OLS}^*}^2 = (\mathbb{E}_{F_{\underline{y}}}[\underline{y}])^2 - 2 \mathbb{E}_{F_{\underline{y}}}[\underline{y}] (\mathbb{E}_{F_{\underline{\vect{x}}_p}}[\underline{\vect{x}}_p])^\mathsf{T} \vect{a}_{p, \mathrm{OLS}^*} + ((\mathbb{E}_{F_{\underline{\vect{x}}_p}}[\underline{\vect{x}}_p])^\mathsf{T} \vect{a}_{p, \mathrm{OLS}^*})^2
\end{equation}
\begin{equation}
\label{eq:2b_mu_a}
2 b_{p, \mathrm{OLS}^*} (\mathbb{E}_{F_{\underline{\vect{x}}_p}}[\underline{\vect{x}}_p])^\mathsf{T} \vect{a}_{p, \mathrm{OLS}^*} = 2 \mathbb{E}_{F_{\underline{y}}}[\underline{y}] (\mathbb{E}_{F_{\underline{\vect{x}}_p}}[\underline{\vect{x}}_p])^\mathsf{T} \vect{a}_{p, \mathrm{OLS}^*} - 2 ((\mathbb{E}_{F_{\underline{\vect{x}}_p}}[\underline{\vect{x}}_p])^\mathsf{T} \vect{a}_{p, \mathrm{OLS}^*})^2
\end{equation}
\begin{equation}
\label{eq:m2b_muy}
-2 b_{p, \mathrm{OLS}^*} \mathbb{E}_{F_{\underline{y}}}[\underline{y}] = -2 (\mathbb{E}_{F_{\underline{y}}}[\underline{y}])^2 + 2 \mathbb{E}_{F_{\underline{y}}}[\underline{y}] (\mathbb{E}_{F_{\underline{\vect{x}}_p}}[\underline{\vect{x}}_p])^\mathsf{T} \vect{a}_{p, \mathrm{OLS}^*}
\end{equation}
In eq.~\eqref{eq:MSE_OLS_lim_with_means}, the quantities from eqs.~\eqref{eq:b2_expand}, \eqref{eq:2b_mu_a}, and \eqref{eq:m2b_muy} with $(\mathbb{E}_{F_{\underline{y}}}[\underline{y}])^2$, $\mathbb{E}_{F_{\underline{y}}}[\underline{y}] (\mathbb{E}_{F_{\underline{\vect{x}}_p}}[\underline{\vect{x}}_p])^\mathsf{T} \vect{a}_{p, \mathrm{OLS}^*}$, and $((\mathbb{E}_{F_{\underline{\vect{x}}_p}}[\underline{\vect{x}}_p])^\mathsf{T} \vect{a}_{p, \mathrm{OLS}^*})^2$ sum to zero. The remaining variance and covariance terms give:
\begin{equation}
\label{eq:MSE_OLS_lim_raw}
\underline{\mathrm{MSE}}_{\mathrm{OLS}}^{(k)} \overset{\mathrm{a.s.}}{\longrightarrow} \mathrm{Var}_{F_{\underline{y}}}[\underline{y}] - 2 \vect{a}_{p, \mathrm{OLS}^*}^\mathsf{T} \mathrm{Cov}_{F_{\underline{\vect{x}}_p, \underline{y}}}(\underline{\vect{x}}_p, \underline{y}) + \vect{a}_{p, \mathrm{OLS}^*}^\mathsf{T} \mathrm{Var}_{F_{\underline{\vect{x}}_p}}[\underline{\vect{x}}_p] \vect{a}_{p, \mathrm{OLS}^*}
\end{equation}
When $\mathrm{Cov}_{F_{\underline{\vect{x}}_p, \underline{y}}}(\underline{\vect{x}}_p, \underline{y}) \neq \0_p$, with the limits for $\vect{a}_{p, \mathrm{OLS}^*}$ and $\sigma_{\mathrm{OLS}^*}^2$ from eqs.~\eqref{eq:OLS_asymptotic_coefficient_estimate} and \eqref{eq:sigma_OLS_asymptotic}, and eq.~\eqref{eq:rho_OLS_asymptotic}, eq.~\eqref{eq:MSE_OLS_lim_raw} becomes:
\begin{equation}
\label{eq:MSE_OLS_lim_preliminary}
\underline{\mathrm{MSE}}_{\mathrm{OLS}}^{(k)} \overset{\mathrm{a.s.}}{\longrightarrow} \mathrm{Var}_{F_{\underline{y}}}[\underline{y}] - \sigma_{\mathrm{OLS}^*}^2 = \mathrm{Var}_{F_{\underline{y}}}[\underline{y}](1 - \rho_{\mathrm{OLS}^*}^2)
\end{equation}
When $\mathrm{Cov}_{F_{\underline{\vect{x}}_p, \underline{y}}}(\underline{\vect{x}}_p, \underline{y}) = \0_p$, $\vect{a}_{p, \mathrm{OLS}^*} = \0_p$ from eq.~\eqref{eq:OLS_asymptotic_coefficient_estimate}, and eq.~\eqref{eq:MSE_OLS_lim_raw} becomes:
\begin{equation}
\label{eq:MSE_OLS_lim_prel_0}
\underline{\mathrm{MSE}}_{\mathrm{OLS}}^{(k)} \overset{\mathrm{a.s.}}{\longrightarrow} \mathrm{Var}_{F_{\underline{y}}}[\underline{y}]
\end{equation}
Together, eqs.~\eqref{eq:MSE_OLS_lim_preliminary} and \eqref{eq:MSE_OLS_lim_prel_0} become:
\begin{equation}
\label{eq:MSE_OLS_lim}
\underline{\mathrm{MSE}}_{\mathrm{OLS}}^{(k)} \overset{\mathrm{a.s.}}{\longrightarrow} \mathrm{Var}_{F_{\underline{y}}}[\underline{y}] - \sigma_{\mathrm{OLS}^*}^2 = \mathrm{Var}_{F_{\underline{y}}}[\underline{y}](1 - \rho_{\mathrm{OLS}_{(0)}^*}^2)
\end{equation}
\noindent\makebox[1cm][l]{(i.2)} \textit{Kling-Gupta predictions}: The test set $\mathrm{MSE}$ of the Kling-Gupta predictions, from eq.~\eqref{eq:MSE_general_1}, is:
\begin{equation}
\label{eq:MSE_KG_finite}
\underline{\mathrm{MSE}}_{\mathrm{KG}}^{(k)} = \mathrm{MSE}(\underline{\vect{z}}_{k, \mathrm{KG}, \mathrm{test}}, \underline{\vect{y}}_{k, \mathrm{test}}) = \frac{1}{k} \left\|\underline{\vect{y}}_{k, \mathrm{test}} - \underline{\mat{X}}_{k \times p, \mathrm{test}} \vect{a}_{p, \mathrm{KG}^*} - b_{p, \mathrm{KG}^*} \1_k \right\|_2^2
\end{equation}
Expanding the squared Euclidean norm in eq.~\eqref{eq:MSE_KG_finite} gives:
\begin{equation}
\label{eq:MSE_KG_expand}
\begin{aligned}
\underline{\mathrm{MSE}}_{\mathrm{KG}}^{(k)} &= \frac{1}{k} (\underline{\vect{y}}_{k, \mathrm{test}}^\mathsf{T} \underline{\vect{y}}_{k, \mathrm{test}} - 2 \underline{\vect{y}}_{k, \mathrm{test}}^\mathsf{T} (\underline{\mat{X}}_{k \times p, \mathrm{test}} \vect{a}_{p, \mathrm{KG}^*} + b_{p, \mathrm{KG}^*} \1_k) \\
&\quad + (\underline{\mat{X}}_{k \times p, \mathrm{test}} \vect{a}_{p, \mathrm{KG}^*} + b_{p, \mathrm{KG}^*} \1_k)^\mathsf{T} (\underline{\mat{X}}_{k \times p, \mathrm{test}} \vect{a}_{p, \mathrm{KG}^*} + b_{p, \mathrm{KG}^*} \1_k)) \\[4pt]
&= \frac{\underline{\vect{y}}_{k, \mathrm{test}}^\mathsf{T} \underline{\vect{y}}_{k, \mathrm{test}}}{k} - 2 \frac{\underline{\vect{y}}_{k, \mathrm{test}}^\mathsf{T} \underline{\mat{X}}_{k \times p, \mathrm{test}}}{k} \vect{a}_{p, \mathrm{KG}^*} - 2 b_{p, \mathrm{KG}^*} \frac{\underline{\vect{y}}_{k, \mathrm{test}}^\mathsf{T} \1_k}{k} \\
&\quad + \vect{a}_{p, \mathrm{KG}^*}^\mathsf{T} \frac{\underline{\mat{X}}_{k \times p, \mathrm{test}}^\mathsf{T} \underline{\mat{X}}_{k \times p, \mathrm{test}}}{k} \vect{a}_{p, \mathrm{KG}^*} + 2 b_{p, \mathrm{KG}^*} \frac{\1_k^\mathsf{T} \underline{\mat{X}}_{k \times p, \mathrm{test}}}{k} \vect{a}_{p, \mathrm{KG}^*} + b_{p, \mathrm{KG}^*}^2
\end{aligned}
\end{equation}
With eqs.~\eqref{eq:lim_y2}, \eqref{eq:lim_yX}, \eqref{eq:lim_XTX}, \eqref{eq:lim_y1}, and \eqref{eq:lim_means}, eq.~\eqref{eq:MSE_KG_expand} becomes:
\begin{equation}
\label{eq:MSE_KG_lim_with_means}
\begin{aligned}
\underline{\mathrm{MSE}}_{\mathrm{KG}}^{(k)} \overset{\mathrm{a.s.}}{\longrightarrow}
& (\mathrm{Var}_{F_{\underline{y}}}[\underline{y}] + (\mathbb{E}_{F_{\underline{y}}}[\underline{y}])^2) - 2 ((\mathrm{Cov}_{F_{\underline{\vect{x}}_p, \underline{y}}}(\underline{\vect{x}}_p, \underline{y}))^\mathsf{T} + \mathbb{E}_{F_{\underline{y}}}[\underline{y}] (\mathbb{E}_{F_{\underline{\vect{x}}_p}}[\underline{\vect{x}}_p])^\mathsf{T}) \vect{a}_{p, \mathrm{KG}^*} \\
& - 2 b_{p, \mathrm{KG}^*} \mathbb{E}_{F_{\underline{y}}}[\underline{y}] + \vect{a}_{p, \mathrm{KG}^*}^\mathsf{T} (\mathrm{Var}_{F_{\underline{\vect{x}}_p}}[\underline{\vect{x}}_p] + \mathbb{E}_{F_{\underline{\vect{x}}_p}}[\underline{\vect{x}}_p] (\mathbb{E}_{F_{\underline{\vect{x}}_p}}[\underline{\vect{x}}_p])^\mathsf{T}) \vect{a}_{p, \mathrm{KG}^*} \\
& + 2 b_{p, \mathrm{KG}^*} (\mathbb{E}_{F_{\underline{\vect{x}}_p}}[\underline{\vect{x}}_p])^\mathsf{T} \vect{a}_{p, \mathrm{KG}^*} + b_{p, \mathrm{KG}^*}^2
\end{aligned}
\end{equation}
The asymptotic Kling-Gupta intercept is $b_{p, \mathrm{KG}^*} = \mathbb{E}_{F_{\underline{y}}}[\underline{y}] - (\mathbb{E}_{F_{\underline{\vect{x}}_p}}[\underline{\vect{x}}_p])^\mathsf{T} \vect{a}_{p, \mathrm{KG}^*}$ (eq.~\eqref{eq:KG_asymptotic_intercept_estimate}). With this expression, the terms with $b_{p, \mathrm{KG}^*}$ expand to:
\begin{equation}
\label{eq:b2_expand_KG}
b_{p, \mathrm{KG}^*}^2 = (\mathbb{E}_{F_{\underline{y}}}[\underline{y}])^2 - 2 \mathbb{E}_{F_{\underline{y}}}[\underline{y}] (\mathbb{E}_{F_{\underline{\vect{x}}_p}}[\underline{\vect{x}}_p])^\mathsf{T} \vect{a}_{p, \mathrm{KG}^*} + ((\mathbb{E}_{F_{\underline{\vect{x}}_p}}[\underline{\vect{x}}_p])^\mathsf{T} \vect{a}_{p, \mathrm{KG}^*})^2
\end{equation}
\begin{equation}
\label{eq:2b_mu_a_KG}
2 b_{p, \mathrm{KG}^*} (\mathbb{E}_{F_{\underline{\vect{x}}_p}}[\underline{\vect{x}}_p])^\mathsf{T} \vect{a}_{p, \mathrm{KG}^*} = 2 \mathbb{E}_{F_{\underline{y}}}[\underline{y}] (\mathbb{E}_{F_{\underline{\vect{x}}_p}}[\underline{\vect{x}}_p])^\mathsf{T} \vect{a}_{p, \mathrm{KG}^*} - 2 ((\mathbb{E}_{F_{\underline{\vect{x}}_p}}[\underline{\vect{x}}_p])^\mathsf{T} \vect{a}_{p, \mathrm{KG}^*})^2
\end{equation}
\begin{equation}
\label{eq:m2b_muy_KG}
-2 b_{p, \mathrm{KG}^*} \mathbb{E}_{F_{\underline{y}}}[\underline{y}] = -2 (\mathbb{E}_{F_{\underline{y}}}[\underline{y}])^2 + 2 \mathbb{E}_{F_{\underline{y}}}[\underline{y}] (\mathbb{E}_{F_{\underline{\vect{x}}_p}}[\underline{\vect{x}}_p])^\mathsf{T} \vect{a}_{p, \mathrm{KG}^*}
\end{equation}
In eq.~\eqref{eq:MSE_KG_lim_with_means}, the quantities from eqs.~\eqref{eq:b2_expand_KG}, \eqref{eq:2b_mu_a_KG}, and \eqref{eq:m2b_muy_KG} with $(\mathbb{E}_{F_{\underline{y}}}[\underline{y}])^2$, $\mathbb{E}_{F_{\underline{y}}}[\underline{y}] (\mathbb{E}_{F_{\underline{\vect{x}}_p}}[\underline{\vect{x}}_p])^\mathsf{T} \vect{a}_{p, \mathrm{KG}^*}$, and $((\mathbb{E}_{F_{\underline{\vect{x}}_p}}[\underline{\vect{x}}_p])^\mathsf{T} \vect{a}_{p, \mathrm{KG}^*})^2$ sum to zero. The remaining variance and covariance terms give the compact expression:
\begin{equation}
\label{eq:MSE_KG_lim_raw}
\underline{\mathrm{MSE}}_{\mathrm{KG}}^{(k)} \overset{\mathrm{a.s.}}{\longrightarrow} \mathrm{Var}_{F_{\underline{y}}}[\underline{y}] - 2 \vect{a}_{p, \mathrm{KG}^*}^\mathsf{T} \mathrm{Cov}_{F_{\underline{\vect{x}}_p, \underline{y}}}(\underline{\vect{x}}_p, \underline{y}) + \vect{a}_{p, \mathrm{KG}^*}^\mathsf{T} \mathrm{Var}_{F_{\underline{\vect{x}}_p}}[\underline{\vect{x}}_p] \vect{a}_{p, \mathrm{KG}^*}
\end{equation}
With the definition of $\vect{a}_{p, \mathrm{KG}^*}$ from Section~\ref{sec:asymptotic_properties}, and eqs.~\eqref{eq:KG_asymptotic_coefficient_estimate}, \eqref{eq:OLS_asymptotic_coefficient_estimate}, and \eqref{eq:sigma_OLS_asymptotic} for the non-zero cross-covariance case:
\begin{equation}
\label{eq:cross_KG_pop}
\vect{a}_{p, \mathrm{KG}^*}^\mathsf{T} \mathrm{Cov}_{F_{\underline{\vect{x}}_p, \underline{y}}}(\underline{\vect{x}}_p, \underline{y}) = \frac{\sqrt{\mathrm{Var}_{F_{\underline{y}}}[\underline{y}]}}{\sigma_{\mathrm{OLS}^*}} \sigma_{\mathrm{OLS}^*}^2 = \sqrt{\mathrm{Var}_{F_{\underline{y}}}[\underline{y}]} \sigma_{\mathrm{OLS}^*}
\end{equation}
In the zero cross-covariance case, $\mathrm{Cov}_{F_{\underline{\vect{x}}_p, \underline{y}}}(\underline{\vect{x}}_p, \underline{y}) = \0_p$, so:
\begin{equation}
\vect{a}_{p, \mathrm{KG}^*}^\mathsf{T} \mathrm{Cov}_{F_{\underline{\vect{x}}_p, \underline{y}}}(\underline{\vect{x}}_p, \underline{y}) = 0
\end{equation}
In the non-zero cross-covariance case, eqs.~\eqref{eq:KG_asymptotic_coefficient_estimate}, \eqref{eq:OLS_asymptotic_coefficient_estimate}, and \eqref{eq:sigma_OLS_asymptotic} give:
\begin{equation}
\label{eq:var_KG_pop}
\vect{a}_{p, \mathrm{KG}^*}^\mathsf{T} \mathrm{Var}_{F_{\underline{\vect{x}}_p}}[\underline{\vect{x}}_p] \vect{a}_{p, \mathrm{KG}^*} = \frac{\mathrm{Var}_{F_{\underline{y}}}[\underline{y}]}{\sigma_{\mathrm{OLS}^*}^2} \sigma_{\mathrm{OLS}^*}^2 = \mathrm{Var}_{F_{\underline{y}}}[\underline{y}]
\end{equation}
In the zero cross-covariance case, eq.~\eqref{eq:a_KG_asymptotic_ellipsoid} is the same equality:
\begin{equation}
\vect{a}_{p, \mathrm{KG}^*}^\mathsf{T} \mathrm{Var}_{F_{\underline{\vect{x}}_p}}[\underline{\vect{x}}_p] \vect{a}_{p, \mathrm{KG}^*} = \mathrm{Var}_{F_{\underline{y}}}[\underline{y}]
\end{equation}
With eqs.~\eqref{eq:cross_KG_pop} and \eqref{eq:var_KG_pop} for the non-zero cross-covariance case, eq.~\eqref{eq:MSE_KG_lim_raw} becomes:
\begin{equation}
\label{eq:MSE_KG_lim}
\underline{\mathrm{MSE}}_{\mathrm{KG}}^{(k)} \overset{\mathrm{a.s.}}{\longrightarrow} \mathrm{Var}_{F_{\underline{y}}}[\underline{y}] - 2 \sqrt{\mathrm{Var}_{F_{\underline{y}}}[\underline{y}]} \sigma_{\mathrm{OLS}^*} + \mathrm{Var}_{F_{\underline{y}}}[\underline{y}] = 2 \mathrm{Var}_{F_{\underline{y}}}[\underline{y}](1 - \rho_{\mathrm{OLS}_{(0)}^*})
\end{equation}
In the zero cross-covariance case, $\rho_{\mathrm{OLS}_{(0)}^*} = 0$, and eq.~\eqref{eq:MSE_KG_lim} is similarly:
\begin{equation}
\underline{\mathrm{MSE}}_{\mathrm{KG}}^{(k)} \overset{\mathrm{a.s.}}{\longrightarrow} \mathrm{Var}_{F_{\underline{y}}}[\underline{y}] - 0 + \mathrm{Var}_{F_{\underline{y}}}[\underline{y}] = 2 \mathrm{Var}_{F_{\underline{y}}}[\underline{y}] = 2 \mathrm{Var}_{F_{\underline{y}}}[\underline{y}] (1 - \rho_{\mathrm{OLS}_{(0)}^*})
\end{equation}
\noindent\makebox[1cm][l]{(ii)} \textbf{Kling-Gupta loss}

The Kling-Gupta loss in eq.~\eqref{eq:Kling_Gupta_loss} is formulated through sample statistics. In the limit as the test set size $k \to \infty$, Appendix~\ref{app:notation} demonstrates that the sample means, sample standard deviations, and the Pearson sample correlation between the predictions and the response variable converge almost surely to their population counterparts (functionals).

\noindent\makebox[1cm][l]{(ii.1)} \textit{OLS predictions}: Under the condition $\mathrm{Cov}_{F_{\underline{\vect{x}}_p, \underline{y}}}(\underline{\vect{x}}_p, \underline{y}) \neq \0_p$, the functionals of the OLS predictions $\underline{z}_{\mathrm{OLS}, \mathrm{test}} = \underline{\vect{x}}_{p, \mathrm{test}}^\mathsf{T} \vect{a}_{p, \mathrm{OLS}^*} + b_{p, \mathrm{OLS}^*}$ are:
\begin{equation}
\label{eq:mean_OLS_pop}
\mathbb{E}_{F_{\underline{\vect{x}}_p}}[\underline{z}_{\mathrm{OLS}, \mathrm{test}}] = \mathbb{E}_{F_{\underline{\vect{x}}_p}}[\underline{\vect{x}}_p]^\mathsf{T} \vect{a}_{p, \mathrm{OLS}^*} + b_{p, \mathrm{OLS}^*} = \mathbb{E}_{F_{\underline{y}}}[\underline{y}]
\end{equation}
\begin{equation}
\label{eq:var_OLS_pop}
\mathrm{Var}_{F_{\underline{\vect{x}}_p}}[\underline{z}_{\mathrm{OLS}, \mathrm{test}}] = \vect{a}_{p, \mathrm{OLS}^*}^\mathsf{T} \mathrm{Var}_{F_{\underline{\vect{x}}_p}}[\underline{\vect{x}}_p] \vect{a}_{p, \mathrm{OLS}^*} = \sigma_{\mathrm{OLS}^*}^2
\end{equation}
\begin{equation}
\label{eq:cov_OLS_pop}
\mathrm{Cov}_{F_{\underline{\vect{x}}_p, \underline{y}}}(\underline{z}_{\mathrm{OLS}, \mathrm{test}}, \underline{y}) = \vect{a}_{p, \mathrm{OLS}^*}^\mathsf{T} \mathrm{Cov}_{F_{\underline{\vect{x}}_p, \underline{y}}}(\underline{\vect{x}}_p, \underline{y}) = \sigma_{\mathrm{OLS}^*}^2
\end{equation}
Thus the population correlation (using eq.~\eqref{eq:rho_OLS_asymptotic}) is
\begin{equation}
\label{eq:corr_OLS_pop}
\mathrm{Corr}_{F_{\underline{\vect{x}}_p, \underline{y}}}(\underline{z}_{\mathrm{OLS}, \mathrm{test}}, \underline{y}) = \frac{\sigma_{\mathrm{OLS}^*}^2}{\sigma_{\mathrm{OLS}^*} \sqrt{\mathrm{Var}_{F_{\underline{y}}}[\underline{y}]}} = \rho_{\mathrm{OLS}^*}
\end{equation}
The three terms of the Kling-Gupta loss therefore converge as follows:
\begin{equation}
\label{eq:bias_OLS_pop}
(1 - \frac{\mu(\underline{\vect{z}}_{k, \mathrm{OLS}, \mathrm{test}})}{\mu(\underline{\vect{y}}_{k, \mathrm{test}})})^2 \overset{\mathrm{a.s.}} {\longrightarrow} (1 - \frac{\mathbb{E}_{F_{\underline{\vect{x}}_p}}[\underline{z}_{\mathrm{OLS}, \mathrm{test}}]}{\mathbb{E}_{F_{\underline{y}}}[\underline{y}]})^2 = 0
\end{equation}
\begin{equation}
\label{eq:varterm_OLS_pop}
(1 - \frac{\sigma(\underline{\vect{z}}_{k, \mathrm{OLS}, \mathrm{test}})}{\sigma(\underline{\vect{y}}_{k, \mathrm{test}})})^2 \overset{\mathrm{a.s.}}{\longrightarrow} (1 - \frac{\sigma_{\mathrm{OLS}^*}}{\sqrt{\mathrm{Var}_{F_{\underline{y}}}[\underline{y}]}})^2 = (1 - \rho_{\mathrm{OLS}^*})^2
\end{equation}
\begin{equation}
\label{eq:corrterm_OLS_pop}
(1 - \rho(\underline{\vect{z}}_{k, \mathrm{OLS}, \mathrm{test}}, \underline{\vect{y}}_{k, \mathrm{test}}))^2 \overset{\mathrm{a.s.}}{\longrightarrow} (1 - \mathrm{Corr}_{F_{\underline{\vect{x}}_p, \underline{y}}}(\underline{z}_{\mathrm{OLS}, \mathrm{test}}, \underline{y}))^2 = (1 - \rho_{\mathrm{OLS}^*})^2
\end{equation}
Combining these components and applying the definition of $\rho_{\mathrm{OLS}_{(0)}^*}$ from eq.~\eqref{eq:rho_OLS_0_asymptotic} under the assumption $\mathrm{Cov}_{F_{\underline{\vect{x}}_p, \underline{y}}}(\underline{\vect{x}}_p, \underline{y}) \neq \0_p$ gives:
\begin{equation}
\label{eq:LKG_OLS_lim}
L_\mathrm{KG}(\underline{\vect{z}}_{k, \mathrm{OLS}, \mathrm{test}}, \underline{\vect{y}}_{k, \mathrm{test}}) \overset{\mathrm{a.s.}}{\longrightarrow} 2 (1 - \rho_{\mathrm{OLS}^*})^2 = 2 (1 - \rho_{\mathrm{OLS}_{(0)}^*})^2
\end{equation}

\noindent\makebox[1cm][l]{(ii.2)} \textit{Kling-Gupta predictions}:
Analogously, the functionals for the Kling-Gupta predictions $\underline{z}_{\mathrm{KG}, \mathrm{test}} = \underline{\vect{x}}_{p, \mathrm{test}}^\mathsf{T} \vect{a}_{p, \mathrm{KG}^*} + b_{p, \mathrm{KG}^*}$ are:
\begin{equation}
\label{eq:mean_KG_pop}
\mathbb{E}_{F_{\underline{\vect{x}}_p}}[\underline{z}_{\mathrm{KG, \mathrm{test}}}] = \mathbb{E}_{F_{\underline{\vect{x}}_p}}[\underline{\vect{x}}_p]^\mathsf{T} \vect{a}_{p, \mathrm{KG}^*} + b_{p, \mathrm{KG}^*} = \mathbb{E}_{F_{\underline{y}}}[\underline{y}]
\end{equation}
\begin{equation}
\label{eq:var_KG_pop2}
\mathrm{Var}_{F_{\underline{\vect{x}}_p}}[\underline{z}_{\mathrm{KG, \mathrm{test}}}] = \vect{a}_{p, \mathrm{KG}^*}^\mathsf{T} \mathrm{Var}_{F_{\underline{\vect{x}}_p}}[\underline{\vect{x}}_p] \vect{a}_{p, \mathrm{KG}^*} = \mathrm{Var}_{F_{\underline{y}}}[\underline{y}]
\end{equation}
In the non-zero cross-covariance case:
\begin{equation}
\label{eq:cov_KG_pop}
\mathrm{Cov}_{F_{\underline{\vect{x}}_p, \underline{y}}}(\underline{z}_{\mathrm{KG, \mathrm{test}}}, \underline{y}) = \vect{a}_{p, \mathrm{KG}^*}^\mathsf{T} \mathrm{Cov}_{F_{\underline{\vect{x}}_p, \underline{y}}}(\underline{\vect{x}}_p, \underline{y}) = \sqrt{\mathrm{Var}_{F_{\underline{y}}}[\underline{y}]} \sigma_{\mathrm{OLS}^*}
\end{equation}
It follows that:
\begin{equation}
\label{eq:corr_KG_pop}
\mathrm{Corr}_{F_{\underline{\vect{x}}_p, \underline{y}}}(\underline{z}_{\mathrm{KG, \mathrm{test}}}, \underline{y}) = \frac{\sqrt{\mathrm{Var}_{F_{\underline{y}}}[\underline{y}]} \sigma_{\mathrm{OLS}^*}}{\sqrt{\mathrm{Var}_{F_{\underline{y}}}[\underline{y}]}\sqrt{\mathrm{Var}_{F_{\underline{y}}}[\underline{y}]}} = \rho_{\mathrm{OLS}^*}
\end{equation}
In the zero cross-covariance case, where $\mathrm{Cov}_{F_{\underline{\vect{x}}_p, \underline{y}}}(\underline{\vect{x}}_p, \underline{y}) = \0_p$:
\begin{equation}
\mathrm{Cov}_{F_{\underline{\vect{x}}_p, \underline{y}}}(\underline{z}_{\mathrm{KG, \mathrm{test}}}, \underline{y}) = \vect{a}_{p, \mathrm{KG}^*}^\mathsf{T} \mathrm{Cov}_{F_{\underline{\vect{x}}_p, \underline{y}}}(\underline{\vect{x}}_p, \underline{y}) = 0
\end{equation}
and therefore:
\begin{equation}
\mathrm{Corr}_{F_{\underline{\vect{x}}_p, \underline{y}}}(\underline{z}_{\mathrm{KG, \mathrm{test}}}, \underline{y}) = 0
\end{equation}
The three component terms converge as eq.~\eqref{eq:corrterm_KG_pop} combines both the non-zero cross-covariance and the zero cross-covariance case, see eq.~\eqref{eq:rho_OLS_0_asymptotic}:
\begin{equation}
\label{eq:bias_KG_pop}
(1 - \frac{\mu(\underline{\vect{z}}_{k, \mathrm{KG}, \mathrm{test}})}{\mu(\underline{\vect{y}}_{k, \mathrm{test}})})^2 \overset{\mathrm{a.s.}}{\longrightarrow} 0
\end{equation}
\begin{equation}
\label{eq:varterm_KG_pop}
(1 - \frac{\sigma(\underline{\vect{z}}_{k, \mathrm{KG}, \mathrm{test}})}{\sigma(\underline{\vect{y}}_{k, \mathrm{test}})})^2 \overset{\mathrm{a.s.}}{\longrightarrow} 0
\end{equation}
\begin{equation}
\label{eq:corrterm_KG_pop}
(1 - \rho(\underline{\vect{z}}_{k, \mathrm{KG}, \mathrm{test}}, \underline{\vect{y}}_{k, \mathrm{test}}))^2 \overset{\mathrm{a.s.}}{\longrightarrow} (1 - \rho_{\mathrm{OLS}_{(0)}^*})^2
\end{equation}
Therefore
\begin{equation}
\label{eq:LKG_KG_lim}
L_\mathrm{KG}(\underline{\vect{z}}_{k, \mathrm{KG}, \mathrm{test}}, \underline{\vect{y}}_{k, \mathrm{test}}) \overset{\mathrm{a.s.}}{\longrightarrow} (1 - \rho_{\mathrm{OLS}_{(0)}^*})^2.
\end{equation}
\noindent\makebox[1cm][l]{(iii)} \textbf{Nash-Sutcliffe efficiency and Kling-Gupta efficiency}

The $\mathrm{NSE}$ limits follow directly from the $\mathrm{MSE}$ limits in eqs.~\eqref{eq:MSE_OLS_lim} and \eqref{eq:MSE_KG_lim}. Given that eq.~\eqref{eq:MSE_constant_x_n} implies $\mathrm{MSE}(\mu(\vect{y}_{k, \mathrm{test}}) \1_k, \vect{y}_{k, \mathrm{test}}) = \sigma^2(\vect{y}_{k, \mathrm{test}})$, which converges almost surely to $\mathrm{Var}_{F_{\underline{y}}}[\underline{y}]$ as $k \to \infty$, eq.~\eqref{eq:Nash_Sutcliffe_Efficiency} gives:
\begin{equation}
\label{eq:NSE_OLS_lim}
\mathrm{NSE}(\underline{\vect{z}}_{k, \mathrm{OLS}, \mathrm{test}}, \underline{\vect{y}}_{k, \mathrm{test}}) \overset{\mathrm{a.s.}}{\longrightarrow} 1 - \frac{\mathrm{Var}_{F_{\underline{y}}}[\underline{y}] (1 - \rho_{\mathrm{OLS}_{(0)}^*}^2)}{\mathrm{Var}_{F_{\underline{y}}}[\underline{y}]} = \rho_{\mathrm{OLS}_{(0)}^*}^2
\end{equation}
\begin{equation}
\label{eq:NSE_KG_lim}
\mathrm{NSE}(\underline{\vect{z}}_{k, \mathrm{KG}, \mathrm{test}}, \underline{\vect{y}}_{k, \mathrm{test}}) \overset{\mathrm{a.s.}}{\longrightarrow} 1 - \frac{2 \mathrm{Var}_{F_{\underline{y}}}[\underline{y}] (1 - \rho_{\mathrm{OLS}_{(0)}^*})}{\mathrm{Var}_{F_{\underline{y}}}[\underline{y}]} = 2 \rho_{\mathrm{OLS}_{(0)}^*} - 1
\end{equation}
By applying the relationship in eq.~\eqref{eq:KGE_LKG_relation} to the results in eqs.~\eqref{eq:LKG_OLS_lim} and \eqref{eq:LKG_KG_lim}, the continuous mapping theorem establishes the following $\mathrm{KGE}$ limits:
\begin{equation}
\label{eq:KGE_OLS_lim}
\mathrm{KGE}(\underline{\vect{z}}_{k,\mathrm{OLS}, \mathrm{test}}, \underline{\vect{y}}_{k, \mathrm{test}}) \overset{\mathrm{a.s.}}{\longrightarrow} 1 - \sqrt{2(1 - \rho_{\mathrm{OLS}_{(0)}^*})^2} = \sqrt{2} \rho_{\mathrm{OLS}_{(0)}^*} + 1 - \sqrt{2}
\end{equation}
\begin{equation}
\label{eq:KGE_KG_lim}
\mathrm{KGE}(\underline{\vect{z}}_{k, \mathrm{KG}, \mathrm{test}}, \underline{\vect{y}}_{k, \mathrm{test}}) \overset{\mathrm{a.s.}}{\longrightarrow} 1 - \sqrt{(1 - \rho_{\mathrm{OLS}_{(0)}^*})^2} = \rho_{\mathrm{OLS}_{(0)}^*}
\end{equation}
These limits align with the corresponding expressions established for the training set in Table~\ref{tab:infinite-training-performance}, assuming the training and test distributions are identical. \qed

\prooftitle{Proof B.5}{(Intercept parameter estimate for the multiple-predictor linear model with fixed regression coefficients estimated with the Kling-Gupta loss, Section~\ref{sec:fixed_regression_coefficient})}{proof:b5}
We determine the intercept estimate for the multiple-predictor linear model $z = \vect{a}_p^\mathsf{T} \vect{x}_p + b_p$ specified in eq.~\eqref{eq:linear_model} when it is trained by minimizing the Kling-Gupta loss function defined in eq.~\eqref{eq:Kling_Gupta_loss} and the coefficient vector $\vect{a}_p$ is taken as fixed. We prove that the estimate $\widehat{b}_p$ is given by eq.~\eqref{eq:intercept_estimate_fixed_a_KG}.

Substituting the sample mean of the predictions from eq.~\eqref{eq:mean_predictions_multi} into the Kling-Gupta loss gives:
\begin{equation}
\begin{aligned}
L_\mathrm{KG}(\vect{z}_n, \vect{y}_n)
&= (1 - \frac{(\vect{\mu}(\mat{X}_{n \times p}))^\mathsf{T} \vect{a}_p + b_p}{\mu(\vect{y}_n)})^2 \\
&+ (1 - \frac{\sigma(\vect{z}_n)}{\sigma(\vect{y}_n)})^2 + (1 - \rho(\vect{z}_n, \vect{y}_n))^2, \mu(\vect{y}_n), \sigma(\vect{y}_n), \sigma(\vect{z}_n) \in \R\backslash\{0\}
\end{aligned}
\end{equation}
The parameter \(b_p\) appears only in the bias term, which is a non-negative squared quantity, while $\sigma(\vect{z}_n)$ and $\rho(\vect{z}_n, \vect{y}_n)$ do not depend on $b_p$ (as established in eqs.~\eqref{eq:sigma_predictions} and \eqref{eq:correlation_predictions_multi}). Therefore, the loss function is minimized with respect to $b_p$ by forcing this term to zero. Solving for $b_p$ gives:
\begin{equation}
\widehat{b}_{p, \mathrm{KG}} = \mu(\vect{y}_n) - (\vect{\mu}(\mat{X}_{n \times p}))^\mathsf{T} \vect{a}_p
\end{equation}
This completes the proof. \qed

\prooftitle{Proof B.6}{(Slope parameter estimate for the single-predictor linear model with fixed intercept estimated with the Kling-Gupta loss, Section~\ref{sec:fixed_intercept_kg_linear_regression})}{proof:b6}
We determine the parameter estimate for the single-predictor linear model $z = ax + b$, as specified in eq.~\eqref{eq:single_pred_model}, when it is estimated by minimizing the Kling-Gupta loss function defined in eq.~\eqref{eq:Kling_Gupta_loss} while treating the intercept parameter $b$ as fixed a priori. We demonstrate that the resulting estimate $\widehat{a}_\mathrm{KG}$ is given by eq.~\eqref{eq:slope_estimate_fixed_b_KG}.

For brevity, we adopt the auxiliary definitions for $u$, $v$, $\rho$, and $w$ established in eq.~\eqref{eq:u_v_w_r_parameters}. Substituting these quantities and the prediction vector $\vect{z}_n = a \vect{x}_n + b \1_n$ from eq.~\eqref{eq:single_pred_model_vec} into the loss definition in eq.~\eqref{eq:Kling_Gupta_loss} gives:
\begin{equation}
\label{eq:KG_with_b_fixed}
L_\mathrm{KG}(a \vect{x}_n + b \1_n, \vect{y}_n) = (w - au)^2 + (1 - |a|v)^2 + (1 - \operatorname{sign}(a) \rho)^2, a \neq 0
\end{equation}
The loss function exhibits a discontinuity at $a = 0$, with one-sided limits given by:
\begin{equation}
\lim_{a \to 0^+} L_\mathrm{KG}(a \vect{x}_n + b \1_n, \vect{y}_n) = (1 - \rho)^2 + w^2 + 1
\end{equation}
\begin{equation}
\lim_{a \to 0^-} L_\mathrm{KG}(a \vect{x}_n + b \1_n, \vect{y}_n) = (1 + \rho)^2 + w^2 + 1
\end{equation}
Furthermore, the limits at infinity are:
\begin{equation}
\lim_{a \to \pm \infty} L_\mathrm{KG}( a \vect{x}_n + b \1_n, \vect{y}_n) = \infty
\end{equation}
By completing the square with respect to $a$, the loss function can be rearranged as:
\begin{equation}
\label{eq:L_ap_completed_squares}
\begin{aligned}
L_\mathrm{KG}(a \vect{x}_n + b \1_n, \vect{y}_n)
&= (u^2 + v^2) (a - \frac{wu + \operatorname{sign}(a) v}{u^2 + v^2})^2 \\
& + (1 - \operatorname{sign}(a) \rho)^2 + \frac{(u - \operatorname{sign}(a) wv)^2}{u^2 + v^2}, a \neq 0
\end{aligned}
\end{equation}
The following inequalities are straightforward to verify and they are essential for the subsequent analysis:
\begin{equation}
\label{eq:u_w_v_inequality_1}
\frac{(u - wv)^2}{u^2 + v^2} \leq w^2 + 1, \text{with equality if and only if } wu = - v
\end{equation}
\begin{equation}
\label{eq:u_w_v_inequality_2}
\frac{(u + wv)^2}{u^2 + v^2} \leq w^2 + 1, \text{with equality if and only if } wu = v
\end{equation}
Because the quadratic term involving $a$ in eq.~\eqref{eq:L_ap_completed_squares} is non-negative, the loss for a fixed sign of $a$ attains its minimum when that term vanishes and the resulting value of $a$ has the correct sign.

\noindent\makebox[1cm][l]{(i)} \textbf{Region} $a > 0$ ($\operatorname{sign}(a) = 1$): Eq.~\eqref{eq:L_ap_completed_squares} simplifies to:
\begin{equation}
L_+(a) = (u^2 + v^2) (a - \frac{wu + v}{u^2 + v^2})^2 + (1 - \rho)^2 + \frac{(u - wv)^2}{u^2 + v^2}
\end{equation}
We define the candidate slope:
\begin{equation}
\label{eq:a_positive_candidate_slope}
a_+ = \frac{wu + v}{u^2 + v^2}
\end{equation}
which makes the quadratic term vanish under the assumption $wu + v > 0$. The loss evaluated at $a_{+}$ is:
\begin{equation}
L_\mathrm{KG}(a_+ \vect{x}_n + b \1_n, \vect{y}_n) = (1 - \rho)^2 + \frac{(u - wv)^2}{u^2 + v^2}
\end{equation}

\noindent\makebox[1cm][l]{(ii)} \textbf{Region} $a < 0$ ($\operatorname{sign}(a) = - 1$): Eq.~\eqref{eq:L_ap_completed_squares} simplifies to:
\begin{equation}
L_-(a) = (u^2 + v^2) (a - \frac{wu - v}{u^2 + v^2})^2 + (1 + \rho)^2 + \frac{(u + wv)^2}{u^2 + v^2}
\end{equation}
We define the candidate slope:
\begin{equation}
a_- = \frac{wu - v}{u^2 + v^2}
\end{equation}
which makes the quadratic term vanish under the assumption $wu - v < 0$. The loss evaluated at $a_-$ is:
\begin{equation}
L_\mathrm{KG}(a_- \vect{x}_n + b \1_n, \vect{y}_n) = (1 + \rho)^2 + \frac{(u + wv)^2}{u^2 + v^2}
\end{equation}
We observe that the condition $wu + v < 0$ together with $wu > v$ is impossible, because $v > 0$. Therefore, we examine the following three mutually exclusive cases:

\noindent\makebox[1cm][l]{(ii.1)} \textbf{Case \#1} ($wu + v > 0$ and $wu \geq v$): which simplifies to $wu \geq v$. In the region $a > 0$, inequality \eqref{eq:u_w_v_inequality_1} together with $wu \neq - v$, implies:
\begin{equation}
L_\mathrm{KG}(a_+ \vect{x}_n + b \1_n, \vect{y}_n) < \lim_{a \to 0^+} L_\mathrm{KG}(a \vect{x}_n + b \1_n, \vect{y}_n)
\end{equation}
\begin{equation}
L_\mathrm{KG}(a_+ \vect{x}_n + b \1_n, \vect{y}_n) < L_\mathrm{KG}(a \vect{x}_n + b \1_n, \vect{y}_n), a > 0
\end{equation}
Hence $a_+$ is a global minimum on the interval $a > 0$.

On the interval $a < 0$, $L_\mathrm{KG}$ is a strictly decreasing function of $a$, thus its infimum is approached as $a \to 0^-$. We therefore evaluate the difference:
\begin{equation}
L_\mathrm{KG}(a_+ \vect{x}_n + b \1_n, \vect{y}_n) - \lim_{a \to 0^-} L_\mathrm{KG}(a \vect{x}_n + b \1_n, \vect{y}_n) = \frac{(u - wv)^2}{u^2 + v^2} - (w^2 + 1) - 4 \rho
\end{equation}
This difference is not guaranteed to be negative; therefore $a_+$ is not automatically a global minimum over the real line. A global minimum is attained at $a_+$, if and only if $\frac{(u - wv)^2}{u^2 + v^2} \leq 4 \rho + w^2 + 1$. Because inequality \eqref{eq:u_w_v_inequality_1} ensures $\frac{(u - wv)^2}{u^2 + v^2} \leq w^2 + 1$, the condition is automatically satisfied when $\rho \geq 0$. If $\rho < 0$, whether $a_+$ is the global minimizer depends on the specific values of $u$, $v$, $w$ and $\rho$.

\noindent\makebox[1cm][l]{(ii.2)} \textbf{Case \#2} ($wu + v > 0, wu < v$): We treat the regions $a > 0$ and $a < 0$ separately:

\textit{Subcase \#2a} ($a_+ = \frac{wu + v}{u^2 + v^2} > 0$): On the interval $a > 0$, the function has a global minimum at $a_+$ because:
\begin{equation}
L_\mathrm{KG}(a_+ \vect{x}_n + b \1_n, \vect{y}_n) < \lim_{a \to 0^+} L_\mathrm{KG}(a \vect{x}_n + b \1_n, \vect{y}_n)
\end{equation}
\begin{equation}
L_\mathrm{KG}(a_+ \vect{x}_n + b \1_n, \vect{y}_n) < L_\mathrm{KG}(a \vect{x}_n + b \1_n, \vect{y}_n),\ a > 0
\end{equation}

\textit{Subcase \#2b} ($a_- = \frac{wu - v}{u^2 + v^2} < 0$): On the interval $a < 0$, the function has a global minimum at $a_-$ because:
\begin{equation}
L_\mathrm{KG}(a_- \vect{x}_n + b \1_n, \vect{y}_n) < \lim_{a \to 0^-} L_\mathrm{KG}(a \vect{x}_n + b \1_n, \vect{y}_n)
\end{equation}
\begin{equation}
L_\mathrm{KG}(a_- \vect{x}_n + b \1_n, \vect{y}_n) < L_\mathrm{KG}(a \vect{x}_n + b \1_n, \vect{y}_n), a < 0
\end{equation}
To determine which of these candidates is the overall global minimum, we compute the difference between their loss values:
\begin{equation}
L_\mathrm{KG}(a_+ \vect{x}_n + b \1_n, \vect{y}_n) - L_\mathrm{KG}(a_{-} \vect{x}_n + b \1_n, \vect{y}_n) = - 4 (\frac{wuv}{u^2 + v^2} + \rho)
\end{equation}
Thus:

(a) If $\rho > - \frac{wuv}{u^2 + v^2}$, the global minimum occurs at $a_+$.

(b) If $\rho < - \frac{wuv}{u^2 + v^2}$, the global minimum occurs at $a_-$.

(c) If $\rho = - \frac{wuv}{u^2 + v^2}$, then there are two global minima at $a_+$ and $a_-$.

\noindent\makebox[1cm][l]{(ii.3)} \textbf{Case \#3} ($wu + v \leq 0$ and $wu < v$): which simplifies to $wu + v \leq 0$. On the region $a < 0$, inequality \eqref{eq:u_w_v_inequality_2} together with $wu \neq v$ implies:
\begin{equation}
L_\mathrm{KG}(a_- \vect{x}_n + b \1_n, \vect{y}_n) < \lim_{a \to 0^-} L_\mathrm{KG}(a \vect{x}_n + b \1_n, \vect{y}_n)
\end{equation}
\begin{equation}
L_\mathrm{KG}(a_- \vect{x}_n + b \1_n, \vect{y}_n) < L_\mathrm{KG}(a \vect{x}_n + b \1_n, \vect{y}_n), a < 0
\end{equation}
Hence $a_-$ is a global minimum on the interval $a < 0$.

On the interval $a > 0$, $L_\mathrm{KG}$ is a strictly increasing function of $a$, thus its infimum is approached as $a \to 0^+$. We therefore evaluate the difference:
\begin{equation}
L_\mathrm{KG}(a_- \vect{x}_n + b \1_n, \vect{y}_n) - \lim_{a \to 0^+} L_\mathrm{KG}(a \vect{x}_n + b \1_n, \vect{y}_n) = \frac{(u + wv)^2}{u^2 + v^2} - (w^2 + 1) + 4 \rho
\end{equation}
This difference is not guaranteed to be negative; therefore $a_-$ is not automatically a global minimum over the real line. A global minimum is attained at $a_-$, if and only if $\frac{(u + wv)^2}{u^2 + v^2} \leq w^2 + 1 - 4\rho$. Because inequality \eqref{eq:u_w_v_inequality_2} ensures $\frac{(u + wv)^2}{u^2 + v^2} \leq w^2 + 1$, the condition is automatically satisfied when $\rho \leq 0$. If $\rho > 0$, whether $a_-$ is the global minimizer depends on the specific values of $u$, $v$, $w$ and $\rho$.

The three cases analyzed above are summarized by eq.~\eqref{eq:slope_estimate_fixed_b_KG}. \qed

\section{Statistical software}
\label{app:software}

All computations were performed in the \texttt{R} programming language (version 4.6.0) using \texttt{RStudio} version 2026.05.0+218. We used the following \texttt{R} packages: \texttt{airGRdatasets} \citep{delaigue2025}, \texttt{data.table} \citep{barrett2026}, \texttt{devtools} \citep{wickham2026, xie2018, xie2020}, \texttt{knitr} \citep{xie2014, xie2015, xie2025}, \texttt{rmarkdown} \citep{allaire2026}, and \texttt{tidyverse} \citep{wickham2023, wickham2019}.

\textbf{Supplementary information:} The data and computations are available at \url{https://doi.org/10.5281/zenodo.20621675}.

\textbf{Conflicts of interest:} The authors declare no conflict of interest.

\textbf{Statement:} During the preparation of this work, the authors used \texttt{DeepSeek-V4} and \texttt{Gemini 3.1 Pro} for \LaTeX{} formatting, language polishing, and readability enhancement. After using these tools, the authors reviewed and edited the content as needed and accept full responsibility for the published article.

\titleformat{\section}
  {\bfseries\fontsize{13}{18}\selectfont}
  {}
  {0pt}
  {}
\titlespacing*{\section}{0pt}{6pt}{6pt}
\renewcommand{\refname}{References}
\begingroup
\fontsize{12}{12}\selectfont
\setlength{\bibsep}{0pt}
\setlength{\itemsep}{0pt}
\setlength{\parsep}{0pt}
\setlength{\parskip}{0pt}
\setlength{\bibhang}{1cm}

\endgroup


\begin{thebibliography}{99}
\bibitem[Allaire et~al.(2026)]{allaire2026} Allaire JJ, Xie Y, Dervieux C, McPherson J, Luraschi J, Ushey K, Atkins A, Wickham H, Cheng J, Chang W, Iannone R (2026) rmarkdown: Dynamic Documents for R. R package version 2.31. \url{https://doi.org/10.32614/CRAN.package.rmarkdown}.

\bibitem[Amemiya(1973)]{amemiya1973} Amemiya T (1973) Regression analysis when the dependent variable is truncated normal. Econometrica 41(6):\mbox{997--1016}. \url{https://doi.org/10.2307/1914031}.

\bibitem[Amemiya(1985)]{amemiya1985} Amemiya T (1985) Advanced Econometrics. Cambridge, MA: Harvard University Press. ISBN: 9780674251991.

\bibitem[Amrhein et~al.(1996)]{amrhein1996} Amrhein M, Srinivasan B, Bonvin D, Schumacher MM (1996) On the rank deficiency and rank augmentation of the spectral measurement matrix. Chemometrics and Intelligent Laboratory Systems 33(1):\mbox{17--33}. \url{https://doi.org/10.1016/0169-7439(95)00086-0}.

\bibitem[Banerjee et~al.(2005)]{banerjee2005} Banerjee A, Guo X, Wang H (2005) On the optimality of conditional expectation as a Bregman predictor. IEEE Transactions on Information Theory 51(7):\mbox{2664--2669}. \url{https://doi.org/10.1109/TIT.2005.850145}.

\bibitem[Barrett et~al.(2026)]{barrett2026} Barrett T, Dowle M, Srinivasan A, Gorecki J, Chirico M, Hocking T, Schwendinger B, Krylov I (2026) data.table: Extension of 'data.frame'. R package version 1.18.4. \url{https://doi.org/10.32614/CRAN.package.data.table}.

\bibitem[Bennett et~al.(2013)]{bennett2013} Bennett ND, Croke BFW, Guariso G, Guillaume JHA, Hamilton SH, Jakeman AJ, Marsili-Libelli S, Newham LTH, Norton JP, Perrin C, Pierce SA, Robson B, Seppelt R, Voinov AA, Fath BD, Andreassian V (2013) Characterising performance of environmental models. Environmental Modelling and Software 40:\mbox{1--20}. \url{https://doi.org/10.1016/j.envsoft.2012.09.011}.

\bibitem[Beven(2025)]{beven2025} Beven KJ (2025) A short history of philosophies of hydrological model evaluation and hypothesis testing. Wiley Interdisciplinary Reviews: Water 12(1):e1761. \url{https://doi.org/10.1002/wat2.1761}.

\bibitem[Biondi et~al.(2012)]{biondi2012} Biondi D, Freni G, Iacobellis V, Mascaro G, Montanari A (2012) Validation of hydrological models: Conceptual basis, methodological approaches and a proposal for a code of practice. Physics and Chemistry of the Earth 42--44:\mbox{70--76}. \url{https://doi.org/10.1016/j.pce.2011.07.037}.

\bibitem[Clark et~al.(2021)]{clark2021} Clark MP, Vogel RM, Lamontagne JR, Mizukami N, Knoben WJM, Tang G, Gharari S, Freer JE, Whitfield PH, Shook KR, Papalexiou SM (2021) The abuse of popular performance metrics in hydrologic modeling. Water Resources Research 57(9):e2020WR029001. \url{https://doi.org/10.1029/2020WR029001}.

\bibitem[Delaigue et~al.(2025)]{delaigue2025} Delaigue O, Brigode P, Thirel G (2025) airGRdatasets: Hydro-Meteorological Catchments Datasets for the 'airGR' Packages. R package version 0.2.3. \url{https://doi.org/10.32614/CRAN.package.airGRdatasets}.

\bibitem[Dimitriadis et~al.(2024)]{dimitriadis2024} Dimitriadis T, Fissler T, Ziegel J (2024) Characterizing $M$-estimators. Biometrika 111(1):\mbox{339--346}. \url{https://doi.org/10.1093/biomet/asad026}.

\bibitem[Gentle(2024)]{gentle2024} Gentle JE (2024) Matrix Algebra. Springer Cham. \url{https://doi.org/10.1007/978-3-031-42144-0}.

\bibitem[Gneiting(2011)]{gneiting2011} Gneiting T (2011) Making and evaluating point forecasts. Journal of the American Statistical Association 106(494):\mbox{746--762}. \url{https://doi.org/10.1198/jasa.2011.r10138}.

\bibitem[Gneiting and Resin(2023)]{gneiting2023} Gneiting T, Resin J (2023) Regression diagnostics meets forecast evaluation: Conditional calibration, reliability diagrams, and coefficient of determination. Electronic Journal of Statistics 17(2):\mbox{3226--3286}. \url{https://doi.org/10.1214/23-EJS2180}.

\bibitem[Gupta et~al.(2009)]{gupta2009} Gupta HV, Kling H, Yilmaz KK, Martinez GF (2009) Decomposition of the mean squared error and NSE performance criteria: Implications for improving hydrological modelling. Journal of Hydrology 377(1--2):\mbox{80--91}. \url{https://doi.org/10.1016/j.jhydrol.2009.08.003}.

\bibitem[Huber(1964)]{huber1964} Huber PJ (1964) Robust estimation of a location parameter. The Annals of Mathematical Statistics 35(1):\mbox{73--101}. \url{https://doi.org/10.1214/aoms/1177703732}.

\bibitem[Huber(1967)]{huber1967} Huber PJ (1967) The behavior of maximum likelihood estimates under nonstandard conditions. In: Le Cam LM, Neyman J (eds) Proceedings of the Fifth Berkeley Symposium on Mathematical Statistics and Probability. Berkeley: University of California Press, Berkeley, pp 221--233.

\bibitem[Jackson et~al.(2019)]{jackson2019} Jackson EK, Roberts W, Nelsen B, Williams GP, Nelson EJ, Ames DP (2019) Introductory overview: Error metrics for hydrologic modelling - A review of common practices and an open source library to facilitate use and adoption. Environmental Modelling and Software 119:\mbox{32--48}. \url{https://doi.org/10.1016/j.envsoft.2019.05.001}.

\bibitem[Klemeš(1986)]{klemes1986} Klemeš V (1986) Operational testing of hydrological simulation models. Hydrological Sciences Journal 31(1):\mbox{13-24}. \url{https://doi.org/10.1080/02626668609491024}.

\bibitem[Kling et~al.(2012)]{kling2012} Kling H, Fuchs M, Paulin M (2012) Runoff conditions in the upper Danube basin under an ensemble of climate change scenarios. Journal of Hydrology 424--425:\mbox{264--277}. \url{https://doi.org/10.1016/j.jhydrol.2012.01.011}.

\bibitem[Knoben et~al.(2019)]{knoben2019} Knoben WJM, Freer JE, Woods RA (2019) Technical note: Inherent benchmark or not? Comparing Nash-Sutcliffe and Kling-Gupta efficiency scores. Hydrology and Earth System Sciences 23(10):\mbox{4323--4331}. \url{https://doi.org/10.5194/hess-23-4323-2019}.

\bibitem[Krause et~al.(2005)]{krause2005} Krause P, Boyle DP, Bäse F (2005) Comparison of different efficiency criteria for hydrological model assessment. Advances in Geosciences 5:\mbox{89--97}. \url{https://doi.org/10.5194/adgeo-5-89-2005}.

\bibitem[Melsen et~al.(2025)]{melsen2025} Melsen LA, Puy A, Torfs PJJF, Saltelli A (2025) The rise of the Nash-Sutcliffe efficiency in hydrology. Hydrological Sciences Journal 70(8):\mbox{1248--1259}. \url{https://doi.org/10.1080/02626667.2025.2475105}.

\bibitem[Montanari and Koutsoyiannis(2012)]{montanari2012} Montanari A, Koutsoyiannis D (2012) A blueprint for process-based modeling of uncertain hydrological systems. Water Resources Research 48(9):W09555. \url{https://doi.org/10.1029/2011WR011412}.

\bibitem[Moriasi et~al.(2007)]{moriasi2007} Moriasi DN, Arnold JG, Van Liew MW, Bingner RL, Harmel RD, Veith TL (2007) Model evaluation guidelines for systematic quantification of accuracy in watershed simulations. Transactions of the ASABE 50(3):\mbox{885--900}. \url{https://doi.org/10.13031/2013.23153}.

\bibitem[Moriasi et~al.(2012)]{moriasi2012} Moriasi DN, Wilson BN, Douglas-Mankin KR, Arnold JG, Gowda PH (2012) Hydrologic and water quality models: Use, calibration, and validation. Transactions of the ASABE 55(4):\mbox{1241--1247}. \url{https://doi.org/10.13031/2013.42265}.

\bibitem[Moriasi et~al.(2015a)]{moriasi2015a} Moriasi DN, Gitau MW, Pai N, Daggupati P (2015a) Hydrologic and water quality models: Performance measures and evaluation criteria. Transactions of the ASABE 58(6):\mbox{1763--1785}. \url{https://doi.org/10.13031/trans.58.10715}.

\bibitem[Moriasi et~al.(2015b)]{moriasi2015b} Moriasi DN, Zeckoski RW, Arnold JG, Baffaut CB, Malone RW, Daggupati P, Guzman JA, Saraswat D, Yuan Y, Wilson BW, Shirmohammadi A, Douglas-Mankin KR (2015b) Hydrologic and water quality models: Key calibration and validation topics. Transactions of the ASABE 58(6):\mbox{1609--1618}. \url{https://doi.org/10.13031/trans.58.11075}.

\bibitem[Murphy(1988)]{murphy1988} Murphy AH (1988) Skill scores based on the mean square error and their relationships to the correlation coefficient. Monthly Weather Review 116(12):\mbox{2417--2424}. \url{https://doi.org/10.1175/1520-0493(1988)116<2417:SSBOTM>2.0.CO;2}.

\bibitem[Murphy and Daan(1985)]{murphy1985} Murphy AH, Daan H (1985) Forecast evaluation. In: Murphy AH, Katz RW (eds) Probability, Statistics and Decision Making in the Atmospheric Sciences. CRC Press, pp 379--437.

\bibitem[Nash and Sutcliffe(1970)]{nash1970} Nash JE, Sutcliffe JV (1970) River flow forecasting through conceptual models part I - A discussion of principles. Journal of Hydrology 10(3):\mbox{282--290}. \url{https://doi.org/10.1016/0022-1694(70)90255-6}.

\bibitem[Newey and McFadden(1994)]{newey1994} Newey WK, McFadden D (1994) Large sample estimation and hypothesis testing. In: Engle RF, McFadden D (eds) Handbook of Econometrics, vol. 4. Elsevier, Amsterdam, pp 2111--2245. \url{https://doi.org/10.1016/S1573-4412(05)80005-4}.

\bibitem[Patton(2011)]{patton2011} Patton AJ (2011) Volatility forecast comparison using imperfect volatility proxies. Journal of Econometrics 160(1):\mbox{246--256}. \url{https://doi.org/10.1016/j.jeconom.2010.03.034}.

\bibitem[Patton(2020)]{patton2020} Patton AJ (2020) Comparing possibly misspecified forecasts. Journal of Business and Economic Statistics 38(4):\mbox{796--809}. \url{https://doi.org/10.1080/07350015.2019.1585256}.

\bibitem[Reichelstein and Osband(1984)]{reichelstein1984} Reichelstein S, Osband K (1984) Incentives in government contracts. Journal of Econometrics 24(2):\mbox{257--270}. \url{https://doi.org/10.1016/0047-2727(84)90029-X}.

\bibitem[Saerens(2000)]{saerens2000} Saerens M (2000) Building cost functions minimizing to some summary statistics. IEEE Transactions on Neural Networks 11(6):\mbox{1263--1271}. \url{https://doi.org/10.1109/72.883416}.

\bibitem[Savage(1971)]{savage1971} Savage LJ (1971) Elicitation of personal probabilities and expectations. Journal of the American Statistical Association 66(336):\mbox{783--801}. \url{https://doi.org/10.1080/01621459.1971.10482346}.

\bibitem[Seber(2007)]{seber2007} Seber GAF (2007) A Matrix Handbook for Statisticians. John Wiley \& Sons, Inc. \url{https://doi.org/10.1002/9780470226797}.

\bibitem[Tyralis and Papacharalampous(2025)]{tyralis2025ioa} Tyralis H, Papacharalampous G (2025) Loss functions arising from the index of agreement. \url{https://doi.org/10.48550/arXiv.2510.14714}.

\bibitem[Tyralis and Papacharalampous(2026)]{tyralis2026nse} Tyralis H, Papacharalampous G (2026) Learning with the Nash-Sutcliffe loss. \url{https://doi.org/10.48550/arXiv.2603.00968}.

\bibitem[Vapnik(1998)]{vapnik1998} Vapnik VN (1998) Statistical Learning Theory. New York: John Wiley \& Sons. ISBN: 978-0-471-03003-4.

\bibitem[Vrugt(2024)]{vrugt2024} Vrugt JA (2024) Distribution-based model evaluation and diagnostics: Elicitability, propriety, and scoring rules for hydrograph functionals. Water Resources Research, 60(6):e2023WR036710. \url{https://doi.org/10.1029/2023WR036710}.

\bibitem[Vrugt and de Oliveira(2022)]{vrugt2022} Vrugt JA, de Oliveira DY (2022) Confidence intervals of the Kling-Gupta efficiency. Journal of Hydrology 612(Part A):\mbox{127968}. \url{https://doi.org/10.1016/j.jhydrol.2022.127968}.

\bibitem[Wickham(2023)]{wickham2023} Wickham H (2023) tidyverse: Easily Install and Load the 'Tidyverse'. R package version 2.0.0. \url{https://doi.org/10.32614/CRAN.package.tidyverse}.

\bibitem[Wickham et~al.(2019)]{wickham2019} Wickham H, Averick M, Bryan J, Chang W, McGowan LD, François R, Grolemund G, Hayes A, Henry L, Hester J, Kuhn M, Pedersen TL, Miller E, Bache SM, Müller K, Ooms J, Robinson D, Seidel DP, Spinu V, Takahashi K, Vaughan D, Wilke C, Woo K, Yutani H (2019) Welcome to the tidyverse. Journal of Open Source Software 4(43):\mbox{1686}. \url{https://doi.org/10.21105/joss.01686}.

\bibitem[Wickham et~al.(2026)]{wickham2026} Wickham H, Hester J, Chang W, Bryan J (2026) devtools: Tools to Make Developing R Packages Easier. R package version 2.5.2. \url{https://doi.org/10.32614/CRAN.package.devtools}.

\bibitem[Williams(2025)]{williams2025} Williams GP (2025) Friends don't let friends use Nash-Sutcliffe Efficiency (NSE) or KGE for hydrologic model accuracy evaluation: A rant with data and suggestions for better practice. Environmental Modelling and Software 194:\mbox{106665}. \url{https://doi.org/10.1016/j.envsoft.2025.106665}.

\bibitem[Willmott(1981)]{willmott1981} Willmott CJ (1981) On the validation of models. Physical Geography 2(2):\mbox{184--194}. \url{https://doi.org/10.1080/02723646.1981.10642213}.

\bibitem[Xie(2014)]{xie2014} Xie Y (2014) knitr: A Comprehensive Tool for Reproducible Research in R. In: Stodden V, Leisch F, Peng RD (eds) Implementing Reproducible Computational Research. New York: Chapman and Hall/CRC. \url{https://doi.org/10.1201/9781315373461-1}.

\bibitem[Xie(2015)]{xie2015} Xie Y (2015) Dynamic Documents with R and knitr. 2nd edition. New York: Chapman and Hall/CRC. \url{https://doi.org/10.1201/9781315382487}.

\bibitem[Xie(2025)]{xie2025} Xie Y (2025) knitr: A General-Purpose Package for Dynamic Report Generation in R. R package version 1.51. \url{https://doi.org/10.32614/CRAN.package.knitr}.

\bibitem[Xie et~al.(2018)]{xie2018} Xie Y, Allaire J, Grolemund G (2018) R Markdown: The Definitive Guide. Boca Raton, Florida: Chapman and Hall/CRC. ISBN: 9781138359338.

\bibitem[Xie et~al.(2020)]{xie2020} Xie Y, Dervieux C, Riederer E (2020) R Markdown Cookbook. Boca Raton, Florida: Chapman and Hall/CRC. ISBN: 9780367563837.
\end{thebibliography}
\end{document}